\documentclass[11pt]{amsart} 
\usepackage{graphicx}
\graphicspath{{images/}}
\usepackage{tikz-cd}
\usepackage{enumerate}
\usepackage{amsmath, amsfonts, amssymb, amscd, amstext, amsthm} 
\usepackage{latexsym}
\usepackage{mathrsfs}
\usepackage{mathscinet}

\usepackage[hidelinks, pagebackref]{hyperref}
\usepackage{url}
\usepackage[margin=3cm]{geometry}

\usepackage{marginnote}
\long\def\@savemarbox#1#2{\global\setbox#1\vtop{\hsize\marginparwidth 
  \@parboxrestore\tiny\raggedright #2}}
\renewcommand*{\backref}[1]{}
\renewcommand*{\backrefalt}[4]{
  \ifcase #1
  [No citations.]
  \or [#2]
  \else [#2]
  \fi }

\AtBeginDocument{%
   \def\MR#1{}
}

\makeatletter
\@namedef{subjclassname@2020}{\textup{2020} Mathematics Subject Classification}
\makeatother

\usepackage{subfiles}

\newcommand{\tts}{\hspace{.05555em}}
\newcommand{\R}{\mathbb{R}}
\newcommand{\Z}{\mathbb{Z}}
\newcommand{\N}{\mathbb{N}}
\newcommand{\V}{\mathbb{V}}
\newcommand{\U}{\mathbb{U}}
\newcommand{\T}{\mathbb{T}}
\newcommand{\F}{\mathcal{F}}
\newcommand{\Sp}{\mathbb{S}}
\newcommand{\abs}[1]{\left| #1 \right|}
\DeclareMathOperator{\stsys}{stabsys} 
\DeclareMathOperator{\ison}{IN} 
\DeclareMathOperator{\sobn}{SN}
\DeclareMathOperator{\pd}{PD}
\DeclareMathOperator{\diam}{diam}

\DeclareMathOperator{\vol}{\operatorname{\mathsf{Vol}}\tts}
\theoremstyle{plain}
\newtheorem{theorem}{Theorem}[section]
\newtheorem{lemma}[theorem]{Lemma}
\newtheorem{proposition}[theorem]{Proposition}
\newtheorem{corollary}[theorem]{Corollary}

\theoremstyle{definition}
\newtheorem{remark}[theorem]{Remark}
\newtheorem{definition}[theorem]{Definition}
\newtheorem{notation}[theorem]{Notation}
\title[Torus Stability]{Stability for a class of three-tori with small negative scalar curvature}

\author[E. Bryden]{Edward Bryden}
\curraddr{Universiteit Antwerpen}
\email{etbryden@gmail.com}
\thanks{Supported by the FWO (grant 12F0223N)}

\author[L. Chen]{Lizhi Chen}
\curraddr{Lanzhou University}
\email{chenzhmath@gmail.com}
\thanks{Partially supported by Youth Scientists Fund of NSFC (award No. 11901261) and NSFC 12271225}

\begin{document}
    \begin{abstract}
        We define a flexible class of Riemmanian metrics on the three-torus. Then, using Stern's
        inequality relating scalar curvature to harmonic one-forms, we show that any sequence of
        metrics in this family whose negative part of the scalar curvature tends to zero in
        $L^2$ norm has a subsequence which converges to some flat metric on the three-torus
        in the sense of Dong-Song.
    \end{abstract}
    \maketitle
	\section{Introduction}
	
	General relationships between geometric quantities are as beautiful as they are useful,
	and they are very useful.
	There are many famous examples of such relationships, all of which have far reaching 
	implications.
	Especially relevant for the current work are the volume growth and isoperimetric
	inequalities implied by lower Ricci curvature bounds,
	the systolic inequality of Gromov giving a lower bound on the volume of a Riemannian
	manifold in terms of its systole, and the relationship between the negative part 
	of the scalar curvature and the integral norm of the Hessian of harmonic maps into 
	$\Sp^{1}$ given by Stern in \cite{stern2019scalar}.

	It is natural to wonder what the extreme geometries are with respect to these relationships,
	and in what sense, if at all, they are unique.
	This is the question of rigidity.
	For example, we have the classical fact that any metric on the three-torus with non-negative
	scalar curvature must be flat, see for example \cite{gromov_lawson_spin_scalar_80} and
	\cite{schoen_yau_79}.
	In fact, for three-tori this follows from Stern's inequality in \cite{stern2019scalar}, mentioned above.

	Once the hard work of establishing such a rigidity result has been done, it is natural to
	wonder what can be said about those metrics which are nearly extremal, which is the
	question of \textit{stability}.
	In the context of this paper, the question is whether metrics on the three-torus whose
	scalar curvature has small negative part must be close to a flat metric in some sense.
	This is a subtle question; for a more in depth discussion of the ideas and difficulties
	involved one can read Sormani's survey article \cite{SormaniIAS}.

	So far there seems to be at least three geometric phenomena that complicate the study of
	three-tori with almost non-negative scalar curvature.
	The first two, \textit{other worlds} and \textit{splines}, have been expected to occur
	since the work of Gromov-Lawson \cite{gromov_lawson_spin_scalar_80} and Schoen-Yau
	\cite{schoen_yau_79}.
	Rigorous examples showing the existence and ubiquity of such objects have recently
	been constructed by Sweeney \cite{sweeney2023examples}.
	The third, \textit{drawstrings}, was first observed in dimensions greater than 3
	by Lee-Naber-Neumayer, see \cite{Lee-Naber-Neumayer-dp_and_scalar}.
	Later, Lee-Topping showed that drawstrings can be used to produce
	counter intuitive convergence results, see \cite{lee2022metric}.
    Drawstrings were shown to exist in three dimensional tori by Kazaras-Xu in 
    \cite{Kazaras-Xu-23}.

	Roughly speaking, the spaces and notions of convergence proposed to study stability
	problems involving scalar curvature lower bounds correspond to which of these three
	phenomena should be considered ``small perturbations", and which 
	are to be eliminated, or controlled, by hypothesis.
	Take for example the amazing result of Dong-Song on the stability of the
	Positive Mass Theorem \cite{Dong-Song-2023}.
	One way to interpret this result is to say that for metrics with nonnegative scalar
    curvature and small mass, such wild geometric phenomena as
	bubbles, splines, and drawstrings must be hidden behind a surface with small area.
	In this way, they are small perturbations of the geometry.

	At the opposite end of the spectrum, we may make an hypothesis which severely controls
	bubbles, splines, and drawstrings.
	For example, one may restrict attention to metrics satisfying a uniform lower bound
	on their Ricci curvature.
	Results of this flavor are the stability of the Positive Mass Theorem proven by 
	Kazaras-Khuri-Lee in \cite{Kazaras-Khuri-Lee}, and the stability for three-tori
	proven by Honda-Ketterer-Mondello-Perales-Rigoni in 
	\cite{Honda-Ketterer-Mondello-Perales-Rigoni-Ricci_lowerbound}.
	Results in a similar vein were obtained in \cite{Allen-Bryden-Kazaras_int_curv_Iso_stab}
	through controlling the geometry by assuming isoperimetric and integral Ricci curvature 
	bounds.

	Another approach has been to assume some geometric condition which eliminates the
	existence of one, or perhaps two, of the wild geometries.
	This approach has seen some success using the Intrinsic Flat distance on
	integral currents, which contain Riemannian manifolds as a subset.
	See for example the work of Allen, Kazaras, and the first named author in
	\cite{Allen-Bryden-Kazaras_Llarull-Stab-dim-3} the specific hypotheses assumed in this  work control bubbles, and eliminate drawstrings,
	but allow splines to exist.
    Using spinor techniques Hirsch and Zhang establish the stability of Llarull's theorem
    in great generality while only needing to control drawstrings, see
	\cite{Hirsch-Zhang-Llarull-stab}.
	Concerning the stability of tori, there is the work of 
	Allen-Vazquez Hernandez-Parise-Payne-Wang \cite{Allen-HernandezVazquez-Parise-Payne-Wang}
    which establishes a stability result for warped product metrics.
    Additionally, stability was established for metrics conformal to elements in a controlled family of flat metrics on the torus \cite{Allen}.

	In a different direction, Lee-Naber-Neumayer \cite{Lee-Naber-Neumayer-dp_and_scalar} give
	conditions under which splines and bubbles are eliminated, but drawstrings are 
	allowed to persist. They are then able to establish a scalar curvature
	stability result using the $d_{p}$ distance, also defined in
	\cite{Lee-Naber-Neumayer-dp_and_scalar}. In \cite{mazurowski2024stability} Mazurowski and Yao use the $d_{p}$
    distance to study the stability of the Yamabe invariant on $\Sp^3$.

	In the present work we will study a family of metrics which allow splines, and a tamer
	version of drawstrings to persist, but which eliminates bubbles. For this family of
	metrics we will prove that metrics with small negative scalar curvature are close to a
	flat metric in the Dong-Song sense \cite{Dong-PMT_Stability}, see Definition \ref{defn:Convergence_Dong-Song}.
    \begin{theorem}
		Fix $V,R,\Lambda,\eta,M>0$ and let $\F(V,R,\Lambda,\eta,M)$ be the family
		of Riemannian metrics on $\T^{3}$ such that
        \begin{enumerate}
				\item $|\T^{3}|_{g}\leq V$;
				\item $\|R^{-}_{g}\|_{L^{2}(g)}\leq R$;
				\item $\ison_{1}(g)\geq\Lambda$;
				\item $\min\bigl\{\stsys_{1}(g),\stsys_{2}(g)\bigr\}\geq\sigma$;
				\item $\kappa(g,\eta)\leq M$, see Definition \ref{defn:covering_constant}.
			\end{enumerate}
		Let $g_{i}$ be a sequence of metrics in $\F$ such that
			\begin{equation}
				\lim_{i\rightarrow{}\infty{}}\|R^{-}_{g_{i}}\|_{L^{2}(g_{i})}=0.
			\end{equation}
		Then, there is a subsequence, also denoted $g_{i}$, and a flat metric $g_{F_{\infty}}$ on
		$\T^{3}$ such that $g_{i}$ converges to $g_{F_{\infty{}}}$ in the sense of Dong-Song.
		That is, for any $\varepsilon{}>0$ there exists an $N\in{}\N$ such that
		for all $i\geq N$ there is an open submanifold $\widetilde{\Omega}_{i}$ with smooth
		boundary such that
			\begin{equation}
				|\widetilde{\Omega}_{i}^{c}|_{g_{i}}+|\partial{}
				\widetilde{\Omega}_{i}|_{g_{i}}\leq
				\varepsilon{},
			\end{equation}
		and
			\begin{equation}
				d_{GH}\left(\left(\widetilde{\Omega}_{i},
					\hat{d}^{g_{i}}_{\widetilde{\Omega}_{i}}\right),
					\left(\T^{3},d^{g_{F_{\infty{}}}}\right)
				\right)\leq\varepsilon{}.
			\end{equation}
	\end{theorem}
	As a consequence of the above, we obtain the following two theorems.
	\begin{theorem}
		Let $\sigma,K,D>0$, and define $\mathcal{R}(\sigma,K,D)$ to be the collection
		of Riemannian metrics $g$ on $\T^3$ such that 
			\begin{enumerate}
				\item $\min\{\mathrm{stsys}_1(g),\mathrm{stsys}_{2}(g)\}\geq\sigma$;
				\item $\mathrm{Ric}_{g}\geq-K$;
				\item $\diam_g(\T^3)\leq D$.
			\end{enumerate}
		Then, for any sequence of metrics $\{g_{i}\}_{i=1}^{\infty{}}\subset{}
		\mathcal{R}(\sigma,K,D)$ such that 
			\begin{equation}
				\lim_{i\rightarrow{}\infty{}}\|R^{-}_{g_{i}}\|_{L^{2}(g_{i})}=0,
			\end{equation}
		there exists a subsequence $\{g_{i_{j}}\}_{j=1}^{\infty{}}$ and a
		flat metric $g_{F_{\infty}}$ on $\T^{3}$ such that $g_{i_{j}}\rightarrow{}g_{F_{\infty}}$ in the
		sense of Dong-Song.
	\end{theorem}
	\begin{theorem}
		Let $g_0$ be a fixed Riemannian metric on $\T^3$, and let $\Lambda,V>0$.
		We denote by $\mathcal{V}(g_{0},\Lambda,V)$ the collection of Riemannian metrics
		$g$ on $\T^{3}$ satisfying the following properties:
			\begin{enumerate}
				\item $g\geq g_{0}$;
				\item $\ison_{1}(g)\geq\Lambda$;
				\item $|\T^{3}|_{g}\leq V$.
			\end{enumerate}
		Then, for any sequence of metrics $\{g_{i}\}_{i=1}^{\infty{}}\subset{}
		\mathcal{V}(g_0,\Lambda,V)$ such that 
			\begin{equation}
				\lim_{i\rightarrow{}\infty{}}\|R^{-}_{g_{i}}\|_{L^{2}(g_{i})}=0,
			\end{equation}
		there exists a subsequence $\{g_{i_{j}}\}_{j=1}^{\infty{}}$ and a
		flat metric $g_{F_{\infty}}$ on $\T^{3}$ such that $g_{i_{j}}\rightarrow{}g_{F_{\infty}}$ in the
		sense of Dong-Song.
	\end{theorem}
 The theorem above has the following corollary for Volume Above Distance Below (VADB) convergence,
 a notion of convergence which implies volume preserving intrinsic flat convergence, see \cite{Allen-Perales-Sormani-VADB-I}.
 \begin{corollary}
     Let $g_0$ be a fixed Riemannian metric on $\T^3$, let $\Lambda, R, V>0$,
     and let $\mathcal{V}=\mathcal{V}(g_0,\Lambda,R,V)$ denote the collection of
     Riemannian metrics $g$ on $\T^3$ which satify the following properties:
     \begin{enumerate}
				\item $g\geq g_{0}$;
				\item $\ison_{1}(g)\geq\Lambda$;
                \item $\|R^{-}\|_{L^2(g)}\leq R$;
				\item $|\T^{3}|_{g}\leq V$.
			\end{enumerate}
   Suppose that $g_i$ is a sequence of metrics in $\mathcal{V}$ such that
   \begin{equation}
       |\T^3|_{g_i}\rightarrow |\T^3|_{g_0},
   \end{equation}
   and
   \begin{equation}
       \lim_{i\rightarrow\infty}\|R^{-}_{g_i}\|_{L^2(g_i)}=0.
   \end{equation}
   Then, the metric $g_0$ must be flat.
 \end{corollary}

	\section{Background}
	
	\subsection{Convergence in the sense of Dong-Song}
	In \cite{Dong-Song-2023} Dong-Song established the stability of the Positive Mass Theorem 
	with respect to a novel notion of convergence.
	Here we will offer a slight modification of this notion.
	Let us begin by recalling the intrinsic length metric associated with a subset of a
	Riemannian manifold.
	\begin{definition}
		Let $(M,g)$ be a Riemannian manifold, and let $\Omega\subset{}M$ be a subset
		of $M$.
		Furthermore, let $\mathrm{L}_{g}(\gamma)$ denote the length of a curve $\gamma$ as 
		measured by $g$.
		For any two points $x,y$ in $\Omega$ let us define
			\begin{equation}
				\hat{d}^{g}_{\Omega}(x,y)=
				\inf\bigl\{\mathrm{L}_{g}(\gamma):\gamma\text{ connects }x\text{ and }y \text{ and }\gamma \subset \Omega\bigr\}.
			\end{equation}
	\end{definition}
	We can now use the above definition to define convergence in the sense of Dong-Song.
	\begin{definition}\label{defn:Convergence_Dong-Song}
		Let $(M,g)$ and $(N,h)$ be two closed Riemannian manifolds.
		We say that $(M,g)$ is $\varepsilon{}$-close to $(N,h)$ in the Dong-Song sense if
		there exists a connected open domain $\Omega_{M}\subset{}M$ with smooth boundary
		such that
			\begin{equation}
				|\Omega^{c}|+|\partial\Omega|\leq\varepsilon{}
			\end{equation}
		and
			\begin{equation}
				d_{GH}\left(
				(N, d^{h}),\left(\Omega_{M},\hat{d}^{g}_{\Omega_{M}}\right)\right)
				\leq\varepsilon{}.
			\end{equation}
	\end{definition}
	\begin{remark}
		Note the asymmetry in the above definition.
		If the results of \cite{lee2022metric} hold for three-tori, then the 
		results of this paper indicate that the above is not in general symmetric.
		However, notions of nearness like the above do seem to be useful for stability
		problems, see \cite{Dong-PMT_Stability} and \cite{dong2024stability}.
	\end{remark}	
	The following lemma is crucial for establishing the type of convergence given
	in Definition \ref{defn:Convergence_Dong-Song}.
	\begin{lemma}
		Let $(\T^{3},h)$ be a flat Riemannian metric on the three-torus.
		For every $\varepsilon{}>0$ there exists a $\delta>0$ such that if
		$\Omega\subset{}\T^{3}$ has smooth boundary and
			\begin{equation}
				|\Omega|\geq|\T^{3}|-\delta;
			\end{equation}
			\begin{equation}
				|\partial{}\Omega|\leq\delta,
			\end{equation}
		then there exists a connected subset $\Omega'\subset{}\Omega$ with smooth boundary
		such that
			\begin{align}
				&|\Omega'|\geq|\T^{3}|-\varepsilon{};
				\\
				&|\partial{}\Omega'|\leq\varepsilon{};
				\\
				&d_{\mathrm{GH}}\bigl((\Omega',\hat{d}^{h}_{\Omega'}),
				(\Omega',d^{h})\bigr)\leq\varepsilon.
			\end{align}			

	\end{lemma}
	\subsection{Harmonic maps and Stern's inequality}
	Let $(M,g)$ be an arbitrary closed and oriented Riemannian manifold.
	We begin by recalling the Hodge star map
	$\star:\Omega^{p}(M)\rightarrow{}\Omega^{n-p}(M)$.
		\begin{definition}
			Let $(M,g)$ be a closed and oriented Riemannian manifold, and let
			$a\in{}\Omega^{p}(M)$.
			Then, we may uniquely define an element $\star\alpha\in{}\Omega^{n-p}(M)$
			as follows.
			Set $\star\alpha$ to be the unique differential form such that
				\begin{equation}
					\int_{M}g(a,b)dV_{g}=\int_{M}\star\alpha\wedge b
				\end{equation}
			for all $b\in{}\Omega^{p}(M)$.
		\end{definition}
	Recall that for each cohomology class in $H^{k}(M;\R)$ there is an unique harmonic
	representative. In particular, for any cohomology class the harmonic representative
	has minimal $L^{2}$ norm among the representatives of the class, see
	\cite{petersen2006riemannian} for a quick introduction to harmonic
	forms and Hodge decomposition.
	In the case that $\alpha\in{}H^{1}(M;\Z)$, then there is a harmonic map
	$u:(M,g)\rightarrow{}\Sp^1$ such that $du$ is the harmonic representative of $\alpha$.

	The existence and uniqueness of harmonic representatives of cohomology classes
	provide $H^{p}(M;\R)$ with an inner product structure, as is defined below.
	\begin{definition}\label{defn:innproduct_on_cohomology}
		Let $(M,g)$ be a closed oriented Riemannian manifold, and let $\alpha$ and $\beta$
		be two cohomology classes in $H^{p}(M;\R)$.
		Furthermore, let $a$ and $b$ be the corresponding unique harmonic representatives.
		Then, we define $g(\alpha,\beta)$ to be
			\begin{equation}
				g(\alpha,\beta)=\int_{M}g(a,b)dV_{g}.
			\end{equation}
	\end{definition}
    \begin{notation}\label{ntn:L2_norm_cohomology}
        Following Hebda \cite{hebda2023stable_systoles}, for $\alpha\in H^{p}(M;\R)$ we let $|\alpha|^{*}_{2}$ denote $\sqrt{g(\alpha,\alpha)}$.
    \end{notation}
	Specializing to three dimensions, Stern \cite{stern2019scalar} connected scalar curvature to 
	harmonic forms
	with the following powerful inequality.
	\begin{theorem}[Stern's inequality]
		Let $(M^{3},g)$ be a closed and oriented 3-manifold, let $u:(M,g)\rightarrow{}\Sp^1$
		be a nontrivial harmonic map, and let $R_{g}$ denote the scalar curvature of $g$.
		For the level sets $\Sigma_{\theta}=u^{-1}\{\theta\}$ we let 
		$\chi(\Sigma_{\theta})$ denote the Euler characteristic of $\Sigma_{\theta}$.
		Then, we have that
			\begin{equation}
				2\pi\int_{\Sp}\chi\left({\Sigma_{\theta}}\right)d\theta\geq
				\frac12\int_{\Sp}\int_{\Sigma_{\theta}}
				(|du|^{-2}|\nabla du|^{2}+R_{g})dA_{g}d\theta.
			\end{equation}
	\end{theorem}
	By inspection, we see that if we could control $\chi(\Sigma_{\theta})$, then we would 
	have a very strong relationship between $R^{-}_{g}$ and the Hessian of non-trivial
	maps into $\Sp^1$.
	The following lemma is crucial in this regard.
	\begin{lemma}[\cite{stern2019scalar}]\label{lem:preimage_nontrivial}
		Let $(M,g)$ be a closed oriented Riemannian manifold, and let 
		$u:(M,g)\rightarrow{}\Sp^1$ be a nontrivial harmonic map.
		Then, for almost every $\theta\in{}\Sp^1$ we have that $\Sigma_{\theta}=u^{-1}(\theta)$
		is smooth, and each component is a non-trivial element of $H_{n-1}(M)$.
	\end{lemma}
	\begin{proof}
		That $\Sigma_{\theta}$ is smooth for almost every $\theta\in{}\Sp^1$ is a 
		standard consequence of Sard's Lemma. In particular, we may also assume that
		$|\nabla u|_{\Sigma_{\theta}}>0$ for almost every $\theta\in{}\Sp^1$.
		Therefore, we will focus on the second statement.
		Let $\theta_{0}\in{}\Sp^1$ be such that $\Sigma=\Sigma_{\theta_{0}}$ is smooth
		and is the disjoint union of its connected components: 
		$\Sigma=\bigsqcup{}\Sigma_{i}$.
		Let $\Sigma_{i_0}$ be an arbitrary component of $\Sigma$, and suppose that
		it is a trivial element in $H_{n-1}(M)$.
		Since $u$ is harmonic, we have that $\star du$ is also harmonic, and in 
		particular is a closed element of $H^{n-1}(M;\R)$.
		Therefore, since $\Sigma_{i_{0}}$ is trivial, we have that
			\begin{equation}
				\int_{\Sigma_{i_{0}}}\star{}du=0.
			\end{equation}
		However, from the definition of $\star{}du$ we see that
			\begin{equation}
				\int_{\Sigma_{i_{0}}}\star{}du=
				\int_{\Sigma_{i_{0}}}|\nabla u|dA_{g}.
			\end{equation}
		Since $|\nabla u|_{\Sigma}>0$, this leads us to a contradiction, and so
		it follows that $\Sigma_{i_{0}}$ could not have been a trivial element of 
		$H_{n-1}(M)$.
	\end{proof}
	With this lemma and Stern's inequality in hand, we immediately get the following
	important corollary.
	\begin{corollary}\label{cor:stern_ineq}
		Let $(M^{3},g)$ be a closed oriented three dimensional Riemannian manifold,
		let $u:(M,g)\rightarrow{}\Sp^1$ be a nontrivial harmonic map, and let $R^{-}_{g}$
		be the negative part of the scalar curvature.
		Finally, suppose that $H_{2}(M)$ has no non-separating 2-spheres.
		Then, we have that
			\begin{equation}
				\|R^{-}_{g}\|_{L^{2}(g)}\|du\|_{L^{2}(g)}\geq
				\int_{M}\frac{|\nabla{}du|^{2}}{du}dV_{g}.
			\end{equation}
	\end{corollary}		
	\begin{proof}
		It follow from our hypotheses, Lemma \ref{lem:preimage_nontrivial}, and the
        classification of surfaces that
		for almost every $\theta\in{}\Sp$ we have
			\begin{equation}
				\chi(\Sigma_{\theta})\leq0.	
			\end{equation}
		Therefore, we may rearrange Stern's inequality to obtain
			\begin{equation}
				-\int_{\Sp}\int_{\Sigma_{\theta}}R_{g}dA_{g}d\theta\geq
				\int_{\Sp}\int_{\Sigma_{\theta}}
				\frac{|\nabla{}du|^{2}}{|du|^{2}}dA_{g}d\theta.
			\end{equation}
		The result now follows from an application of the coarea formula on both sides,
		the definition of $R^{-}_{g}$ and an application of H{\"o}lder's inequality
		on the left hand side.
	\end{proof}
    At this stage it is convenient to introduce the following notation.
    \begin{notation}
	   For $\Sp^{1}$ let $d\theta^{2}$ denote the metric for which $\Sp^{1}$ has length 1.
	   Then, we let $h=(d\theta^{1})^2+\cdots+(d\theta^{n})^{2}$ be the product metric
	   on $\T^{n}$.
    \end{notation}
	Recall that to each element of $H^{1}(\T^{3};\Z)$ we may associate a map to $\Sp^{1}$,
	and so we see that every element of $H^{1}(\T^{3};\Z)_{\R}\subset{}H^{1}(\T^{3};\R)$
	corresponds to an harmonic map from $(M,g)$ to $\Sp^1$, which is unique up to translation.
	Therefore, to any three elements $\alpha^{i}\in{}H^{1}(M;\R)$ we may find three
	harmonic maps $u^{i}$ such that $[du^{i}]=\alpha^{i}$, and from these three maps
	we get an harmonic map $\U:(\T^{3},g)\rightarrow{}(\T^{3},h)$ defined by
		\begin{equation}
			\U(x)=\bigl(u^{1}(x),u^{2}(x),u^{3}(x)\bigr)
		\end{equation}	
	for $x\in{}T^{3}$.
    From this expression we see that
        \begin{equation}
            d\U=\left(du^1,du^2,du^3\right).
        \end{equation}
    This immediately leads us to the following proposition, which we will often make
    use of without further comment.
	\begin{proposition}\label{prop:form_of_differential}
		Let $\U:(\T^{3},g)\rightarrow{}(\T^{3},h)$, let $x\in{}M$, and let
		$\nu\in{}T_{x}\T^{3}$.
		Then, we have that
			\begin{equation}
				d\U(\nu)=du^{i}(\nu)\frac{\partial{}}{\partial{}\theta^{i}}=
				g(\nabla{}u^{i},\nu)\frac{\partial{}}{\partial{}\theta^{i}}
			\end{equation}
        and
            \begin{equation}
                \nabla d\U=\left(\nabla du^1,\nabla du^2,\nabla du^3\right).
            \end{equation}
	\end{proposition}
    Due to the above proposition, we see that the Hessian of $\U$ will be controlled in terms of Stern's inequality in the form
	of Corollary \ref{cor:stern_ineq}, which is the heart of this paper.
	\subsection{Lattices and Successive Minima}
	As may be suspected from the above section, it is important to analyze how
	$H^{1}(\T^{n};\Z)$ sits inside of $H^{1}(\T^{n};\R)$ as a lattice.
	Let us now fix some notation and terminology which will be helpful in this pursuit.
	\begin{notation}
		Let $(M,g)$ be a closed and oriented Riemannian manifold.
        We denote by $H^{p}(M;\Z)_{\R}$ the lattice in $H^{p}(M;\R)$
        generated by $H^{p}(M;\Z)$.
		If $H^{p}(M;\Z)$ is free abelian, then $H^{p}(M;\Z)\simeq H^{p}(M;\Z)_{\R}$.
        Let us give $H^{p}(M;\R)$ the inner product structure mentioned in Definition \ref{defn:innproduct_on_cohomology}.
		Then, we define $\mathrm{det}\Bigl(H^{p}(M;\Z)_{\R}\Bigr)$ to be the determinant of
		the lattice $H^{p}(M;\Z)_{\R}\subset{}H^{p}(M;\R)$ with respect to the inner product
		on $H^{p}(M;\R)$.
	\end{notation}
	To begin analyzing $H^{p}(\T^{n};\Z)_{\R}$ as a lattice of $H^{p}(\T^{n};\R)$, one may use 
	Poincare Duality and the free-ness of $H^{p}(\T^{n};\Z)$ to show that tori have the 
	following property.
	\begin{lemma}[Berger \cite{berger1972ombre}]
		For any $p$ and any $\alpha\in{}H^{p}(\T^{n};\R)$ we have that if
		$(\alpha\cup{}\beta)[\T^{n}]$ is in $\Z$ for all $\beta$ in $H^{n-p}(\T^n;\Z)$,
		then $\alpha$ is in $H^{p}(\T^n;\Z)$.
	\end{lemma}
	\begin{definition}
		Let $M$ be any smooth closed manifold, and fix $p\in{}\mathbb{N}$.
		We say that $M$ satisfies the dual lattice condition in degree $p$ if the
		conclusion of the above lemma holds for $p$-forms.
	\end{definition}
	We now come to an important result relating the determinants of different cohomology
	groups to each other.
	\begin{lemma}[Berger \cite{berger1972ombre}]\label{lem:det_dual_lat}
		Let $(M,g)$ be a closed oriented $n$-dimensional Riemannian manifold which
        satisfies the dual lattice condition in degree $p$.
		Then, we have that
			\begin{equation}
				1=\mathrm{det}\Bigl(H^{p}(M;\Z)_{\R}\Bigr)\mathrm{det}\Bigl(H^{n-p}(M;\Z)_{\R}
                \Bigr).
			\end{equation}
	\end{lemma}
	In order to make good use of Lemma \ref{lem:det_dual_lat}, we need a few results
	from the Geometry of Numbers and Systolic Geometry.
    Let us begin by reviewing the Geometry of Numbers.
	In what follows, most definitions and results have generalizations which are not needed
	for this paper, but can be found in the relevant sections of \cite{cassels2012introduction}.
	\begin{definition}
		Let $L$ be a lattice in $\R^{n}$, and let $F_{0}(x)$ denote the Euclidean norm.
		Then, we define
			\begin{equation}
				F_{0}(L)=\inf\{F_{0}(a):a\in{}L\}.
			\end{equation}
	\end{definition}
	We can go a step further to define the following quantity associated with $F_{0}$.
	\begin{definition}
		For $F_{0}$ as above, we define $\delta_{0}$ as follows:
			\begin{equation}
				\delta_{0}=
				\sup\left\{\frac{F_{0}(L)}{\mathrm{det}(L)}:L
					\text{ is a lattice in }\R^{n}
				\right\}
			\end{equation}
	\end{definition}
	We have the following result bounding $\delta_{0}$.
	\begin{lemma}[\cite{cassels2012introduction}]\label{lem:delta0_bound}
		Let $|B(0,1)|$ be the volume of the unit ball in $\R^{n}$.
		Then, we have the following upper bound on $\delta_{0}$:
		\begin{equation}
			\delta_{0}\leq\frac{2^{n}}{|B(0,1)|}.	
		\end{equation}
	\end{lemma}
	\begin{proof}
		From \cite[Theorem 1 Chapter 3]{cassels2012introduction}, we know that if
		$|B(0,r)|$ is greater than $\mathrm{det}(L)$, then there exists two points in
		$B(0,r)$, say $x_{1}$ and $x_{2}$, such that $x_{1}-x_{2}\in{}B(0,r)\cap{}L$,
		and so we see that $F_{0}(L)\leq2r$.
		Therefore, if we choose $r_{L}=\sqrt{n}{\frac{\mathrm{det}(L)}{|B(0,1)|}}$, then
		we get that $F_{0}(L)\leq2r_{L}$ for all lattices $L$.
		This gives the desired bound.
	\end{proof}
	For any given lattice $L\subset{}\R^{n}$ the quantity $F_{0}(L)$ has the following useful
	generalization.
	\begin{definition}
		Let $L\subset{}\R^{n}$ be a lattice, and let $F_{0}$ denote the Euclidean norm.
		We denote by $\lambda_{k}$ the following quantity:
			\begin{equation}
				\lambda_{k}=\inf\bigl\{\lambda:\exists
				\text{ linearly independent }
				\{\nu_{i}\}_{i=1}^{k}\subset{}L
				\text{ such that }
				F_{0}(\nu_{i})\leq\lambda \text{ for any }i
			\bigr\}.
			\end{equation}
			We refer to $\lambda_{k}$ as the $k^{th}$ \textit{successive minima}
			of $L$ with respect to $F_{0}$. Observe that 
            $\lambda_1=F_0(L)$.
	\end{definition}
	The importance of successive minima for this paper is contained in the following lemma.
	\begin{lemma}\label{lem:bounded_lat_basis}
		Let $L$ be a lattice in $\R^{n}$, and let $\lambda_{1},\dots,\lambda_{n}$ be the
		successive minima of $L$ with respect to $F_{0}$.
		Then, there is a basis of $L$, say $\{b_{i}\}_{i=1}^{n}$ such that
			\begin{align}
				&|b_{1}|=\lambda_{1}
				\\
				&|b_{j}|\leq\frac12j\lambda_{j}\quad(2\leq j\leq n).
			\end{align}
	\end{lemma}
	\begin{proof}
		One may apply \cite[Chapter 8 Lemma 1]{cassels2012introduction} to find $n$
		linearly independent elements of the lattice $a_{1},\dots,a_{n}$ such that
		$|a_{j}|=\lambda_{j}$ for all $j=1,\dots,n$.
		Then, we may apply \cite[Chapter 5 Lemma 1]{cassels2012introduction} to find
		a basis $b_{1},\dots,b_{n}$ such that $|b_{1}|=\lambda_{1}$ and for each
		$j\geq2$ we have
			\begin{equation}
				|b_{j}|\leq\max\left\{|a_{j}|,\frac12\sum_{i=1}^{j}|a_{i}|\right\}.
			\end{equation}	
	\end{proof}
	\begin{remark}
		In fact, more can be said since we are working with the Euclidean norm,
		see \cite{remak1938minkowskische} and \cite{van1956reduktionstheorie}.
	\end{remark}
	The above shows that we may always find a basis for a lattice whose norms are bounded
	by the successive minima of the lattice.
	The following result is the key to estimating these successive minima.
	\begin{theorem}[Chapter 8 Theorem 1 \cite{cassels2012introduction}]
		Let $L$ be a lattice in $\R^{n}$, and let $\lambda_{1},\dots,\lambda_{n}$ be its
		successive minima with respect to $F_{0}$, the Euclidean norm.
		Then, we have that
			\begin{equation}
				\mathrm{det}(L)\leq\lambda_{1}\cdots\lambda_{n}\leq\delta_{0}
				\cdot\mathrm{det}(L).
			\end{equation}
	\end{theorem}
	Combining this theorem with Lemma \ref{lem:delta0_bound} gives us the following
	useful corollary.
	\begin{corollary}\label{cor:lat_minima_det_ineq}
		Let $L\subset{}\R^{n}$ be a lattice, and let $\lambda_{1},\cdots,\lambda_{n}$ be
		its successive minima with respect to the Euclidean norm.
		Then, we have that
			\begin{equation}
				\lambda_{1}\cdots\lambda_{n}\leq\frac{2^{n}}{|B(0,1)|}
				\mathrm{det}(L).
			\end{equation}
	\end{corollary}

\subsection{Stable Systoles}
We now turn to a quick review of some concepts in Systolic Geometry.
\begin{definition}[Stable Norm of a Real Homology Class]
    The volume of a real $k$-dimensional Lipschitz cycle $c= \sum_i r_i \sigma_i$ is given by
\begin{equation*}
 \vol_k(c) = \sum_i |r_i| \vol_k (\triangle^k, \sigma_i^{*} g). 
\end{equation*}
The stable norm $\| \alpha \|$ of a real homology class $\alpha \in H_k(M; \R)$ is defined as the infimum of the volumes of all real Lipschitz cycles representing $\alpha$.
\end{definition}
\begin{definition}[Stable Systoles]
    The \textit{stable k-systole}, denoted $\stsys_k(M, g)$, is defined to be the minimum of the stable norm on the nonzero classes of the integral lattice $H_k(M; \Z)_{\R}$ in $H_k(M; \R)$. 
\end{definition}
\begin{notation}\label{ntn:Poincare_Dual}
    Let $\alpha\in H^{n-p}(M;\R)$, we shall denote by $\pd(\alpha)$ its Poincar{\'e} dual in
    $H_{p}(M;\R)$.
\end{notation}
The following proposition is \cite[Corollary 3]{Hebda-Collars}, see \cite[Proposition 3.2]{hebda2023stable_systoles} for its statement using the stable norm. 
\begin{proposition} \label{Federa:stable}
 Let $\pd(\alpha) \in H_p(M; \R)$ be the Poincar{\'e} dual of the cohomology class $\alpha \in H^{n-p}(M; \R)$, and let $\|\pd(\alpha)\|$ be the stable norm of $\pd(\alpha)$.
 Then
 \begin{equation*}
  \| \pd(\alpha)\| \leq |M|_{g}^{1/2} C(n, p) |\alpha|_2^{*},
 \end{equation*}
 where $C(n, p)$ is a constant depending only on $n$ and $p$, and $|\alpha|_2^{*}$ is the $L^2$
 norm of $\alpha$, see Definition \ref{defn:innproduct_on_cohomology} and Notation \ref{ntn:L2_norm_cohomology}.
\end{proposition}
The above proposition has a simple, but important corollary:
\begin{corollary}\label{cor:latice_constant_systole_bound}
    Let $(M,g)$ be a closed Riemannian manifold.
    Then we have that
    \begin{equation}
        \stsys_{p}(M,g)\leq|M|_{g}^{\frac12}\min\Bigl\{|\alpha|^{*}_{2}:\alpha\neq0;
        \alpha\in H^{n-p}(M;\Z)_{\R}
        \Bigr\}=|M|^{\frac12}_{g}F_0\left(H^{n-p}(M;\Z)_{\R}\right).
    \end{equation}
\end{corollary}
\begin{proof}
 According to Proposition~\ref{Federa:stable}, for $\alpha \in H^p(M;\Z)_{\R}$,
 $1 \leqslant p \leq n-1$,
 \begin{align*}
  \| \pd(\alpha) \| \leq  |M|_{g}^{1/2} C(n, p) |\alpha|_2^{*},
 \end{align*}
 Since the Poincar{\'e} dual map is an isomorphism, we have that
 $\pd(\alpha)\neq0$ in $H_{p}(M;\Z)_{\R}$.
 Therefore, by the
 definition of the stable $p$-systole, we have that
 \begin{equation}
     \stsys_{p}(M,g)\leq\|\pd(\alpha)\|.
 \end{equation}
 This gives the result.
\end{proof}
	\subsection{Isoperimetric constants}
    On a smooth Riemannian manifold, functions are bounded in terms of their gradients. The character and quality of this bound can be determined by the character and quality of isoperimetric bounds.
	Here we recall the definition of the Cheeger and Sobolev constants, and the fact
	that they are closely related.
	\begin{definition}\label{defn:Ison}
		Let $(M,g)$ be a given $n$-dimensional Riemannian manifold, and let
		$\alpha\in{}[1,\tfrac{n}{n-1}]$.
		We denote by $\ison_{\alpha}(M,g)$ the following quantity:
			\begin{equation}
				\ison_{\alpha}(M,g)=\inf\left\{
				\frac{|\partial{}\Omega|}{\min\{|\Omega|,|\Omega^{c}|\}}:
				\Omega\subset{}M
				\right\}.
			\end{equation}
		When $\alpha=1$, we call $\ison_{1}(M,g)$ the Cheeger constant of $(M,g)$.
	\end{definition}
	Next, we have the Sobolev constant of a Riemannian manifold.
	\begin{definition}
		Let $(M,g)$ be a given $n$-dimensional Riemannian manifold, and let
		$\alpha\in{}[1,\frac{n}{n-1}{}]$.
		Let us denote by $\sobn_{\alpha}(M,g)$ the following quantity:
			\begin{equation}
				\sobn_{\alpha}(M,g)=\inf\left\{
				\frac{\int_{M}|\nabla f|dV_{g}}
				{\inf_{k\in{}\R}\|f-k\|_{L^{\alpha}(g)}}{}:f\in{}W^{1,1}(M,g)
			\right\}.
			\end{equation}
	\end{definition}
	It is standard, see \cite[Theorem 9.6]{Li_2012}, that the Cheeger constant and $\sobn_{1}$ are equivalent:
	\begin{proposition}
		Let $(M,g)$ be a given $n$-dimensional Riemannian manifold, then we have that
			\begin{equation}
				\ison_{1}(M,g)=\sobn_{1}(M,g).
			\end{equation}
	\end{proposition}

	\section{Stable systole bounds and $L^{2}$ estimates}\label{sec:stab_systole_bounds_L2}
	
In this section we will show that lower bounds on $\stsys_{1}(M,g)$ and $\stsys_{n-1}(M,g)$
guarantee the existence of a good basis for $H^{1}(\T^{n};\Z)_{\R}$.
The following lemma and its proof are modeled on \cite[Proposition 6]{Hebda-Collars}.
\begin{lemma}\label{lem:prod_min_upper_bound}
	Let $g$ be a Riemannian metric on $\T^{n}$ such that
		\begin{equation}
			\stsys_{1}(g)\geq\sigma>0
		\end{equation}
	and let $\lambda_{1},\dots,\lambda_{n}$ be the successive minima of the lattice
	$H^{1}(\T^{n};\Z)_{\R}$ in $H^{1}(\T^{n};\R)$ with respect to the inner product on
	$H^{1}(\T^{n};\R)$ induced by $g$, see Definition \ref{defn:innproduct_on_cohomology}.
	Then, we have that
		\begin{equation}
			\lambda_{1}\cdots\lambda_{n}\leq2^{2n}|B(0,1)|^{2}
			\sigma^{-n}|\T^{n}|^{\frac{n}{2}{}}_{g}.
		\end{equation}
	In particular,
		\begin{equation}\label{eq:lambda1_upper_bound}
			\lambda_{1}\leq4|B(0,1)|^{\frac{2}{n}{}}\sigma^{-1}|\T^{n}|^{\frac12}_{g}.
		\end{equation}
\end{lemma}
\begin{proof}
	From Corollary \ref{cor:lat_minima_det_ineq} we see that
		\begin{equation}
			\frac{|B(0,1)|}{2^{n}}{}\lambda_{1}\cdots\lambda_{n}\leq
			\mathrm{det}\Bigl(H^{1}(\T^{n};\Z)_{\R}\Bigr).
		\end{equation}
	Letting $\mu_1,\dots,\mu_{n}$ be the successive minima of $H^{n-1}(\T^n;\Z)_{\R}$,
    we get the following as well:
		\begin{equation}
			\frac{|B(0,1)|}{2^{n}}{}\mu_{1}\cdots\mu_{n}\leq
			\mathrm{det}\Bigl(H^{n-1}(\T^{n};\Z)_{\R}\Bigr).
		\end{equation}
	Since $\T^n$ satisfies the dual lattice condition in all degrees,
    it follows from Lemma \ref{lem:det_dual_lat} that
		\begin{equation}
			2^{-2n}|B(0,1)|^{2}\lambda_{1}\cdots\lambda_{n}\cdot
			\mu_{1}\cdots\mu_{n}\leq1.
		\end{equation}
	Since $\mu_{1}\leq\mu_{j}$ for all $j\geq1$, it follows that
		\begin{equation}
			\lambda_{1}\cdots\lambda_{n}\leq2^{2n}|B(0,1)|^{2}\mu_{1}^{-n}.
		\end{equation}
	In fact, since $\lambda_{1}\leq\lambda_{j}$ for all $j\geq1$ as well, we have
		\begin{equation}
			\lambda_{1}^{n}\leq2^{2n}|B(0,1)|^{2}\mu_{1}^{-n}.
		\end{equation}
	Therefore, a lower bound on $\mu_{1}$ will give us our desired upper bound on
	$\lambda_{1}$.
    From Corollary \ref{cor:latice_constant_systole_bound}, we have that
    \begin{equation}
        \stsys_1(\T^n,g)\leq|\T^n|_{g}^{\frac12}\min\Bigl\{|\alpha|^{*}_{2}:\alpha\neq0;
        \alpha\in H^{n-1}(\T^n;\Z)_{\R}\Bigr\}
        =|\T^n|_{g}^{\frac12}\mu_1.
    \end{equation}
    This gives the desired lower bound on $\mu_1$, and so the result follows.
\end{proof}
We now know that $\lambda_{1}=\inf\bigl\{|b|^*_{2}:b\neq0;b\in{}H^{1}(\T^{n};\Z)_{\R}\bigr\}$ is
bounded above in terms of $|\T^{n}|_{g}$ and $\sigma^{-1}$.
However, in order to construct a useful harmonic map, it seems reasonable to suppose that
we need to use $n$ one-forms, which together form a basis for $H^{1}(\T^n;\Z)_{\R}$.
It turns out that a lower bound on $\min\{\stsys_{1}(g),\stsys_{n-1}(g)\}$ is sufficient to
ensure the existence of such a basis.
\begin{lemma}
	Let $g$ be a Riemannian metric on $\T^{n}$ such that
		\begin{equation}
			\min\{\stsys_{1}(g),\stsys_{n-1}(g)\}\geq\sigma>0.	
		\end{equation}
	Then, there exists a basis $\alpha_{1},\dots,\alpha_{n}$ of $H^{1}(\T^{n};\Z)_{\R}$ with
	harmonic representatives $a_{j}$ such that
		\begin{equation}
			\|a_{j}\|_{L^{2}(g)}\leq j\times\sqrt[(n-j+1)]
			{2^{2n}|B(0,1)|\sigma^{-(n+j-1)}
			|\T^{n}|^{\frac{1}{2}{}(n+j-1)}_{g}}.
		\end{equation}
\end{lemma}
\begin{proof}
    We wish to apply Lemma \ref{lem:bounded_lat_basis}, however first we must bound 
    each successive minima.
    Equation \eqref{eq:lambda1_upper_bound} gives the desired upper bound on
    $\lambda_1$.
    Therefore we must focus our attention on the higher successive minima.
	Using that $\lambda_{j}\leq\lambda_{i}$ for $j \leq i$, we may apply 
	Lemma \ref{lem:prod_min_upper_bound} to obtain
		\begin{equation}
			\lambda_{j}^{n-j+1}\leq\lambda_{j}\cdots\lambda_{n}\leq
			\lambda_{1}^{-(j-1)}2^{2n}|B(0,1)|^{2}\sigma^{-n}|\T^{n}|^{\frac{n}{2}{}}_{g}.
		\end{equation}
    From Corollary \ref{cor:latice_constant_systole_bound}, we have that
    \begin{equation}\label{eq:lower_bound_min_H1}
        \stsys_{n-1}(\T^n,g)\leq|\T^n|_{g}^{\frac12}\min\Bigl\{|\alpha|^*_{2}:
        \alpha\neq0;
        H^1(\T^n;\Z)_{\R}
        \Bigr\}
        =|\T^n|_{g}^{\frac12}\lambda_1.
    \end{equation}
	Thus, using Equation \eqref{eq:lower_bound_min_H1}, we see that
		\begin{equation}
			\lambda_{j}^{n-j+1}\leq2^{2n}|B(0,1)|^{2}\sigma^{-(n+j-1)}
			|\T^{n}|^{\frac{1}{2}{}(n+j-1)}_{g}.
		\end{equation}
	Taking $(n-j+1)^{th}$ roots and then applying Lemma \ref{lem:bounded_lat_basis}
	gives the desired result.
\end{proof}

Using the above results and notation, we can now establish the existence of a harmonic map
$\U:(\T^{n},g)\rightarrow{}(\T^{n},h)$ with several desirable properties.
\begin{corollary}\label{cor:L2-bounded_deg-1-harmonic_map}
	Let $g$ be a Riemannian metric on $\T^{n}$ such that
		\begin{equation}
			\min\{\stsys_{1}(g),\stsys_{n-1}(g)\}\geq\sigma>0.	
		\end{equation}
	Then, there exists a surjective harmonic function $\U:(\T^{n},g)\rightarrow{}(\T^{n},h)$ such that
	$\mathrm{deg}(\U)=1$ and
		\begin{equation}
			\|d\U\|_{L^2(g)}\leq\sum_{j=1}^{n}j\times
			\sqrt[(n-j+1)]{4^{n}|B(0,1)|^{2}\sigma^{-(n+j-1)}|\T^{n}|^{
			\frac{1}{2}{}(n+j-1)}_{g}}.
		\end{equation}  
\end{corollary}
\begin{proof}
	Let $\alpha^{i}$ be a basis for the lattice $H^{1}(\T^{n};\Z)_{\R}$ with harmonic representatives
	$a^{i}$ as in Lemma \ref{lem:bounded_lat_basis}, and let 
	$u^{i}:(\T^{n},g)\rightarrow{}\Sp^1$ be the harmonic map such that $du^{i}=a^{i}$.
	Then, if we let $\U=\bigl(u^{1},\dots,u^{n}\bigr)$, the estimate on
	$\|d\U\|_{L^{2}(g)}$ follows from Proposition \ref{prop:form_of_differential}.

	Next, since the $\alpha^{i}$ form a basis of $H^{1}(\T^{n};\Z)_{\R}$, it follows that the wedge product
	$\bigwedge^{n}_{i=1} a^{i}$ of their representatives $a^i$ forms a basis of $H^{n}(\T^{n};\Z)_{\R}$.
	Let $\theta^{i}$ be the $i^{th}$ coordinate function for $\T^{n}$. Then we also have
	that the wedge product $\bigwedge_{i=1}^{n}d\theta^{i}$ represents a basis for $H^{n}(\T^{n};\Z)_{\R}$.
	As such, we may calculate the degree of $\U$ as follows:
		\begin{equation}
			\mathrm{deg}(\U)=\int_{\T^{n}}\U^{*}\bigwedge^{n}_{i=1}d\theta^{i}
			=\int_{\T^{n}}\bigwedge^{n}_{i=1}a^{i}=\pm1.
		\end{equation}
  We may take $\U=(-u^1,u^2,\dots,u^n)$ as necessary to ensure that $\mathrm{deg}(\U)=1$.
  The surjectivity of $\U$ follows from the general fact that degree one maps between
  Riemannian manifolds are surjective \cite{epstein1966degree}.
\end{proof}
For the sake of completeness, let us give a simple proof of the fact that degree one maps between tori are surjective.
\begin{proposition}\label{prop:surjectivity_degree_1_maps_between_tori}
	Let $f:\T^{n}\rightarrow{}\T^{n}$ be a degree one map.
	Then, it must be true that $f$ is surjective.
\end{proposition}
\begin{proof}
	Since the result is purely topological, we may assume that $\T^{n}$ has the product metric
	$\sum_{i=1}^{n}\bigl(d\theta^{i}\bigr)^{2}$.
	Let $a$ be an arbitrary point in $\T^{n}$.
	Using the flat metric, we can see that we may give $\T^{n}$ the structure of a CW-complex
	with one $n$-cell which contains $a$ in its interior.

	With the above reduction in place, let $\alpha\in{}H_{n}(\T^{n};\Z)$ be a fundamental
	class with respect to which $f:\T^{n}\rightarrow{}\T^{n}$ is a degree one map:
		\begin{equation}
			f_{*}\alpha=\alpha.
		\end{equation}
	Suppose that $f$ is not surjective, and so there is a point $a\in{}\T^{n}$ such that
	$a$ is not in the image of $\T^{n}$ under $f$.
	Since $f(\T^{n})$ is compact, it follows that there is an $\varepsilon{}>0$ such that
		\begin{equation}
			d\Bigl(a,f\bigl(\T^{n}\bigr)\Bigr)\geq\varepsilon{}.
		\end{equation}
	As stated above, we may give $\T^{n}$ the structure of a CW complex with one $n$
	cell, which contains $a$ in its interior.
	Then, by shrinking $\varepsilon{}$ as necessary, we may assume that
	$B(a,\varepsilon{})\subset{}\mathrm{int}D^{n}$, where $D^{n}$ is the $n$-cell.
	We also have that $f(\T^{n})\subset{}B(a,\varepsilon{})^{c}$.
	As such, we see that there is a homotopy of the map, say $f_{t}$, such that
	$f_{0}=f$ and $f_{1}(\T^{n})$ is contained in the $n-1$ skeleton of $\T^{n}$.
	
	Let $X^{n-1}$ denote the $n-1$ skeleton of $\T^{n}$ and let 
	$\iota:X^{n-1}\rightarrow{}\T^{n}$ denote the inclusion map; we have that 
	$H_{n}(X^{n-1};\Z)=0$.
	In particular, we see that $(f_{1})_{*}\alpha=0$ in $X^{n-1}$.
	Using the inclusion map, we have that $f_{1}:\T^{n}\rightarrow{}\T^{n}$ has the
	following decomposition, which is almost tautological:
	\begin{equation}
		\T^{n}\xrightarrow{f_{1}}X^{n-1}\xrightarrow{\iota}\T^{n}. 	
	\end{equation}
	This shows that $(f_{1})_{*}\alpha=0\in{}H_{n}(\T^{n};\Z)$.
	However, since $f_{1}$ is homotopic to $f_{0}=f$, this shows that
	$\alpha=f_{*}\alpha=0$.
	This is a contradiction, since $\alpha$ was assumed to be a fundamental class.
\end{proof}

	\section{From $L^{2}$ to $L^{3}$ bounds}\label{sec:L2_to_L3}
	
	From the previous section we know that given a Riemannian metric $g$ on $\T^{n}$,
	we may find an harmonic map $\U:(\T^{n},g)\rightarrow{}(\T^{n},h)$ whose $L^{2}$
	energy is controlled in terms of the volume of $|\T^{n}|_{g}$ and 
	$\min\{\stsys_{1}(g),\stsys_{n-1}(g)\}$.
	In this section we will ultimately focus on the case that $n=3$, and study the relationship
	between the universal cover $\R^{3}$ with its pullback metric $\pi^{*}g$ and 
	$(\T^{3},g)$.
	The goal is to use this relationship along with Stern's inequality, in particular
	Corollary \ref{cor:stern_ineq}, to improve the $L^{2}$ bounds of the previous 
	section to $L^{3}$ bounds.

	Of fundamental importance to this discussion is the notion of fundamental domain,
	which we recall now.
	\begin{definition}\label{defn:fund_domain}
		Let $\pi:\R^{n}\rightarrow{}\T^{n}$ be the covering map.
		Suppose that $\V\subset{}\R^{n}$ is a closed subset such that
		$\pi\bigl(\V\bigr)=\T^{n}$ and $\left.\pi\right|_{\mathrm{int}\V}$ is 
		injective.
		If in addition $\V$ is path connected, and $\partial{}\V$ has measure
		zero, then we call $\V$ a fundamental domain of $\T^{n}$.
		Given a fundamental domain $\V$, for $\nu\in{}\Z^{n}$ we let 
		$\V^{\nu}$ denote the image of $\V$ under the deck transformation
		associated to $\nu\in{}\Z^{n}=\pi_{1}(\T^{n})$.
	\end{definition}
	The following proposition lists some of the basic properties of fundamental domains.
	\begin{proposition}\label{prop:lift_of_map_and_estimates}
		Let $g$ be a Riemannian metric on $\T^{n}$, let 
		$\U:\T^{n}\rightarrow{}\T^{n}$ be a map, and let $u^{k}$ denote the components
		of $\U$.
		Then, there exists lifts $\hat{u}^{k}$, and so a lift $\hat{\U}$, such that
		the following diagram commutes:
			\begin{center}\label{diag:lift_of_uk}
			    \begin{tikzcd}
				\R^n \arrow[d,"\pi"] \arrow[r,"\hat{u}^{k}"] & \R \arrow[d,"\pi"]
				\\ 
				\T^n \arrow[r,"u^{k}"] & \Sp^{1}
			    \end{tikzcd}
		    \end{center}
		Furthermore, for any integrable function $f:\T^{n}\rightarrow{}\R$ and fundamental
		domain $\V$ we have that
			\begin{equation}
				\int_{\V}f\circ\pi dV_{\pi^{*}g}=\int_{\T^{n}}fdV_{g}.
			\end{equation}
	\end{proposition}
	\begin{proof}
	Note that $\T^n = \R^n / \Z^n$ and $\Sp^1 = \R/\Z$.
    For the moment, let us denote elements of $\T^n$ by $[x]$, where $x\in\R^n$:
    we have that $[x]=\Z^n x$.
    We denote elements in $\Sp^1$ similarly.
    We can define a lift $\hat{u}^{k}$ as follows.
    First, denote by $u^k_{*}:\Z^n\rightarrow\Z$ the homomorphism between fundamental
    groups induced by $u^k$.
    Next, fix $0$ in $\R^n$, and consider the points
    $[0]$ in $\T^n$ and $u^{k}\bigl([0]\bigr)$.
    Pick any point $y_0\in\R$ such that $u^{k}\bigl([0]\bigr)=\Z y_0$.
    Abusing notation, we may use $u^k_{*}(\alpha)$ to define $\hat{u}^{k}(x)$ for any $\alpha\in\Z^n.$
    Every element $x\in\R^n$ can be written as $x=\Tilde{x}+\alpha$, where $\Tilde{x}$ is in
    the same fundamental domain as $0$, and $\alpha\in\Z^n$.
    Then, we may set $\hat{u}^{k}(x)=u^k_*(\alpha)+y_0,$
    $\hat{u}^{k}(\alpha(x)) = u^{k}_{*}(\alpha) + \hat{u}^{k}(x)$ for $\alpha \in \Z^n.$ 
    The second part follows from the fact that $\partial{}\V$ has measure zero,
	that $\left.\pi\right|_{\mathrm{int}\V}$ is injective, and from the fact that
    the image of measure zero sets under $\pi$ have measure zero, since 
    $\pi$ is a smooth map.
	\end{proof}
	Working with fundamental domains allows us to treat maps $\U:\T^{3}\rightarrow{}\T^{3}$
	as maps $\hat{\U}:\R^{3}\rightarrow{}\R^{3}$.
	The next lemma is an example of this, and will play a vital role in strengthening
	our $L^{2}$ bounds to $L^{3}$ bounds.
	\begin{lemma}\label{lem:u_bounded_fund_domain}
		Let $g$ be a Riemannian metric on $\T^{3}$ such that $\stsys_{2}(g)\geq\sigma>0$,
		let $u:(\T^{3},g)\rightarrow{}\Sp^{1}$ be a nontrivial harmonic map,
		let $\hat{u}$ denote its lift, and let $\V$ be any fundamental domain in
		$\R^{3}$.
		Then, we have that
			\begin{equation}
				\max_{\V}\hat{u}-\min_{\V}\hat{u}\leq
				\sigma^{-1}|\T^{3}|_{g}^{\frac12}\|du\|_{L^{2}(g)}.
			\end{equation}
	\end{lemma}
	\begin{proof}
		Here $\chi_{\V}$ will denote the indicator function for the set $\V$.
		From the coarea formula, we have that
			\begin{equation}
				\int_{\V}|d\hat{u}|dV_{\pi^{*}g}=
				\int_{\min_{\V}\hat{u}}^{\max_{\V}\hat{u}}
				\int_{\hat{u}^{-1}\{t\}}\chi_{\V}dA_{\pi^{*}g}dt.
			\end{equation}
		For every $t$, let $\theta(t)=t\mod 1$, then it follows from the 
		commutativity of the Diagram \ref{diag:lift_of_uk} and the fact that
		$\pi\bigl(\V\bigr)=\T^{3}$ that for all $t$ we have
			\begin{equation}
				\int_{\hat{u}^{-1}\{t\}}\chi_{\V}dA_{\pi^{*}g}\geq
				\int_{u^{-1}\{\theta(t)\}}dA_{g}.
			\end{equation}
		In particular, observe that for almost every $t$ we have that
		$t$ is a regular value for $\hat{u}$ and $\theta(t)$ is a regular
		value for $u$.
		From Lemma \ref{lem:preimage_nontrivial} it follows that for almost every
		$t$ the surface $u^{-1}\{\theta(t)\}$ is a nontrivial element of
		$H_{2}(\T^{3})$.
		Therefore, by hypothesis we have that
			\begin{equation}
				|u^{-1}\{\theta(t)\}|_{g}\geq\sigma.
			\end{equation}
        Furthermore, since $\pi:\V\rightarrow\T^3$ is surjective, it follows from
        the above inequality that for almost every $t$ we have
            \begin{equation}
                |\hat{u}^{-1}\{t\}\cap\V|_{\pi^*g}\geq|u^{-1}\{\theta(t)\}|_g\geq\sigma.
            \end{equation}
		It now follows from the coarea formula that
			\begin{equation}
				\int_{\V}|d\hat{u}|dV_{p^{*}g}\geq
				\sigma\bigl(\max_{\V}\hat{u}-\min_{\V}\hat{u}\bigr).
			\end{equation}
		We observe that $|d\hat{u}|=|d(u\circ\pi)|=|du|\circ\pi$, and so from Proposition
		\ref{prop:lift_of_map_and_estimates} we have that
			\begin{equation}
				\int_{\T^{3}}|du|dV_{g}\geq\sigma
				\bigl(\max_{\V}\hat{u}-\min_{\V}\hat{u}\bigr).
			\end{equation}
		After applying H{\"o}lder's inequality to the left hand side,
		we get the desired result.
	\end{proof}
	We will actually need to understand $\sup_{\V_{\eta}}\hat{u}-\inf_{\V_{\eta}}\hat{u}$,
	where $\eta>0$ and $\V_{\eta}$ denotes the $\eta$ neighborhood of $\V$ with respect
	to the distance function induced by the metric $\pi^{*}g$.
	In order to obtain such information, we are led to consider the following quantity.
	\begin{definition}\label{defn:covering_constant}
		Let $g$ be a Riemannian metric on $\T^{3}$, let $\pi:\R^{3}\rightarrow{}\T^{3}$
		be the covering map, and let $\eta>0$.
		We define the constant $\kappa(g,\eta)$ to be the smallest integer such that
		there exists a fundamental domain $\V$ such that
			\begin{equation}
				\sup_{x\in{}\T^{3}}|\pi^{-1}\{x\}\cap{}\V_{\eta}|\leq
				\kappa(g,\eta).
			\end{equation}
		Let us refer to $\kappa(g,\eta)$ as the $\eta$-covering constant of $g$,
		and refer to $\V$ as a test domain for $\kappa(g,\eta)$.
		See Figure \ref{fig:covering_constant} below.
	\end{definition}
	\begin{figure}[ht]
		\centering
		\def\svgwidth{5cm}
\begingroup%
  \makeatletter%
  \providecommand\color[2][]{%
    \errmessage{(Inkscape) Color is used for the text in Inkscape, but the package 'color.sty' is not loaded}%
    \renewcommand\color[2][]{}%
  }%
  \providecommand\transparent[1]{%
    \errmessage{(Inkscape) Transparency is used (non-zero) for the text in Inkscape, but the package 'transparent.sty' is not loaded}%
    \renewcommand\transparent[1]{}%
  }%
  \providecommand\rotatebox[2]{#2}%
  \newcommand*\fsize{\dimexpr\f@size pt\relax}%
  \newcommand*\lineheight[1]{\fontsize{\fsize}{#1\fsize}\selectfont}%
  \ifx\svgwidth\undefined%
    \setlength{\unitlength}{396.8503937bp}%
    \ifx\svgscale\undefined%
      \relax%
    \else%
      \setlength{\unitlength}{\unitlength * \real{\svgscale}}%
    \fi%
  \else%
    \setlength{\unitlength}{\svgwidth}%
  \fi%
  \global\let\svgwidth\undefined%
  \global\let\svgscale\undefined%
  \makeatother%
  \begin{picture}(1,0.75714286)%
    \lineheight{1}%
    \setlength\tabcolsep{0pt}%
    \put(0,0){\includegraphics[width=\unitlength,page=1]{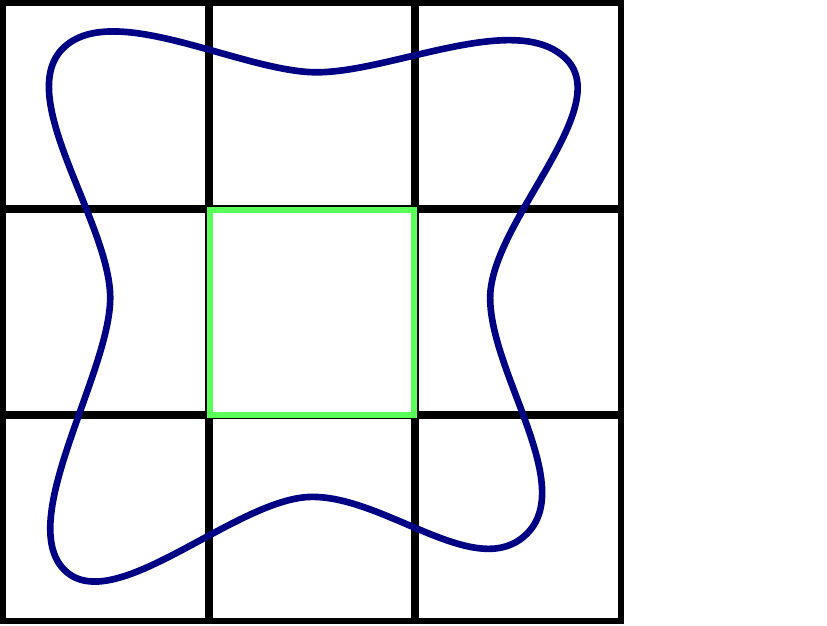}}%
    \put(0.35267018,0.29669631){\color[rgb]{0,0,0}\makebox(0,0)[lt]{\lineheight{1.25}\smash{\begin{tabular}[t]{l}$\mathbb{V}$\end{tabular}}}}%
    \put(0.33474798,0.55063491){\color[rgb]{0,0,0}\makebox(0,0)[lt]{\lineheight{1.25}\smash{\begin{tabular}[t]{l}$\mathbb{V}_{\eta}$\end{tabular}}}}%
  \end{picture}%
\endgroup%
	
		\caption[covering constant]{A fundamental domain $\V$ with neighborhood $\V_{\eta}$.
		In the case loosely depicted here, we have $\kappa(g,\eta)\leq9$.}
		\label{fig:covering_constant}
	\end{figure}
	Before we can estimate $\bigl(\sup_{\V_{\eta}}\hat{u}-\inf_{\V_{\eta}}\hat{u}\bigr)$ 
	more precisely, we need to understand how the different copies of $\V$ cover
	$\V_{\eta}$.
	\begin{figure}[ht]
		\centering
		\def\svgwidth{5cm}
\begingroup%
  \makeatletter%
  \providecommand\color[2][]{%
    \errmessage{(Inkscape) Color is used for the text in Inkscape, but the package 'color.sty' is not loaded}%
    \renewcommand\color[2][]{}%
  }%
  \providecommand\transparent[1]{%
    \errmessage{(Inkscape) Transparency is used (non-zero) for the text in Inkscape, but the package 'transparent.sty' is not loaded}%
    \renewcommand\transparent[1]{}%
  }%
  \providecommand\rotatebox[2]{#2}%
  \newcommand*\fsize{\dimexpr\f@size pt\relax}%
  \newcommand*\lineheight[1]{\fontsize{\fsize}{#1\fsize}\selectfont}%
  \ifx\svgwidth\undefined%
    \setlength{\unitlength}{300.47244094bp}%
    \ifx\svgscale\undefined%
      \relax%
    \else%
      \setlength{\unitlength}{\unitlength * \real{\svgscale}}%
    \fi%
  \else%
    \setlength{\unitlength}{\svgwidth}%
  \fi%
  \global\let\svgwidth\undefined%
  \global\let\svgscale\undefined%
  \makeatother%
  \begin{picture}(1,1)%
    \lineheight{1}%
    \setlength\tabcolsep{0pt}%
    \put(0,0){\includegraphics[width=\unitlength,page=1]{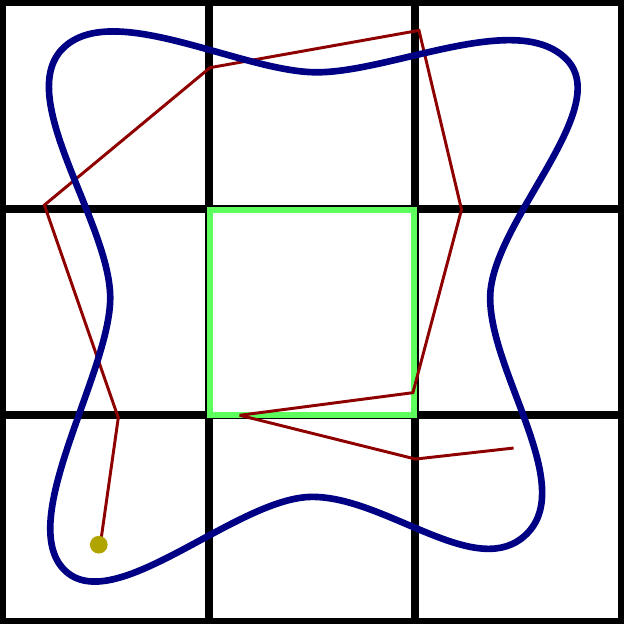}}%
    \put(0.46579081,0.39186305){\color[rgb]{0,0,0}\makebox(0,0)[lt]{\lineheight{1.25}\smash{\begin{tabular}[t]{l}$\mathbb{V}$\end{tabular}}}}%
    \put(0.44211997,0.72725364){\color[rgb]{0,0,0}\makebox(0,0)[lt]{\lineheight{1.25}\smash{\begin{tabular}[t]{l}$\mathbb{V}_{\eta}$\end{tabular}}}}%
    \put(0,0){\includegraphics[width=\unitlength,page=2]{fundamental_domain_neighborhood_int_bound.pdf}}%
    \put(0.18100625,0.15588315){\color[rgb]{0,0,0}\makebox(0,0)[lt]{\lineheight{1.25}\smash{\begin{tabular}[t]{l}$x_0$\end{tabular}}}}%
    \put(0.76783349,0.21365666){\color[rgb]{0,0,0}\makebox(0,0)[lt]{\lineheight{1.25}\smash{\begin{tabular}[t]{l}$x_9$\end{tabular}}}}%
  \end{picture}%
\endgroup%

		\caption[curve-bound]{An example of a curve whose existence is established
		in Lemma \ref{lem:neighborhood_curve_lem}.}
		\label{fig:neighborhood_curve_fig}
    \end{figure}
	This is the content of the following lemma.
	\begin{lemma}\label{lem:neighborhood_curve_lem}
		Let $g$ be a Riemannian metric on $\T^{3}$ and let $\V$ be a test domain
		for $\kappa(g,\eta)$, see Definition \ref{defn:covering_constant} above.
		Then there are $\kappa(g,\eta)$ copies of $\V$, say $\V^{\nu_{i}}$,
		generated by $\pi_{1}(\T^{3})$ such that the following statements are
		true.
			\begin{equation}
				\V_{\eta}\subset{}\bigcup{}^{m}_{i=1}\V^{\nu_{i}}.
			\end{equation}
			\begin{equation}
				\V_{\eta}\cap{}\V^{\nu_{i}}\neq\emptyset\forall i.
			\end{equation}
		Finally, for any $x_{0}$ and $x_{1}$ in $\V_{\eta}$ there exists a curve
		$c$ in $\bigcup{}\V^{i}$ and times $\{t_{j}\}_{1}^{l}$ with $l\leq\kappa(g,\eta)$ 
		satisfying the following properties:
			\begin{enumerate}
				\item $c(0)=x_{0}$;
				\item $c(t_{j})\in{}\V^{\nu_{i_{j}}}\cap{}\V^{\nu_{i_{j-1}}}$
				\item once $c$ leaves a domain $\V^{\nu_{i}}$, it
				      does not re-enter it.
			\end{enumerate}
		See Figure \ref{fig:neighborhood_curve_fig} 
	\end{lemma}
	\begin{proof}
		By definition, the fundamental domain $\V$ is path connected, and so we have
		that $\V_{\eta}$ is as well.
		Let $c_{0}:[0,1]\rightarrow{}\R^{3}$ be a curve in $\V_{\eta}$ connecting 
		$x_{0}$ to $x_{1}$.
		We will describe a process for modifying this curve to fit the criteria 
		laid out above.
		Let $i_{0}$ be the smallest index such that $\V^{i_{0}}=\V^{\nu_{i_{0}}}$
		contains $x_{0}=c_{0}(0)$.
		We define $t_0$ as follows:
			\begin{equation}
				t_{0}=\sup\bigl\{t:c_{0}(t)\in{}\V^{i_{0}}\bigr\}
			\end{equation}
		Since $c_{0}$ is continuous and $\V^{i_{0}}$ is closed, it follows that
			\begin{equation}
				c_{0}(t_{0})\in{}\V^{i_{0}}.
			\end{equation}
		Since $\V^{i_{0}}$ is path connected, we may replace 
		$\left.c_{0}\right|_{[0,t_{0}]}$ with a curve which lies entirely
		in $\V^{i_{0}}$.
		Let $c_{1}$ be the curve which results from this substitution,
		and let $\widetilde{t}_{1}$ be as follows
			\begin{equation}
				\widetilde{t}_{1}=\inf\bigl\{t:c_{1}(t)\in{}\V^{i};
				i\neq i_{0}\bigr\}
			\end{equation}
		Since there are only finitely many $\V^{i}$, we find a smallest $i_{1}$
		such that $c_{1}(\widetilde{t}_{1})\in{}\V^{i_{1}}$.
		By continuity, and the fact that all $\V^{\nu}$ are closed, we also have that
		$c_{1}(\widetilde{t}_{1})\in{}\V^{i_{0}}$, and actually $\widetilde{t}_{1}$
		is equal to $t_{0}$.
		Continuing as before, we let $t_{1}$ be defined as follows:
			\begin{equation}
				t_{1}=\sup\bigl\{t:c_{1}(t)\in{}\V^{i_{1}}\bigr\}
			\end{equation}	
		Since $\V^{i_{1}}$ is path connected, we may replace 
		$\left.c_{1}\right|_{[t_{0},t_{1}]}$ by a path contained entirely in
		$\V^{i_{1}}$.
		Let $c_{2}$ be the resulting curve, and continue in this manner.
		Since there are only finitely many $\V^{i}$, this process will terminate.
		The resulting curve has the desired properties.
	\end{proof}
	We are now in a position to obtain a bound on 
	$\sup{\V_{\eta}}\hat{u}-\inf{\V_{\eta}}\hat{u}$.
	\begin{lemma}\label{lem:sup-inf_nhbd_bound}
		Let $g$ be a Riemannian metric on $\T^{3}$ such that $\stsys_{2}(g)\geq\sigma$,
		let $\V$ be a test domain for $\kappa(g,\eta)$, and let 
		$u:(\T^{3},g)\rightarrow{}\Sp^{1}$ be a nontrivial harmonic map.
		Then, we have that 
			\begin{equation}
				\sup_{\V_{\eta}}\hat{u}-\inf_{\V_{\eta}}\hat{u}\leq
				\kappa(g,\eta)\sigma^{-1}|\T^{3}|^{\frac12}\|du\|_{L^{2}(g)}.
			\end{equation}
	\end{lemma}
	\begin{proof}
		By definition, the subset $\V_{\eta}$ is contained in $\kappa(g,\eta)$ copies
		$\V^{\nu_{i}}$ of $\V$.
		Let $x_{1}$ and $x_{0}$ be any two points in $\V_{\eta}$, and connect them
		by a curve $c\subset\bigcup{}_{i=1}^{\kappa(g,\eta)}$ such as in
		Lemma \ref{lem:neighborhood_curve_lem}.
		Let $\{t_{i_{j}}\}_{j=1}^{l}$ be the times such that
		$c(t_{i_{j}})\in{}\V^{\nu_{i_{j}}}\cap{}\V^{\nu_{i_{j-1}}}$, and let
		$x_{j}=c(t_{i_{j}})$.
		Then, we may calculate as follows:
			\begin{equation}
				\hat{u}(x_{1})-\hat{u}(x_{0})\leq
				\hat{u}(x_{1})-\hat{u}(x_{l})+
				\sum_{j=1}^{l}\hat{u}(x_{j})-\hat{u}(x_{j-1}).
			\end{equation}
		From Lemma \ref{lem:u_bounded_fund_domain} each element in the sum
		on the right hand side is bounded above by 
		$\sigma^{-1}|\T^{3}|_{g}^{\frac12}\|du\|_{L^{2}(g)}$.
		Therefore, we see that
			\begin{equation}
				\hat{u}(x_{1})-\hat{u}(x_{0})\leq\kappa(g,\eta)
				\sigma^{-1}|\T^{3}|_{g}^{\frac12}\|du\|_{L^{2}(g)}
			\end{equation}
	\end{proof}
	The above sup bound is the key to the integration by parts argument in the
	proof of the following lemma, as it
	allows us to avoid a H{\"o}lder like inequality, and so obtain control over higher $L^p$
	norms of non-trivial harmonic maps.
	\begin{lemma}\label{lem:L2_to_L3}
		Let $g$ be a Riemannian metric on $\T^{3}$, and suppose that
			\begin{equation}
				\min\{\stsys_{1}(g),\stsys_{2}(g)\}\geq\sigma.
			\end{equation}
		Given any nontrivial harmonic map, say $u:(\T^{3},g)\rightarrow{}\Sp^{1}$
		we have that
			\begin{align}
				\|du\|_{L^{3}(g)}\leq&1+\left(\left(1+\kappa(g,\eta)\sigma^{-1}
				|\T^{3}|^{\frac12}\|du\|_{L^{2}(g)}\right)
				\eta^{-1}\kappa(g,\eta)\|du\|^{2}_{L^{2}(g)}
				\right)^{\frac23}
				\\
				&+\left(\left(1+\kappa(g,\eta)\sigma^{-1}
				|\T^{3}|^{\frac12}\|du\|_{L^{2}(g)}\right)
				\kappa(g,\eta)\|du\|^{\frac32}_{L^{3}(g)}
				\|R^{-}_{g}\|^{\frac12}_{L^{2}(g)}
				\right)^{\frac23}.
			\end{align}
			
	\end{lemma}
	\begin{proof}
		Let $\V$ be a test domain for $\kappa(g,\eta)$, let $\V_{\eta}$ be the $\eta$
		neighborhood of $\V$, let $\widetilde{u}$ be a lift of $u$.
		From Lemma \ref{lem:sup-inf_nhbd_bound} we know that
			\begin{equation}
				\sup_{\V_{\eta}}\widetilde{u}-\inf_{\V_{\eta}}\widetilde{u}\leq
				\kappa(g,\eta)\sigma^{-1}|\T^{3}|^{\frac12}\|du\|_{L^{2}(g)}.
			\end{equation}
		Let $N=\left\lceil\inf_{\V_{\eta}}\widetilde{u}\right\rceil$, and let
		$\hat{u}=\widetilde{u}-N$.
		Then, $\hat{u}$ covers $u$, $d\hat{u}=d\widetilde{u}$, and
			\begin{equation}\label{eq:u-hat-00-bound}
				\|u\|_{L^{\infty{}}(\V_{\eta})}\leq
				\kappa(g,\eta)\sigma^{-1}|\T^{3}|^{\frac12}\|du\|_{L^{2}(g)}
				+1.
			\end{equation}
		Finally, let 
		$f:\R^{3}\rightarrow{}\R$ be a cutoff function such that
			\begin{enumerate}
				\item $\chi_{\V}\leq f\leq 1$;
				\item $\mathrm{Lip}(f)\leq\eta^{-1}$;
				\item $\mathrm{supp}(f)\subset{}\V_{\eta}$.
			\end{enumerate}
		With these elements we can calculate as follows:
			\begin{align}
				\int_{\T^{3}}|du|^{3}dV_{g}&\leq\int_{\V_{\eta}}
				g\bigl(f|d\hat{u}|d\hat{u},d\hat{u}\bigr)dV_{\pi^{*}g}
				\\
				&=-\int_{\V_{\eta}}\hat{u}\Bigl(g\bigl(df,|d\hat{u}|d\hat{u}\bigr)
					+fg\bigl(d|d\hat{u}|,d\hat{u}\bigr)\Bigr)dV_{\pi^{*}g},
			\end{align}
        where in the second line we integrated by parts and used the fact that
        $\hat{u}$ is an harmonic function.
		Taking the absolute value of both sides, and then using the Cauchy-Schwarz and
		Kato inequalities, we obtain
			\begin{equation}
				\|du\|^{3}_{L^{3}(g)}\leq \|\hat{u}\|_{L^{\infty}(\V_{\eta})}
				\left(\eta^{-1}\|d\hat{u}\|^2_{L^2(\V_{\eta})}+\int_{\V_{\eta}}
				|d\hat{u}||\nabla d\hat{u}|dV_{\pi^*g}\right).
			\end{equation}
		The right most term of the above can be estimated as follows
			\begin{align}
				\int_{\V_{\eta}}|d\hat{u}||\nabla d\hat{u}|dV_{\pi^*g}
				&=\int_{\V_{\eta}}|d\hat{u}|^{\frac32}
				\frac{|\nabla d\hat{u}|}{|d\hat{u}|^{\frac12}}dV_{\pi^*g}
				\\
				&\leq\|d\hat{u}\|^{\frac32}_{L^3(\V_{\eta})}
				\left(\int_{\V_{\eta}}\frac{|\nabla d\hat{u}|^2}{|d\hat{u}|}
				dV_{\pi^*g}\right)^{\frac12}.
			\end{align}
		Since $\V_{\eta}\subset{}\bigcup{}_{i=1}^{\kappa(g,\eta)}\V^{i}$, it follows from
        Proposition \ref{prop:lift_of_map_and_estimates} that
			\begin{equation}
				\|d\hat{u}\|^{2}_{L^{2}(\pi^{*}g,\V_{\eta})}\leq
				\kappa(g,\eta)\|du\|^{2}_{L^{2}(g)},
			\end{equation}
			\begin{equation}
				\|d\hat{u}\|^{3}_{L^{3}(\pi^{*}g,\V_{\eta})}\leq
				\kappa(g,\eta)\|du\|^{3}_{L^{3}(g)},
			\end{equation}
		and
			\begin{equation}\label{eq:fund_domain_stern_term}
				\int_{\V_{\eta}}\frac{|\nabla d\hat{u}|^{2}}{|d\hat{u}|}{}
				dV_{\pi^{*}g}\leq\kappa(g,\eta)\int_{\T^{3}}
				\frac{|\nabla du|^{2}}{|du|}{}dV_{g}.
			\end{equation}
		Therefore, we may apply Stern's inequality to the right hand side of Equation
        \eqref{eq:fund_domain_stern_term},
        see Corollary \ref{cor:stern_ineq},
		and H{\"o}lder's inequality to obtain
			\begin{equation}
				\int_{\V_{\eta}}|d\hat{u}||\nabla d\hat{u}|dV_{\pi^{*}g}\leq
				\kappa(g,\eta)\|du\|^{\frac32}_{L^{3}(g)}\|du\|^{\frac12}_{L^{2}(g)}
				\|R^{-}_{g}\|^{\frac12}_{L^{2}(g)}.
			\end{equation}
		Putting everything together gives us that
			\begin{align}
				\|du\|^{3}_{L^{3}(g)} \leq&\|\hat{u}\|_{L^{\infty{}}(\V_{\eta})}
				\eta^{-1}\kappa(g,\eta)\|du\|^{2}_{L^{2}(g)}
				\\
				&+\|\hat{u}\|_{L^{\infty{}}(\V_{\eta})}\kappa(g,\eta)
				\|du\|^{\frac32}_{L^{3}(g)}\|du\|^{\frac12}_{L^{2}(g)}
				\|R^{-}_{g}\|^{\frac12}_{L^{2}(g)}.
			\end{align}
		If $\|du\|_{L^{3}(g)}\geq1$, then we may divide both sides by 
		$\|du\|^{\frac32}_{L^{3}(g)}$, otherwise we already have a good bound.
		Therefore, we see that
			\begin{align}
				\|du\|^{\frac32}_{L^{3}(g)}\leq&1+
				\|\hat{u}\|_{L^{\infty{}}(\V_{\eta})}
				\eta^{-1}\kappa(g,\eta)\|du\|^{2}_{L^{2}(g)}
				\\
				&+\|\hat{u}\|_{L^{\infty{}}(\V_{\eta})}\kappa(g,\eta)
				\|du\|^{\frac12}_{L^{2}(g)}
				\|R^{-}_{g}\|^{\frac12}_{L^{2}(g)}.
			\end{align}
		Including the bound on $\|\hat{u}\|_{L^{\infty{}}(\V_{\eta})}$ from
        Equation \eqref{eq:u-hat-00-bound}, and then taking the $\frac23$ root gives the result.
	\end{proof}

	\section{Approximation by constant matrices}\label{sec:approx_const_matrices}
	
	The importance of obtaining $L^{3}$ estimates on non-trivial harmonic functions is
	twofold.
    The second reason is in some sense the heart of Lemma \ref{lem:well_approximating_sets_good_properties} below.
	In this section, we will explore the first reason for their importance.
	Consider two nontrivial harmonic functions from $(\T^{3},g)$ to $\Sp^{1}$, say
	$u^{j}$ and $u^{k}$. Then, we can control 
	$\|\nabla g\bigl(du^{j},du^{k}\bigr)\|_{L^{1}(g)}$ in terms of 
	$\max\{\|du^{j}\|_{L^{3}(g)},\|du^{k}\|_{L^{3}(g)}\}$ and 
	$\max\left\{\left\|\frac{|\nabla du^{j}|^{2}}{|du^{j}|}{}\right\|_{L^{1}(g)},
	\left\|\frac{|\nabla du^{k}|^{2}}{|du^{k}|}{}\right\|_{L^{1}(g)}\right\}$.
	This will be carried out for smooth maps in general, and will be applied to 
	harmonic maps in particular in Section \ref{sec:convergence}.
	We begin with the following lemma.
	\begin{lemma}\label{lem:good_approx_matrix}
		Let $g$ be a Riemannian metric on $\T^{3}$ such that $\ison_{1}(g)\geq\Lambda$,
		let $\U:(\T^{3},g)\rightarrow{}(\T^{3},h)$ be a smooth map, and let $u^{k}$ denote
		the components of $\U$.
		Then, we have that
			\begin{equation}\label{eq:metric_grad_bound}
				\left\|\nabla g\bigl(du^{j},du^{k}\bigr)\right\|_{L^{1}(g)}\leq
				2\sup_{jk}\|du^{j}\|^{\frac12}_{L^{3}(g)}
				\|du^{k}\|_{L^{3}(g)}
				\left(\int_{\T^{3}}
				\frac{|\nabla du^{j}|^{2}}{|du^{j}|}{}dV_{g}\right)^{\frac12}.
			\end{equation}
		Furthermore, letting $g^{jk}$ denote $g\bigl(du^{j},du^{k}\bigr)$, there
		exists a constant symmetric and non-negative matrix $a=a^{jk}$ such that
			\begin{equation}\label{eq:g_int_close_to_const}
				\|g^{jk}-a^{jk}\|_{L^{1}(g)}\leq
				2\Lambda^{-1}\sup_{jk}\|du^{j}\|^{\frac12}_{L^{3}(g)}
				\|du^{k}\|_{L^{3}(g)}
				\left(\int_{\T^{3}}
				\frac{|\nabla du^{j}|^{2}}{|du^{j}|}{}dV_{g}\right)^{\frac12}.
			\end{equation}
		Finally, we have that
			\begin{equation}\label{eq:const_matrix_bounded}
				\begin{split}
					\sup_{jk}|a^{jk}|\leq
					&2|\T^{3}|_{g}^{-1}\Lambda^{-1}\sup_{jk}
					\|du^{j}\|^{\frac12}_{L^{3}(g)}\|du^{k}\|_{L^{3}(g)}
					\left(\int_{\T^{3}}
						\frac{|\nabla du^{j}|^{2}}{|du^{j}|}{}
						dV_{g}
					\right)^{\frac12}
					\\
					&+|\T^3|_{g}^{-1}\sup_{jk}\|du^{j}\|_{L^{2}(g)}\|du^{k}\|_{L^{2}(g)}
					.
				\end{split}
			\end{equation}
	\end{lemma}
	\begin{proof}
		Once we establish \eqref{eq:metric_grad_bound}, we have that 
		\eqref{eq:g_int_close_to_const} follows from the definition of
		$\ison_{1}(g)$. Furthermore \eqref{eq:const_matrix_bounded} also follows by writing
		$|a^{jk}|=\frac{1}{|\T^{3}|_{g}}{}\int_{\T^{3}}|a^{jk}|dV_{g}$, then 
		adding and subtracting $g^{jk}$, applying the triangle inequality, and then
		using \eqref{eq:g_int_close_to_const}.

		Therefore, we need only establish \eqref{eq:metric_grad_bound}.
		Using the Cauchy-Schwarz inequality and rearranging terms, we estimate as follows
			\begin{equation}
				\begin{split}
					\int_{\T^{3}}|\nabla g\bigl(du^{j},du^{k}\bigr)|dV_{g}\leq&
					\int_{\T^{3}}|du^{k}||du^{j}|^{\frac12}
					\frac{|\nabla du^{j}|}{|du^{j}|^{\frac12}}{}dV_{g}
					\\
					&+\int_{\T^{3}}|du^{j}||du^{k}|^{\frac12}
					\frac{|\nabla du^{k}|}{|du^{k}|^{\frac12}}{}dV_{g}.
				\end{split}
			\end{equation}
		Both integrals may be estimated using H{\"o}lder's inequality for three
		terms with exponents $3,6$, and $2$, respectively.
		This gives \eqref{eq:metric_grad_bound}, and so the other results as well.
	\end{proof}
	The above result can be used to show that on a subset of $\T^{3}$ the sup norm of
	$|g^{jk}-a^{jk}|$ is controlled.
	\begin{corollary}\label{cor:first_approximating_set}
		Let $g$ be a metric on $\T^{3}$ such that $\ison_{1}(g)\geq\Lambda$, let
		$\U:\T^{3}\rightarrow{}\T^{3}$ be a smooth map, let $u^{j}$ be the components
		of $\U$, let $g^{jk}$ denote $g\bigl(du^{j},du^{k}\bigr)$, and let
		$a^{jk}$ be as in Lemma \ref{lem:good_approx_matrix}.
		Finally, let us denote by $E^{1}(g,\tau)$ the following set
			\begin{equation}
				E^{1}(g,\tau)=\left\{x:\sum_{jk}|g^{jk} - a^{jk}|<\tau\right\}.
			\end{equation}
		Set $\tau$ to be as follows
			\begin{equation}
				\tau=
				\left(\frac{36}{|\T^{3}|_{g}\Lambda}{}
				\sup_{jk}\|du^{j}\|^{\frac12}_{L^{3}(g)}
				\|du^{k}\|_{L^{3}(g)}
				\left(\int_{\T^{3}}
				\frac{|\nabla du^{j}|^{2}}{|du^{j}|}{}dV_{g}\right)^{\frac12}
                \right)^{\frac12}.
			\end{equation}
        We have that if $\tau<1$, then
            \begin{equation}
                |E^1(g,\tau)|_{g}\geq\frac12|\T^3|_{g}.
            \end{equation}
	\end{corollary}
	\begin{proof}
		The proof is shorter than the statement.
		After summing \eqref{eq:g_int_close_to_const} over the indices
		$j,k$ we can apply Chebyshev's inequality to the result to obtain
        \begin{align}
            |E^1(g,\tau)^c|_{g}&\leq\frac{1}{\tau}\frac{18}{|\T^{3}|_{g}\Lambda}{}
				\sup_{jk}\|du^{j}\|^{\frac12}_{L^{3}(g)}
				\|du^{k}\|_{L^{3}(g)}
				\left(\int_{\T^{3}}
				\frac{|\nabla du^{j}|^{2}}{|du^{j}|}{}dV_{g}\right)^{\frac12}
            \\
            &=\frac{\tau}{2}|\T^3|_{g}.
        \end{align}
        Thus, if $\tau<\tfrac12$, we get the result.
	\end{proof}
	Once we have a good subset to begin with, we can use it to show that there exists
	an open connected submanifold with smooth boundary on which $|g^{jk}-a^{jk}|$ is controlled.
	\begin{lemma}\label{lem:well_approximating_set}
		Let $g$ be a metric on $\T^{3}$ be such that $\ison_{1}(g)\geq\Lambda$,
		let $\U:\T^{3}\rightarrow{}\T^{3}$ be a smooth map, let $u^{j}$ be the 
		components of $\U$, let $g^{jk}$ denote $g\bigl(du^{j},du^{k}\bigr)$,
		and let $a^{jk}$ be as in Lemma \ref{lem:good_approx_matrix}.
		Next, let us denote by $E^{2}(g,\tau)$ the following set
			\begin{equation}
				E^{2}(g,\tau)=\left\{
				x:\sum_{jk}|g^{jk}-a^{jk}|^{2}<\tau^2
				\right\}	
			\end{equation}
		Finally let $\tau$ be given by
			\begin{equation}
				\tau=\left(
				\frac{36}{|\T^{3}|_{g}\Lambda}{}
				\sup_{jk}\|du^{j}\|^{\frac12}_{L^{3}(g)}
				\|du^{k}\|_{L^{3}(g)}
				\left(\int_{\T^{3}}
				\frac{|\nabla du^{j}|^{2}}{|du^{j}|}{}dV_{g}\right)^{\frac12}
				\right)^{\frac12}.
			\end{equation}
		If $\tau\leq\frac18$, then there exists an open connected submanifold
		$\Omega(g,\tau)$ with smooth boundary and the following properties:
			\begin{enumerate}
				\item $\Omega(g,\tau)\subset{}
					\{x:\sup_{jk}|g^{jk}-a^{jk}|\leq2\tau\}$;
				\item $|\Omega(g,\tau)|\geq\frac12|\T^{3}|_{g}$;
				\item $|\partial{}\Omega|\leq2|\T^{3}|_{g}\Lambda\tau$;
				\item $|\Omega(g,\tau)^{c}|\leq2|\T^{3}|_{g}\tau$.
			\end{enumerate}
	\end{lemma}
	\begin{proof}
		From the coarea formula, we have that
			\begin{equation}
				\int_{\tau^{2}}^{4\tau^{2}}|\partial{}E^{2}(g,\sqrt{s})|ds
				\leq
				2\int_{E^{2}(g,2\tau)}\sum_{jk}|g^{jk}-a^{jk}||\nabla g^{jk}|dV_{g}.
			\end{equation}
		Since we are working inside of $E^{2}(g,2\tau)$, it follows that for all $j,k$
		we have that $|g^{jk}-a^{jk}|\leq2\tau$.
		Therefore, the righthand side of the above is bounded by
			\begin{equation}
				36\tau\sup_{jk}\int_{\T^{3}}|\nabla g^{jk}|dV_{g}.
			\end{equation}
		Therefore, we may apply \eqref{eq:metric_grad_bound} to see that
			\begin{equation}
				\int_{\tau^2}^{4\tau^{2}}|\partial{}E^{2}(g,\sqrt{s})|ds
				\leq
				2\tau|\T^{3}|_{g}\Lambda\tau^{2}=2|\T^{3}|\Lambda\tau^{3}.
			\end{equation}
		Using Chebyshev's inequality on the above equation shows us that
			\begin{equation}
				\left|\left\{s:
					|\partial{}E^{2}(g,\sqrt{s}) | \leq2|\T^{3}|_{g}
					\Lambda\tau
				\right\}\cap{}[\tau^{2},4\tau^2]\right|\geq2\tau^{2}.
			\end{equation}
		Since $\sum_{jk}|g^{jk}-a^{jk}|^{2}$ is smooth, we may apply Sard's Lemma
		to conclude that almost every value is regular.
		In particular, if the image of this function contains all of 
		$[\tau^2,4\tau^2]$, then we may find a regular value 
		$t_{0}\in{}[\tau^2,4\tau^{2}]$ such that
			\begin{equation}
				|\partial{}E^{2}(g,\sqrt{t_{0}})|\leq2|\T^{3}|_{g}\Lambda\tau.
			\end{equation}
		So, suppose that the image of $\sum_{jk}|g^{jk}-a^{jk}|$ does not contain
		$[\tau^2,4\tau^2]$.
		Since the function is continuous, and so has connected image, 
		there is an $\varepsilon{}>0$ such that the image lies either in
		$[0,4\tau^{2}-\varepsilon{}]$ or in $[4\tau^{2},\infty{})$.
		We now observe that 
		$\sum_{jk}|g^{jk}-a^{jk}|^{2}\leq\left(\sum_{jk}|g^{jk}-a^{jk}|\right)^{2}$,
		and so
			\begin{equation}
				E^{1}(g,\tau)\subset{}E^{2}(g,\tau)	
			\end{equation}
		Since $\tau\leq\tfrac18$ by assumption, we know from Corollary \ref{cor:first_approximating_set} that
		$|E^{1}(g,\tau)|\geq\frac12|\T^{3}|_{g}$, and so we must have that
		$|E^{2}(g,\tau)|\geq\frac12|\T^{3}|_{g}$.
		In this case, we must have that the image of $\sum_{jk}|g^{jk}-a^{jk}|^{2}$
		must lie in $[0,4\tau^{2}-\varepsilon{}]$.
		Then, we see that $E^{2}(g,2\tau)=\T^{3}$, and the result is clear.

		So, let us continue by supposing that the image of 
		$\sum_{jk}|g^{jk}-a^{jk}|^2$ does indeed contain $[\tau^{2},4\tau^{2}]$, and
		use Sard's Lemma to find a $t_{0}$ in $[\tau^{2},4\tau^{2}]$ such that
		$\partial{}E^{2}(g,\sqrt{t_{0}}$ is smooth and satisfies
			\begin{equation}
				|\partial{}E^{2}(g,\sqrt{t_{0}})|\leq2|\T^{3}|_{g}\Lambda\tau.	
			\end{equation}
		$E^{2}(g,\sqrt{t_{0}})$ can be decomposed into a union of disjoint
		connected components, say $G_{i}$.
		Since
        \begin{equation}
            \nabla\sum_{jk}|g^{jk}-a^{jk}|^{2}
        \end{equation}
        is nonzero and outward pointing
		everywhere on $\partial{}E^{2}(g,\sqrt{t_{0}})$, it follows that
			\begin{equation}
				\partial{}E^{2}(g,\sqrt{t_{0}})=
				\bigsqcup_{i}\partial{}G_{i}.
			\end{equation}
		From this, we immediately see that if there exists an $i_{0}$ such that
		$|G_{i_{0}}|\geq\frac12|\T^{3}|_{g}$, then we may take $\Omega(g,\tau)$
		to be $G_{i_{0}}$.

		So, suppose that no such index exists.
		Then, we have that
			\begin{equation}
				\frac12|\T^{3}|_{g}\leq|E(g,\sqrt{t_{0}})|=
				\sum_{i}|G_{i}|\leq
				\Lambda^{-1}\sum_{i}|\partial{}G_{i}|=
				\Lambda^{-1}|\partial{}E^{2}(g,\sqrt{t_{0}})|.
			\end{equation}
		However, by the work done above, we know that the last term on the right
		has the bound
			\begin{equation}
				\Lambda^{-1}|\partial{}E^{2}(g,\sqrt{t_{0}})|\leq
				2|\T^{3}|\tau.
			\end{equation}
		Since we assumed that $\tau\leq\frac18$, this last inequality gives us a
		contradiction.
		Therefore, there does indeed exist some $i_{0}$ such that 
		$|G_{i_{0}}|\geq\frac12|\T^{3}|_{g}$, and we let $\Omega(g,\tau)=G_{i_{0}}$.
	\end{proof}

	\section{Convergence}\label{sec:convergence}
	
	In this section we will combine the estimates for smooth functions coming from control
	on $\ison_{1}(g)$ with the estimates on harmonic maps coming from the $L^{2}$ bounds
	obtained in terms of stable systoles, and $L^{3}$ bounds which result from
	control on the covering constant of $g$, see Definition \ref{defn:covering_constant}.
	This combined control will lead to convergence in the sense of Dong-Song for
	certain sequences of Riemannian metrics whose negative part of their scalar curvature
	tends to zero.

	We begin by defining the family of metrics on $\T^{3}$ for which the convergence
	result will hold.
	\begin{definition}\label{defn:family_of_metrics}
		Let $V,R,\Lambda,\sigma,\eta,M>0$, and define
		$\F=\F(V,R,\Lambda,\sigma,\eta,M)$ to be the family of Riemannian metrics
		on $\T^{3}$ such that for all $g\in{}\F$ we have that
			\begin{enumerate}
				\item $|\T^{3}|_{g}\leq V$;
				\item $\|R^{-}_{g}\|_{L^{2}(g)}\leq R$;
				\item $\ison_{1}(g)\geq\Lambda$;
				\item $\min\bigl\{\stsys_{1}(g),\stsys_{2}(g)\bigr\}\geq\sigma$;
				\item $\kappa(g,\eta)\leq M$.
			\end{enumerate}
	\end{definition}
	The results of the previous section imply strong controls on metrics in $\F$.
	In order to make this clear, we will summarize the results obtained so far
	as they apply to $\F$.
	\begin{notation}
		In order to avoid ever expanding equations and terms, we shall denote by
		$B$ any constant which depends only on $V,R,\Lambda,\sigma,\eta$ and $M$.
		It may be that from line to line $B$ will change, increasing to be as large as
        necessary. This will only happen a finite number of times.
	\end{notation}	
	Let us begin with a volume lower bound and the existence of well controlled harmonic maps
    for the metrics in $\F$. 
	\begin{proposition}\label{prop:vol_lower_bound-harmonic_maps}
		There exists a constant $B>0$ depending only on $V,R,\Lambda,\sigma,\eta$ and
		$M$ such that for any $g\in{}\F$ we have
			\begin{equation}
				|\T^{3}|_{g}\geq B^{-1}>0.
			\end{equation}
		Furthermore, if $\U_{g}:(\T^{3},g)\rightarrow{}(\T^{3},h)$ is the harmonic
		map produced in Corollary \ref{cor:L2-bounded_deg-1-harmonic_map} with
		components $u^{j}$, then for all $j$
		we have that 
			\begin{align}
				&\|du^{j}\|_{L^{2}(g)}\leq B;
				\\
				&\int_{\T^{3}}
				\frac{|\nabla du^{j}|^{2}}{|du^{j}|}{}dV_{g}\leq 
				B\|R^{-}\|_{L^{2}(g)};
				\\
				&\|du^{j}\|_{L^{3}(g)}\leq B.
			\end{align}
	\end{proposition}
	\begin{proof}
		Since $\stsys_{1}(g)\leq \mathrm{sys}_{1}(g)$, the first result follows from
		Gromov's systolic inequality.
		The next inequality follows from applying 
		Corollary \ref{cor:L2-bounded_deg-1-harmonic_map} to metrics in $\F$.
		Then, we have from Stern's inequality in the form of Corollary
		\ref{cor:stern_ineq} that
			\begin{equation}
				\int_{\T^{3}}
				\frac{|\nabla du^{j}|^{2}}{|du^{j}|}{}dV_{g}\leq
				\int_{\T^{3}}|du^{j}|R^{-}_{g}dV_{g}\leq
				\|du^{j}\|_{L^{2}(g)}\|R^{-}_{g}\|_{L^{2}(g)}.
			\end{equation}
        \noindent
		The final bound on the $L^3$ norms of the coordinate functions
        follows from substituting the above into
		Lemma \ref{lem:L2_to_L3}.
	\end{proof}
	Next, we see that for any $g\in{}\F$ there is a constant matrix which approximates
	$g$ in an integral sense, and the quality of the approximation depends on
	$\|R^{-}_{g}\|_{L^{2}(g)}$.
	\begin{proposition}\label{prop:const_matrix_approx}
		Fix $V,R,\Lambda,\sigma,\eta$, and $M$ greater than zero, and for every
		$g\in{}\F$ let $\U=\U_{g}$ be the degree 1 harmonic map
		$\U_{g}:(\T^{3},g)\rightarrow{}(\T^{3},h)$ given in Corollary
		\ref{cor:L2-bounded_deg-1-harmonic_map}, let $u^{j}$ denote its
		components,
		and let $g^{jk}=g\bigl(du^{j},du^{k}\bigr)$.
		Then there exists a constant $B$ depending only on
		$V,R,\Lambda,\sigma,\eta$ and $M$ such that
		for every $g\in{}\F$ there exists a symmetric and non-negative
		matrix $a^{jk}$ with the following properties:
			\begin{equation}
				\int_{\T^{3}}|g^{jk}-a^{jk}|dV_{g}\leq
				B\|R^{-}_{g}\|^{\frac12}_{L^{2}(g)}
			\end{equation}
		and
			\begin{equation}
				\sup_{jk}|a^{jk}|\leq
				B.
			\end{equation}
	\end{proposition}
	\begin{proof}
		This result will follow from Lemma \ref{lem:good_approx_matrix} if we can
		control $\sup_{j}\|du^{j}\|_{L^{3}(g)}$ and
		$\sup_{j}\int_{\T^{3}}\frac{|du^{j}|^{2}}{|du^{j}|}{}dV_{g}$ uniformly
		in terms of $\|R^{-}_{g}\|_{L^{2}(g)}$ and some constant $B$.
		Luckily, this is the content of Proposition 
		\ref{prop:vol_lower_bound-harmonic_maps}.
	\end{proof}
	In fact, as we show in the following proposition, for any metric $g\in{}\F$ whose
	negative part of the scalar curvature is small enough we can find a large connected open
	sub-manifold with smooth and small boundary on which $g^{jk}$ is uniformly well
	approximated by a constant matrix.
    Before diving in, let us recall an elementary estimate on determinants, which
    will be a useful tool in the proof of Lemma \ref{lem:well_approximating_sets_good_properties}
    below.
	\begin{proposition}\label{prop:det_estimate}
		Let $a$ and $b$ be any two $n\times{}n$ matrices.
		Then, there exists a constant $C=C(n)$, depending only on $n$, such that
			\begin{equation}
				|\mathrm{det}(a)-\mathrm{det}(b)|\leq
				C(\|a\|+\|b\|)^{n-1}\|a-b\|.
			\end{equation}
	\end{proposition}
	\begin{proof}
		Let $c(t)=(1-t)a+tb$, and use the cofactor expansion to calculate
			\begin{equation}
				\frac{d}{dt}{}\mathrm{det}\bigl(c(t)\bigr)=
				\sum_{j=1}^{n}\Bigl((a-b)_{ij}\mathrm{Cof}\bigl(c(t)\bigr)_{ij}+
				c(t)_{ij}\frac{d}{dt}{}\mathrm{Cof}\bigl(c(t)\bigr)_{ij}\Bigr).
			\end{equation}
		Therefore, we may use an inductive argument to see that
			\begin{equation}
				\left|\frac{d}{dt}{}\mathrm{det}\bigl(c(t)\bigr)\right|\leq
				C(n)(\|a\|+\|b\|)^{n-1}\|a-b\|.
			\end{equation}
	\end{proof}
	We are now in the position to prove the following Lemma.
	\begin{lemma}\label{lem:well_approximating_sets_good_properties}
		Fix $V,R,\Lambda,\sigma,\eta$, and $M$ greater than zero, and for every
		$g\in{}\F$ let $\U=\U_{g}$ be the degree 1 harmonic map
		$\U:(\T^{3},g)\rightarrow{}(\T^{3},h)$ given in Corollary
		\ref{cor:L2-bounded_deg-1-harmonic_map}, let $u^{j}$ denote its
		components,
		and let $g^{jk}=g\bigl(du^{j},du^{k}\bigr)$.
		For every $g\in{}\F$ let $a^{jk}$ be the constant symmetric non-negative
		matrix given in Proposition \ref{prop:const_matrix_approx}.
		Then, there is a $B$ such that for any $g\in{}\F$ with 
			\begin{equation}
				\|R^{-}_{g}\|_{L^{2}(g)}\leq\frac{1}{B}{}
			\end{equation}
		there is a connected open submanifold $\Omega(g)$ with smooth boundary,
		and which satisfies the following properties:
			\begin{align}
				&\Omega(g)\subset{}
					\Bigl\{x:\sup_{jk}|g^{jk}-a^{jk}|\leq 
					B\|R^{-}_{g}\|^{\frac14}_{L^{2}(g)}\Bigr\};
				\label{eq:const_approximation}
				\\
				&|\Omega(g)^{c}|\leq B\|R^{-}_{g}|_{L^{2}(g)}^{\frac14};
				\label{eq:small_compliment}
				\\
				&|\partial{}\Omega|\leq B\|R^{-}_{g}\|^{\frac14}_{L^{2}(g)};
				\label{eq:small_boundary}
				\\
				&\|d\U\|^{3}_{L^{3}\left(g;\Omega(g)^{c}\right)}\leq B\|R^{-}_{g}\|
				^{\frac1{12}}_{L^{2}(g)};
				\label{eq:bound_on_L3_compl_Omega}
				\\
				&1+B\|R^{-}_{g}\|^{\frac1{12}}_{L^{2}(g)}\geq
				\int_{\Omega(g)}\mathrm{det}(d\U)dV_{g}\geq
				1-B\|R^{-}_{g}\|^{\frac1{12}}_{L^{2}(g)};
				\label{eq:up_lower_int_bound_detU_Omega}
				\\
				&\int_{\Omega(g)}|\mathrm{det}(d\U)|dV_{g}= 
				\int_{\Omega(g)}\mathrm{det}(d\U)dV_{g};
				\label{eq:detdU_has_sign_Omega}
				\\
				&|\U\left(\Omega(g)^{c}\right)|_{h}\leq B\|R^{-}_{g}\|
				^{\frac1{12}}_{L^2(g)};
				\label{eq:U-image_compl_small_in_h}
				\\
				&\left.\mathrm{det}\left(g^{jk}\right)\right|_{\Omega(g)}\geq
					\frac1{|\Omega(g)|^{2}}\left(1-B\|R^{-}_{g}\|^{\frac1{12}}
						\right)^{2}-B\|R^{-}_{g}\|^{\frac14}_{L^{2}(g)};
				\label{eq:det_g_bounded_below_Omega}
				\\
				&\mathrm{det}\left(a^{jk}\right)\geq
					\frac1{|\Omega(g)|^{2}}\left(1-B\|R^{-}_{g}\|^{\frac1{12}}
						\right)^{2}-B\|R^{-}_{g}\|^{\frac14}_{L^{2}(g)}.
				\label{eq:det_a_bounded_below}
			\end{align}
	\end{lemma}
	\begin{proof}
        We begin by proving the first three results.
        In order to do this, we will apply Lemma \ref{lem:well_approximating_set} to
	    metrics in $\F$.
	    Let us recall that this lemma was stated in terms of a parameter $\tau$ given as follows:
		\begin{equation}
			\tau=\left(
			\frac{36}{|\T^{3}|_{g}\Lambda}{}
			\sup_{jk}\|du^{j}\|^{\frac12}_{L^{3}(g)}
			\|du^{k}\|_{L^{3}(g)}
			\left(\int_{\T^{3}}
			\frac{|\nabla du^{j}|^{2}}{|du^{j}|}{}dV_{g}\right)^{\frac12}
			\right)^{\frac12}.
		\end{equation}
        We know from Proposition \ref{prop:vol_lower_bound-harmonic_maps} that for $\tau$
        as above, we have that
	    $\tau\leq B\|R^{-}_{g}\|_{L^{2}(g)}^{\frac14}$.
        Thus, the requirement that
	    $\tau\leq\frac18$ for Lemma \ref{lem:well_approximating_set} follows from
    	assuming that $\|R^{-}_{g}\|_{L^{2}(g)}\leq
    	\tfrac{1}{B}{}$, where $B$ is chosen to be sufficiently large.
        The first three results now follow from Lemma \ref{lem:well_approximating_set} and
		the expression of $\tau$ in terms of $B\|R^{-}_{g}\|_{L^{2}(g)}^{\frac14}$.

		We will begin by establishing an estimate for 
		$\|du^{j}\|_{L^{3}\left(g;\Omega(g)^{c}\right)}$.
		This will involve another integration by parts argument, very similar to
		the one in the proof of Lemma \ref{lem:L2_to_L3}.
		Let $\V$ be a test domain for $\kappa(g,\eta)$ as in Definition 
		\ref{defn:covering_constant}.
		Using the $L^{2}$ bound in Proposition \ref{prop:vol_lower_bound-harmonic_maps}
		together with Lemma \ref{lem:sup-inf_nhbd_bound}, for any $j$ we find a lift of
		$u^{j}$, say $\widetilde{u}^{j}$, such that
			\begin{equation}
				\sup_{\V_{\eta}}\widetilde{u}^{j}-\inf_{\V_{\eta}}\widetilde{u}^{j}
				\leq B.
			\end{equation}
		Let $\hat{u}^{j}=\widetilde{u}^{j}-
		\left\lceil\inf_{\V_{\eta}}\widetilde{u}^{j}\right\rceil$.
		Then, $\hat{u}^{j}$ is a lift of $u^{j}$ such that
			\begin{equation}
				\|\hat{u}^{j}\|_{L^{\infty{}}(\V_{\eta})}\leq B. 
			\end{equation}
		Finally, let 
		$f:\R^{3}\rightarrow{}\R$ be a cutoff function such that
			\begin{enumerate}
				\item $\chi_{\V}\leq f\leq 1$;
				\item $\mathrm{Lip}(f)\leq\eta^{-1}$;
				\item $\mathrm{supp}(f)\subset{}\V_{\eta}$.
			\end{enumerate}
		Then, letting $\vec{n}$ denote the outwards unit normal to
        $\partial\pi^{-1}\Omega(g)^c$,
        we may use integration by parts to calculate as follows:
			\begin{align}
				\int_{\Omega(g)^{c}}|du^{j}|^{3}dV_{g}&\leq
				\int_{\pi^{-1}(\Omega(g)^{c})}\pi^{*}g
				\bigl(f|d\hat{u}^{j}|d\hat{u}^{j},d\hat{u}^{j}\bigr)dV_{\pi^{*}g}
				\\
				&=\int_{\partial{}\pi^{-1}\Omega(g)^{c}}\hat{u}^{j}f
				|d\hat{u}^{j}|d\hat{u}^{j}(\vec{n})dA_{\pi^{*}g}
				\\
				&\quad-\int_{\pi^{-1}\Omega(g)^{c}}\hat{u}^{j}
				\pi^{*}g\bigl(df,|d\hat{u}^{j}|d\hat{u}^{j}\bigr)dV_{\pi^{*}(g)}
				\\
				&\quad-\int_{\pi^{-1}\Omega(g)^{c}}\hat{u}^{j}f
				\pi^{*}g\bigr(d|d\hat{u}^{j}|,d\hat{u}^{j}\bigr)dV_{\pi^{*}g}.
			\end{align}
		We now observe that from the properties of the covering map and $\Omega(g)$,
		we have that
			\begin{equation}
				\partial{}\pi^{-1}\Omega(g)^{c}=\pi^{-1}\partial{}\Omega(g),
			\end{equation}
        where we are looking at the full preimage.
		Therefore, taking the absolute value of these terms gives us that
			\begin{align}
				\|du^{j}\|^{3}_{L^{3}\bigl(g;\Omega(g)^{c}\bigr)}\leq&
				\|\hat{u}^{j}\|_{L^{\infty{}}(\V_{\eta})}
				\int_{\V_{\eta}\cap{}\pi^{-1}\partial{}\Omega(g)^{c}}
				|d\hat{u}^{j}|^{2}dA_{\pi^{*}g}
				\\
				&+\eta^{-1}\|\hat{u}^{j}\|_{L^{\infty{}}(\V_{\eta})}
				\int_{\V_{\eta}\cap{}\pi^{-1}\Omega(g)^{c}}|d\hat{u}^{j}|^{2}
				dV_{\pi^{*}g}
				\\
				&+\|\hat{u}^{j}\|_{L^{\infty{}}(\V_{\eta})}
				\int_{\V_{\eta}\cap{}\pi^{-1}\Omega(g)^{c}}
				|d\hat{u}^{j}||\nabla d\hat{u}^{j}|dV_{\pi^{*}(g)}.
			\end{align}
		It now follows from the definition of
		$\kappa(g,\eta)$ and the fact that $\hat{u}^{j}$ covers $u^{j}$ that we 
		have the following inequalities:
			\begin{align}
				&\int_{\V_{\eta}\cap{}\pi^{-1}\partial{}\Omega(g)}
				|d\hat{u}^{j}|^{2}dA_{\pi^{*}g}
				\leq\kappa(g,\eta)\int_{\partial{}\Omega(g)}
				|du^{j}|^{2}dA_{g};
				\\
				&\int_{\V_{\eta}\cap{}\pi^{-1}\Omega(g)^{c}}
				|d\hat{u}^{j}|^{2}dV_{\pi^{*}g}\leq
				\kappa(g,\eta)\int_{\Omega(g)^{c}}|du^{j}|^{2}dV_{g};
				\\
				&\int_{\V_{\eta}\cap{}\pi^{-1}\Omega(g)^{c}}
				|d\hat{u}^{j}||\nabla d\hat{u}^{j}|dV_{\pi^{*}g}\leq
				\kappa(g,\eta)\int_{\T^{3}}
				|du^{j}||\nabla du^{j}|dV_{g}.
			\end{align}
		We now recall that $\kappa(g,\eta)\leq M$ for all $g\in{}\F$ by assumption.
		
		Since $|g^{jj}-a^{jj}|$ is small, see the first property of $\Omega(g)$ in Lemma 
		\ref{lem:well_approximating_sets_good_properties}, and $a^{jj}$ is bounded, see the second conclusion of
		Proposition \ref{prop:const_matrix_approx}, we have that $|du^{j}|\leq B$
        on $\partial\Omega(g)^c$.
		As the area of $\partial\Omega(g)$ is bounded in terms of $\|R^{-}_{g}\|^{\frac14}_{L^2(g)}$,
        from the second property of $\Omega(g)$ in Lemma
		\ref{lem:well_approximating_sets_good_properties}, we have that
			\begin{equation}
				\int_{\partial{}\Omega(g)}|du^{j}|^{2}dA_{g}\leq
				|\partial{}\Omega(g)|B\leq B\|R^{-}_{g}\|^{\frac14}_{L^{2}(g)}.
			\end{equation}
		Next, we may apply H{\"o}lder's inequality to obtain
			\begin{equation}
				\int_{\Omega(g)^{c}}|du^{j}|^{2}dV_{g}\leq
				|\Omega(g)^{c}|^{\frac13}\|du^{j}\|^{2}_{L^{3}(g)}.
			\end{equation}
		So, from the third property of $\Omega(g)$ appearing in Lemma
		  \ref{lem:well_approximating_sets_good_properties} and Proposition
		\ref{prop:vol_lower_bound-harmonic_maps}, we see that we have the bound
			\begin{equation}
				\int_{\Omega(g)^{c}}|du^{j}|^{2}dV_{g}\leq
				B\|R^{-}_{g}\|^{\frac{1}{12}}.	
			\end{equation}
		Finally, rewriting $|du^{j}||\nabla du^{j}|$ as $|du^{j}|^{\frac32}
		|\frac{|\nabla du^{j}|}{|du^{j}|^{\frac12}}{}|$ and
		using Proposition \ref{prop:vol_lower_bound-harmonic_maps}
		we have
			\begin{equation}
				\int_{\T^{3}}|du^{j}||\nabla du^{j}|dV_{g}\leq
				B\|R^{-}_{g}\|_{L^{2}(g)}^{\frac12}.
			\end{equation}
		Thus, putting everything together results in the following inequality:
			\begin{align}
				\|du^{j}\|^{3}_{L^{3}\bigl(g;\Omega(g)^{c}\bigr)}&\leq
				\kappa(g,\eta)B\left(\|R^{-}_{g}\|^{\frac1{12}}_{L^{2}(g)}+
				2\|R^{-}_{g}\|^{\frac12}_{L^{2}(g)}\right)
				\\
				&\leq B\|R^{-}_{g}\|^{\frac1{12}}_{L^{2}(g)},
			\end{align}
		where we got the last inequality by absorbing the bound $\kappa(g,\eta)\leq M$
		and the fact that  $\|R^{-}\|^{\frac{5}{12}}_{L^{2}(g)}$ is less than 
		$R^{\frac{5}{12}}$ into the constant $B$. This gives us
		\eqref{eq:bound_on_L3_compl_Omega}.

		We will now use \eqref{eq:bound_on_L3_compl_Omega} to estimate
		$\int_{\Omega(g)}|\mathrm{det}\bigl(d\U\bigr)|dV_{g}$.
		To do this, recall that because $\U$ is a degree 1 map, we have that
			\begin{equation}
				\int_{\T^{3}}\mathrm{det}\bigl(d\U\bigr)dV_{g}=1.
			\end{equation}
		So, we see that
			\begin{equation}
				\int_{\Omega(g)}\mathrm{det}\bigl(d\U\bigr)dV_{g}=
				1-\int_{\Omega(g)^{c}}\mathrm{det}(d\U)dV_{g}.
			\end{equation}
		Furthermore, we have the estimate $|\mathrm{det}\bigl(d\U\bigr)|\leq|d\U|^{3}$,
		so
			\begin{equation}\label{eq:bound_on_det_compl_Omega}
				\int_{\Omega(g)^{c}}|\mathrm{det}\bigl(d\U\bigr)|dV_{g}\leq
				\|d\U\|^{3}_{L^{3}(g;\Omega(g)^{c})}\leq
				B\|R^{-}_{g}\|^{\frac1{12}}_{L^{2}(g)}.
			\end{equation}
		Combined with the above, we get
			\begin{equation}
				1+B\|R^{-}_{g}\|^{\frac1{12}}_{L^{2}(g)}\geq
				\int_{\Omega(g)}\mathrm{det}(d\U)dV_{g}\geq
				1-B\|R^{-}_{g}\|^{\frac1{12}}_{L^{2}(g)},
			\end{equation}
		which is \eqref{eq:up_lower_int_bound_detU_Omega}.

		Next, we can use \eqref{eq:bound_on_det_compl_Omega} and the area formula to see that
			\begin{equation}
				\int_{\U\left(\Omega(g)^{c}\right)}|\U^{-1}\{y\}|dV_{h}=
				\int_{\Omega(g)^{c}}|\mathrm{det}\bigl(d\U\bigr)|dV_{g}
				\leq
				B\|R^{-}_{g}\|^{\frac1{12}}_{L^{2}(g)}.
			\end{equation}
		Thus, it follows that
			\begin{equation}
				|\U\bigl(\Omega(g)^{c}\bigr)|_{h}\leq
				\int_{\U\bigl(\Omega(g)^{c}\bigr)}|\U^{-1}\{y\}|dV_{h}\leq
				B\|R^{-}_{g}\|^{\frac1{12}}_{L^{2}(g)}.
			\end{equation}
		This gives us \eqref{eq:U-image_compl_small_in_h}.

		Since $\mathrm{det}(g^{jk})=\mathrm{det}(d\U)^{2}$, we have that
			\begin{equation}
				\int_{\Omega(g)}\mathrm{det}(g^{jk})dV_{g}\geq
				\frac{1}{|\Omega(g)|}{}\left(\int_{\Omega(g)}|
					\mathrm{det}(d\U)|dV_{g}\right)^{2}.
			\end{equation}
		As $\int_{\Omega(g)}|\det(d\U)|dV_{g}\geq\int_{\Omega(g)}\det(d\U)dV_{g}$, it
        follows that
			\begin{equation}
				\int_{\Omega(g)}\mathrm{det}(g^{jk})dV_{g}\geq
				\frac{1}{|\Omega(g)|}{}
				\left(1-B\|R^{-}_{g}\|^{\frac1{12}}_{L^{2}(g)}\right)^{2}.
			\end{equation}
		From the mean value inequality, we know that there exists a point
		$x_{0}\in{}\Omega(g)$ such that
			\begin{equation}
				\mathrm{det}\left(g^{jk}(x_{0})\right)\geq
				\frac{1}{|\Omega(g)|^{2}}{}
				\left(1-B\|R^{-}_{g}\|^{\frac1{12}}_{L^{2}(g)}\right)^{2}.
			\end{equation}
		In order to turn this into a lower bound for all $x\in{}\Omega(g)$, we recall
		the first property of $\Omega(g)$, namely that
		$|g^{jk}-a^{jk}|\leq B\|R^{-}_{g}\|_{L^{2}(g)}^{\frac14}$ for all 
		$x\in{}\Omega(g)$.
		We may combine this with Proposition \ref{prop:det_estimate} to see that
			\begin{align}
				|\mathrm{det}\left(g^{jk}(x_{0})\right)-
				\mathrm{det}\left(a^{jk}\right)|&\leq
				C(n)\left(\|g^{jk}\|+\|a^{jk}\|\right)^{n-1}\|g^{jk}-a^{jk}\|
				\\
				&\leq
				B\|R^{-}_{g}\|^{\frac14}_{L^{2}(g)}.
			\end{align}
		One may use the triangle inequality to see that
			\begin{equation}
				\mathrm{det}\left(a^{jk}\right)\geq
				\frac{1}{|\Omega(g)|^{2}}{}
				\left(1-B\|R^{-}_{g}\|^{\frac1{12}}_{L^{2}(g)}\right)^{2}
				-B\|R^{-}_{g}\|^{\frac14}_{L^{2}(g)},
			\end{equation}
		and
			\begin{equation}
				\mathrm{det}\left(g^{jk}(x)\right)\geq
				\frac{1}{|\Omega(g)|^{2}}{}
				\left(1-B\|R^{-}_{g}\|^{\frac1{12}}_{L^{2}(g)}\right)^{2}
				-2B\|R^{-}_{g}\|^{\frac14}_{L^{2}(g)}.
			\end{equation}
		Thus, we have \eqref{eq:det_a_bounded_below} and
		\eqref{eq:det_g_bounded_below_Omega}, respectively.
		
		Now, recall that $\Omega(g)$ is connected, and so if
		$\mathrm{det}(d\U)$ changes signs on $\Omega(g)$, then there must be a point
		$x_{0}\in{}\Omega(g)$ such that
			\begin{equation}
				0=\mathrm{det}\bigl(d\U(x_{0})\bigr)^{2}=
				\mathrm{det}\left(g^{jk}(x_{0})\right).
			\end{equation}
		Since $|\Omega(g)|^2\leq |\T^{3}|_{g}\leq V$, there is a $B$ depending only
		on $V,R,\Lambda,\sigma,\eta$ and $M$ such that if 
			\begin{equation}
				\|R^{-}_{g}\|_{L^{2}(g)}\leq\frac1{B},
			\end{equation}
		then
			\begin{equation}
				\mathrm{det}\left(g^{jk}(x)\right)\geq
				\frac{1}{|\Omega(g)|^{2}}{}
				\left(1-B\|R^{-}_{g}\|^{\frac1{12}}_{L^{2}(g)}\right)^{2}
				-2B\|R^{-}_{g}\|^{\frac14}_{L^{2}(g)}
				\geq
				\frac12
			\end{equation}
		for all $x\in{}\Omega(g)$.
		In particular, the point $x_{0}\in{}\Omega(g)$ cannot exist in this case.
		Therefore, we see that $|\mathrm{det}(d\U(x))|>0$ for all $x\in{}\Omega(g)$.
		As a consequence, we have that
			\begin{equation}
				\int_{\Omega(g)}|\mathrm{det}(d\U)|dV_{g}=
				\pm\int_{\Omega(g)}\mathrm{det}(d\U)dV_{g}.
			\end{equation}
			It now follows from \eqref{eq:up_lower_int_bound_detU_Omega} that,
			once again, there is a $B$ depending only on $V,R$,$\Lambda,\sigma$,$\eta$
			and $M$ such that if
			\begin{equation}
				\|R^{-}_{g}\|_{L^{2}(g)}\leq\frac1{B},
			\end{equation}
			then
			\begin{equation}
				\int_{\Omega(g)}|\mathrm{det}(d\U)|dV_{g}=
				\int_{\Omega(g)}\mathrm{det}(d\U)dV_{g}.
			\end{equation}
		Thus, we have established \eqref{eq:detdU_has_sign_Omega}.
		This finishes the proof of the result.
	\end{proof}
	Now that we have established that the set $\Omega(g)$ has quite a few good properties
	with respect to the metrics $g\in{}\F$ and harmonic maps 
	$\U:(\T^{3},g)\rightarrow{}(\T^{3},h)$,
	we can begin to show that sequences of metrics whose negative part of their scalar
	curvatures tend to zero have sub-sequences converging to flat metrics.
	\begin{lemma}\label{lem:inj_on_well_approximating_set}
		Fix $V,R,\Lambda,\eta,M>0$, and for every $g\in{}\F(V,R,\Lambda,\sigma,\eta,M)$
		let $\U=\U_{g}$ be the harmonic map $\U:(\T^{3},g)\rightarrow{}(\T^{3}, h)$ given
		in Corollary \ref{cor:L2-bounded_deg-1-harmonic_map}, let
		$u^{j}$ be the components of $\U$, and let $g^{jk}$ denote
		$g\bigl(du^{j},du^{k}\bigr)$.
		Finally, for each $g$ in $\F$, let $a$ be the associated symmetric nonnegative
		matrix as in Proposition \ref{prop:const_matrix_approx}, and let
		$\Omega(g)$ be the set described in Lemma
		\ref{lem:well_approximating_sets_good_properties}.	
		Then, there exists a $B$ depending only on $V,R,\Lambda,\eta,M>0$ such that for
		any $g\in{}\F$ with
			\begin{equation}
				\|R^{-}_{g}\|_{L^{2}(g)}\leq\frac1{B}
			\end{equation}
		we have that there is an open neighborhood $W(g)$ containing
		$\overline{\Omega}(g)$ such that $\U$ restricted to $\overline{W}(g)$ is injective.
		In particular, we have that
			\begin{equation}
				\U\bigl(\partial{}\Omega(g)\bigr)
				=
				\partial{}\U\bigl(\Omega(g)\bigr).
			\end{equation}
	\end{lemma}
	\begin{proof}
		From the fact that $\mathrm{det}\left(g^{jk}\right)=\mathrm{det}\left(d\U\right)$
		and \eqref{eq:det_g_bounded_below_Omega}, we know that $\U$ is a local
		diffeomorphism around points of $\Omega(g)$ for all $g$ with
		$\|R^{-}_{g}\|_{L^{2}(g)}$ small enough, by continuity this is actually true
		of a small connected neighborhood $W(g)$ of $\overline{\Omega}(g)$, and in fact
		true of the compact closure $\overline{W}(g)$ of $W(g)$.
		We will begin by showing that $y\mapsto |\U^{-1}\{y\}\cap{}\overline{W}(g)|$
		is continuous on $\U\bigl(\overline{W}(g)\bigr)$, and so locally constant.
		Fix $y_{0}\in{}\U\bigl(\overline{W}(g)\bigr)$. Since $\overline{W}(g)$ is 
		compact and $\U$ is local diffeomorphism about every point in $\overline{W}(g)$,
		it must be that $\U^{-1}\{y_{0}\}\cap{}\overline{W}(g)$ is finite.
		As such, there is an open set $V_{y_{0}}$ about $y_{0}$ and open sets 
		$G_{i}$ for $i=1,\dots,n_{y_{0}}$ such that $\U$ is a diffeomorphism from
		$G_{i}$ to $V_{y_{0}}$, and the $G_{i}$ are pairwise disjoint and cover
		$\U^{-1}\{y_{0}\}\cap{}\overline{W}(g)$.
		This shows that for any $y\in{}V_{y_{0}}$ we have that $n(y)\geq n(y_{0})$.
		
		We will now argue by contradiction that $n$ is actually continuous at
		$y_{0}$.
		Suppose this were not the case, then we would be able to find a sequence of
		$y_{i}$ converging to $y_{0}$ such that $n(y_{i})>n(y_{0})$ for all $i$.
		In particular, this implies that for each $i$ there exists an $x_{i}$
		in $\U^{-1}\{y_{i}\}\cap{}\overline{W}(g)$ which is not in
		$\bigcup{}_{i=1}^{n(y_{0})}G_{i}$.
		Since $\overline{W}(g)$ is compact, a subsequence of the $x_{i}$, say
		$x_{i_{j}}$, converges to some element $x_{0}$.
		The continuity of $\U$ implies that $\U(x_{0})=y_{0}$.
		However, this implies that $x_{0}\in{}\bigcup_{i=1}^{n(y_{0})}G_{i}$.
		Since $x_{i_{j}}\rightarrow{}x_{0}$, for all $j$ big enough we must have
		that $x_{i_{j}}\in{}\bigcup_{i=1}^{n(y_{0})}G_{i}$.
		This is a contradiction to how we chose the sequence $x_{i}$ in the first place,
		and so we see that $y\mapsto n(y)$ is in fact continuous on $\overline{W}(g)$.
		Therefore, since $W(g)$, and so $\overline{W}(g)$, is connected,
		we see that $y\mapsto n(y)$ is in fact constant on $\overline{W}(g)$,
		and so constant on $W(g)$.

		Let $n=n(y_{0})$ be the value of $n$ on $\overline{W}(g)$.
		From the area formula we have that 
			\begin{equation}
				n|\U\bigl(\Omega(g)\bigr)|_{h}=
				\int_{\U\bigl(\Omega(g)\bigr)}n(y)dV_{h}(y)
				=\int_{\Omega(g)}|\mathrm{det}(d\U)|dV_{g}.
			\end{equation}
		Recalling \eqref{eq:detdU_has_sign_Omega} and
		\eqref{eq:up_lower_int_bound_detU_Omega}, we get that
			\begin{equation}
				n|\U\bigl(\Omega(g)\bigr)|_{h}\leq
				1+B\|R^{-}_{g}\|^{\frac1{12}}_{L^{2}(g)}.
			\end{equation}
		Furthermore, from \eqref{eq:U-image_compl_small_in_h}, we have that
			\begin{equation}
				|\U\bigl(\Omega(g)\bigr)|_{h}\geq
				1-B\|R^{-}_{g}\|^{\frac1{12}}_{L^{2}(g)},
			\end{equation}
		since the fact that $\mathrm{deg}(\U)=1$ implies that it is surjective, and so
		$\U\bigl(\Omega(g)\bigr)^{c}\subset{}\U\bigl(\Omega(g)^{c}\bigr)$.
		Thus, rearranging terms show us that
			\begin{equation}
				n\leq
				\frac{1+B\|R^{-}_{g}\|^{\frac1{12}}_{L^{2}(g)}}
				{1-B\|R^{-}_{g}\|^{\frac1{12}}_{L^{2}(g)}}{}.
			\end{equation}
		Therefore, there is a constant $B$ depending only on $V,R,\Lambda,\eta,M>0$
		such that if $\|R^{-}_{g}\|_{L^{2}(g)}\leq\frac1B$, then
		$n=1$, and $\left.\U\right|_{\overline{W}(g)}$ is injective.

		Since $\U$ is a local diffeomorphism around every point of $\overline{W}(g)$,
		it follows from the fact that $\U$ is injective when resticted to $\overline{W}(g)$
		that $\U\bigl(W(g)\bigr)$ is an open subset of $\T^{3}$ and
		$\left. \U\right|_{W(g)}$ is a diffeomorphism.
		Since $\overline{\Omega}(g)\subset{}W(g)$, we see that
		$\U\bigl(\overline{\Omega}(g)\bigr)=\mathrm{cl}\left(\U\bigl(\Omega(g)\bigr)\right)$.
		As such, we have that
		$\partial{}\U\bigl(\Omega(g)\bigr)\subset{}\U\bigl(W(g)\bigr)$.
		The result now follows from the fact that $\left. \U\right|_{W(g)}$ is a 
		diffeomorphism.
	\end{proof}
    Consider a sequence of metrics $g_i\in\F$ such that $\|R^{-}_{g_{i}}\|_{L^{2}(g_{i})}$
    converges to zero, and let $a_i$ be the corresponding sequence of matrices
    which approximate $g_i$ as in Proposition \ref{prop:const_matrix_approx}.
    As may have been guessed, the constant matrices $a_i$ will be used to show
    that a subsequence of the metrics $g_i$ converge to a flat metric in the
    sense of Dong-Song.
	The following Corollary gives this idea a clearer form.
	\begin{corollary}\label{cor:C0_close_to_flat_metric}
		Fix $V,R,\Lambda,\sigma,\eta, M>0$. For every $g\in{}\F$ let 
		$\U=\U_{g}:(\T^{3},g)\rightarrow{}(\T^{3},h)$ be the harmonic map given in
		Corollary \ref{cor:L2-bounded_deg-1-harmonic_map}, let $u^{j}$ denote its
		components, let $g^{jk}=g\bigl(du^{j},du^{k}\bigr)$, let $a$ be the symmetric
		and non-negative matrix approximating $g^{jk}$ as in Proposition 
		\ref{prop:const_matrix_approx}, and let $\Omega(g)$ be as in Lemma
		\ref{lem:well_approximating_sets_good_properties}.
		Then, there is a $B>0$ depending only on $V,R,\Lambda,\sigma,\eta, M$ such that
		for all $g\in{}\F$ with $\|R^{-}_{g}\|_{L^{2}(g)}\leq\frac1{B}$ we may find
		a flat metric $g_{F}$ on $\T^{3}$ such that
			\begin{equation}
				\|g-\U^{*}g_{F}\|_{g;\Omega(g)}\leq B\|R^{-}_{g}\|^{\frac1{4}}_{L^2(g)}.
			\end{equation}
	\end{corollary}
	\begin{proof}
		From \eqref{eq:det_a_bounded_below} of Lemma \ref{lem:well_approximating_sets_good_properties},
		we know that there
		is a constant $B$ depending only on $V,R,\Lambda,\sigma,\eta, M$ such that
		if $\|R^{-}_g\|_{L^{2}(g)}\leq\frac1B$, then the approximating symmetric
		matrix has the upper bound $\|a\|\leq B$ and lower bound
			\begin{equation}
				\mathrm{det}(a)\geq
				\frac{1}{|\Omega(g)|^{2}}{}
				\left(1-B\|R^{-}_{g}\|^{\frac1{12}}_{L^{2}(g)}\right)^{2}-
				B\|R^{-}_{g}\|^{\frac14}_{L^{2}(g)}.
			\end{equation}
		Therefore, from \eqref{eq:small_compliment} and the lower volume bound in
		Proposition \ref{prop:vol_lower_bound-harmonic_maps},
		there is a $B$ depending only on $V,R,\Lambda,\sigma,\eta$ and $M$
		such that if $\|R^{-}_{g}\|_{L^{2}(g)}\leq\frac1B$, then we have that
			\begin{equation}
				\mathrm{det}(a)\geq\frac1B.
			\end{equation}
		As such, we see that for some $B$ not depending on $g\in{}F$, the approximating
		matrix $a$ for $g$ is invertible with a uniform bound on its inverse:
			\begin{equation}
				\max\{\|a\|,\|a^{-1}\|\}\leq B.
			\end{equation}
		Letting $\theta^{j}$ denote the standard coordinates on $\T^{3}$, we may use the
		fact that $a$ is invertible to define the following flat metric on $\T^{3}$:
			\begin{equation}
				g_{F}=\left(a^{-1}\right)_{st}d\theta^{s}d\theta^{t}.
			\end{equation}
		We see that
			\begin{equation}
				\U^{*}g_{F}=\left(a^{-1}\right)_{st}du^{s}du^{t}.
			\end{equation}
		At this point, we want to estimate $\|g-\U^{*}g_{F}\|_{g;\Omega(g)}$.
		To do this, observe that from the bound on $\max \{\|a,\|,\|a^{-1}\|\}$
		and from \eqref{eq:const_approximation}
		in Lemma \ref{lem:well_approximating_sets_good_properties}, we have that
			\begin{equation}
				\max \left\{\|g^{jk}\|_{L^{\infty{}}\bigl(\Omega(g)\bigr)},
					\|\left(g^{jk}\right)^{-1}\|_{L^{\infty{}}\bigl(\Omega(g)\bigr)}
				\right\}
				\leq B.
			\end{equation}
		Thus, for any $x\in{}\Omega(g)$ and any two-form $\beta$ over  the point $x$,
		we have that
			\begin{equation}
				|\beta|_{g}^{2}\leq 
				B\sum_{jk}\beta\bigl(\nabla u^{j},\nabla u^{k}\bigr).
			\end{equation}
		In particular, we see that for any $x\in{}\Omega(g)$, we have that
			\begin{align}
				|g-\U^{*}g_{F}|^{2}_{g}&\leq
				B\sum_{jk}\bigl(g(x)-\U^{*}g_{F}\bigr)
				\bigl(\nabla u^{j},\nabla u^{k}\bigr)
				\\
				&\leq
				B\sum_{jk}\left(g^{jk}-(a^{-1})_{st}g^{sj}g^{tk}\right)
			\end{align}
		Adding and subtracting $a^{sj}$ gives
			\begin{equation}
				\sum_{jk}\left(g^{jk}-\left(\delta^{j}_{t}+
					(a^{-1})_{st}\left(g^{sj}-a^{sj}\right)\right)g^{tk}\right)
				=
				-\sum_{jk}g^{tk}(a)^{-1}_{st}\left(g^{sj}-a^{sj}\right).
			\end{equation}
		Taking the absolute value of the right hand side, and using the bounds
		on $\|a^{-1}\|$ and $\|g^{jk}\|_{L^{\infty{}}\bigl(\Omega(g)\bigr)}$, along with the
		estimate in \eqref{eq:const_approximation}, we see that
			\begin{equation}
				\|g-\U^{*}g_{F}\|_{g;\Omega(g)}\leq 
				B\|R^{-}_{g}\|^{\frac14}_{L^{2}(g)}.
			\end{equation}
	\end{proof}
	The above shows that there is a large set with a small boundary on which metrics
	in $\F$ with small negative part of their scalar curvature are close to some flat metric
	in the $C^{0}$ sense.
	This is precisely the setting of Dong-Song's convergence and approximation result
	\cite{Dong-Song-2023},
	which we recall here with minor modifications to suite the present situation a little
	better.
	\begin{lemma}\label{lem:Dong-Song_approx_lem}
		Let $(M,g)$ be any smooth closed three dimensional Riemannian manifold. For any $\varepsilon{}>0$
		there exists a $\delta>0$ such that if $\Omega$ is any connected open
		submanifold of $M$ with smooth boundary such that
			\begin{equation}
				|\Omega^{c}|+|\partial{}\Omega|\leq\delta,
			\end{equation}
		then we may find another open connected submanifold $\widetilde{\Omega}$
		with smooth boundary which satisfies the following properties:
			\begin{equation}
				\widetilde{\Omega}\subset{}\Omega;
			\end{equation}
			\begin{equation}
				|\widetilde{\Omega}^{c}|+|\partial{}\widetilde{\Omega}|\leq
				\varepsilon{};
			\end{equation}
		for every $z\in{}M$ we have that
			\begin{equation}
				d(z,\widetilde{\Omega})\leq\varepsilon{};
			\end{equation}
		and for every $x,y\in{}\widetilde{\Omega}$ there exists a curve
		$\gamma\subset{}\widetilde{\Omega}$ connecting them such that
			\begin{equation}
				\mathrm{L}(\gamma)\leq d^{g}(x,y)+\varepsilon{}.
			\end{equation}
		The last two conditions imply the following:
			\begin{equation}
				d_{GH}\left(
					\left(\widetilde{\Omega},\hat{d}^{g}_{\widetilde{\Omega}}
					\right),
					(M,d^{g})
				\right)\leq 2\varepsilon{}.
			\end{equation}
	\end{lemma}
	With all of the results up to now in hand, we can finally establish the main stability
	result of this paper.
	\begin{theorem}\label{thm:Dong-Song_conv_for_F}
		Fix $V,R,\Lambda,\eta,M>0$ and let $\F(V,R,\Lambda,\eta,M)$ be the family
		of Riemannian metrics on $\T^{3}$ given in Definition \ref{defn:family_of_metrics}.
		Let $g_{i}$ be a sequence of metrics in $\F$ such that
			\begin{equation}
				\lim_{i\rightarrow{}\infty{}}\|R^{-}_{g_{i}}\|_{L^{2}(g_{i})}=0.
			\end{equation}
		Then, there is a subsequence, also denoted $g_{i}$, and a flat metric $g_{F_{\infty}}$ on
		$\T^{3}$ such that $g_{i}$ converges to $g_{F_{\infty{}}}$ in the sense of Dong-Song.
		That is, for any $\varepsilon{}>0$ there exists an $N\in{}\N$ such that
		for all $i\geq N$ there is an open submanifold $\widetilde{\Omega}_{i}$ with smooth
		boundary such that
			\begin{equation}
				|\widetilde{\Omega}_{i}^{c}|_{g_{i}}+|\partial{}
				\widetilde{\Omega}_{i}|_{g_{i}}\leq
				\varepsilon{},
			\end{equation}
		and
			\begin{equation}
				d_{GH}\left(\left(\widetilde{\Omega}_{i},
					\hat{d}^{g_{i}}_{\widetilde{\Omega}_{i}}\right),
					\left(\T^{3},d^{g_{F_{\infty{}}}}\right)
				\right)\leq\varepsilon{}.
			\end{equation}
	\end{theorem}
	\begin{proof}
 		Consider any sequence of metrics $g_{i}\in{}\F$ such that
			\begin{equation}
				\lim_{i\rightarrow{}\infty{}}\|R^{-}_{g_{i}}\|_{L^{2}(g_{i})}=0.
			\end{equation}
		Let $\U_{i}$ denote the harmonic maps 
		$\U_{i}:(\T^{3},g_{i})\rightarrow{}(\T^{3},h)$ and let $a_{i}$ denote the
		non-negative symmetric matrices which approximate $g_{i}$ as in Proposition
		\ref{prop:const_matrix_approx}.
		Since the terms $\|R_{g_{i}}^{-}\|_{L^{2}(g)}$ are tending towards zero,
		we may always assume without loss of generality that $\|R^{-}_{g_{i}}\|_{L^{2}(g)}$
		is small enough so that Lemma \ref{lem:well_approximating_sets_good_properties},
		Lemma \ref{lem:inj_on_well_approximating_set}, and Corollary 
		\ref{cor:C0_close_to_flat_metric} apply to $g_{i}$, $\U_{i}$, $a_{i}$, $g_{F_{i}}$,
		and $\Omega_{i}=\Omega(g_{i})$.

		Then, we see that on $\Omega_{i}=\Omega(g_{i})$ we have that
			\begin{equation}
				\|g_{i}-\U^{*}_{i}g_{F_{i}}\|_{g_{i},\Omega_{i}}\leq
				\|R^{-}_{g_{i}}\|^{\frac14}_{L^{2}(g)}.
			\end{equation}
		Furthermore, we have for all $i$ that $\max\{\|a_{i}\|,\|a^{-1}_{i}\|\}\leq B$.
		As such, we see that there is a subsequence $a_{i(m)}$ and a symmetric positive
		definite matrix $a_{\infty{}}$ such that
			\begin{equation}\label{eq:a_i_conv_to_a_00}
				\lim_{m\rightarrow{}\infty{}}a_{i(m)}=a_{\infty{}}.
			\end{equation}
		Letting $g_{F_{\infty{}}}$ be defined by
			\begin{equation}
				g_{F_{\infty{}}}=\left(a^{-1}_{\infty{}}\right)_{st}
				d\theta^{s}d\theta^{t},
			\end{equation}
		we see that $g_{F_{i}}$ converges to $g_{F_{\infty{}}}$ in the $C^{0}$ sense.
		Now, we will estimate
			\begin{equation}
				\|\U^{*}_{i(m)}(g_{F_{i(m)}}-g_{F_{\infty{}}})\|
				_{g_{i(m)},\Omega_{i(m)}}.
			\end{equation}
		From the bound on $\max\{\|a_{i}\|,\|a^{-1}_{i}\|\}$ and 
		\eqref{eq:const_approximation}, we see that
		\begin{equation}\label{eq:bound_on_g_and_g_inv}
				\max\left\{\|g_{i(m)}^{jk}\|_{L^{\infty{}}\left(\Omega_{i(m)}\right)},
					\|(g^{jk}_{i(m)})^{-1}\|_{L^{\infty{}}
					\left(\Omega_{i(m)}\right)}
				\right\}\leq B
			\end{equation}
		So, we have that for all $x\in{}\Omega_{i}$ that
			\begin{align}
				\|\U^{*}_{i(m)}(g_{F_{i(m)}}-g_{F_{\infty{}}})\|
				_{g_{i(m)}}^{2}&\leq
				B\sum_{jk}\U^{*}_{i(m)}(g_{F_{i}}-g_{F_{\infty{}}})
				\left(\nabla u_{i(m)}^{j},\nabla u_{i(m)}^{k}\right)
				\\
			       &\leq B\sum_{jk}\left(\left(a^{-1}_{i(m)}\right)_{st}-
				       \left(a^{-1}_{\infty{}}\right)_{st}
			       \right)
			       g_{i(m)}^{sj}g_{i(m)}^{tk}.
			\end{align}
		So, combined with the bounds on $\max\{\|a_{i}\|,\|a^{-1}_{i}\|\}$ and
		\eqref{eq:bound_on_g_and_g_inv}, the above shows that
			\begin{equation}
				\lim_{m\rightarrow{}\infty{}}\|\U^{*}_{i(m)}
				(g_{F_{i(m)}}-g_{F_{\infty{}}})
				\|_{g_{i(m)},\Omega_{i(m)}}=0.
			\end{equation}
		Therefore, using the triangle inequality, we see that
			\begin{equation}\label{eq:subseq_C0_convergence}
				\lim_{m\rightarrow{}\infty{}}
				\|g_{i(m)}-\U^{*}_{i(m)}
				g_{F_{\infty{}}}\|_{g_{i(m)},\Omega_{i(m)}}=0.
			\end{equation}

		Now, we recall that because 
		$\lim_{i\rightarrow{}\infty{}}\|R^{-}_{g}\|_{L^{2}(g)}=0$, we have that
		Lemma \ref{lem:inj_on_well_approximating_set} applies to $g_{i}$ for all
		$i$ large enough.
		In particular, we may assume without loss of generality that
		$\left.\U_{i(m)}\right|_{\Omega_{i(m)}}$ is injective for all $m$.
		Thus, it follows from \eqref{eq:subseq_C0_convergence} that
		\begin{equation}\label{eq:Lipschitz_Convergence_Length}
			\lim_{m\rightarrow{}\infty{}}
			\max\left\{\mathrm{Lip}_{\hat{d}^{g_{i(m)}}_{\Omega_{i(m)}}}
				\left(\U_{i(m)}\right),
				\mathrm{Lip}_{\hat{d}^{g_{F_{\infty{}}}}_{
					\U_{i(m)}\left(\Omega_{i(m)}\right)}}
				\left(\U_{i(m)}^{-1}\right)
				\right\}=1
		\end{equation}	

		From the fact that $g_{F_{\infty{}}}$ and $h$ are two fixed Riemannian metrics on
		$\T^{3}$, they are uniformly comparable in the sense that there is a constant $C$
		such that for any $p\in{}\T^{3}$ and any $\nu\in{}T_{p}\T^{3}$ we have that
			\begin{equation}
				\frac1C\leq\frac{g_{F_{\infty{}}}(\nu,\nu)}{h(\nu,\nu)}{}
				\leq C
			\end{equation}
		In particular, for any subset $W\subset{}\T^{3}$ we have that
			\begin{equation}
				\frac1{C^{3}}\leq\frac{|W|_{g_{F_{\infty{}}}}}{|W|_{h}}{}
				\leq C^{3}.
			\end{equation}
		Since each $\U_{i(m)}$ has degree 1, and is therefore necessarily surjective
		by Proposition \ref{prop:surjectivity_degree_1_maps_between_tori},
		we have that $\U_{i(m)}\left(\Omega_{i(m)}\right)^{c}
		\subset{}\U_{i(m)}(\Omega_{i(m)}^{c})$.
		Therefore, equation \eqref{eq:small_compliment} of Lemma
		\ref{lem:well_approximating_sets_good_properties} shows that
			\begin{equation}
				|\U_{i(m)}\left(\Omega_{i(m)}\right)|_{g_{F_{\infty{}}}}\leq 
				B\|R^{-}_{g}\|^{\frac14}_{L^{2}\left(g_{i(m)}\right)}.
			\end{equation}
		Similarly, for any $2$ dimensional smooth submanifold $\Sigma\subset{}\T^{3}$
		we have that
			\begin{equation}
				\frac1{C^{2}}\leq
				\frac{|\Sigma|_{g_{F_{\infty{}}}}}{|\Sigma|_{h}}{}
				\leq C^{2},
			\end{equation}
		and so Lemma \ref{lem:inj_on_well_approximating_set}, along with
		the bounds on $\left|du_{i(m)}^{j}\right|_{g_{i(m)}}$ for points in $\Omega_{i(m)}$,
		implies that
			\begin{equation}
				|\partial{}\U\left(\Omega_{i(m)}\right)|_{g_{F_{\infty{}}}}=
				|\U\left(\partial{}\Omega_{i(m)}\right)|_{g_{F_{\infty{}}}}\leq
				B|\partial{}\Omega_{i(m)}|_{g_{i(m)}}
				\leq
				B\|R^{-}_{g_{i(m)}}\|^{\frac14}_{L^{2}(g_{i(m)}}.
			\end{equation}
		Therefore, we may apply Lemma \ref{lem:Dong-Song_approx_lem} to conclude that
		for all $m$ large enough we may find open submanifolds
		$\Omega'_{i(m)}\subset{}\U_{i(m)}\left(\Omega_{i(m)}\right)$
		with smooth boundary so that
			\begin{align}
				&\lim_{m\rightarrow{}\infty{}}
				|\Omega^{'c}_{i(m)}|_{g_{F_{\infty{}}}}=0
				\\
				&\lim_{m\rightarrow{}\infty{}}
				|\partial{}\Omega'_{i(m)}|_{g_{F_{\infty{}}}}=0
				\\
				&\lim_{m\rightarrow{}\infty{}}d_{GH}\left(
					\left(\Omega'_{i(m)},
						\hat{d}^{g_{F_{\infty{}}}}_{
							\Omega'_{i(m)}}\right),	
					\left(\T^{3},d^{g_{F_{\infty{}}}}\right)
				\right)=0.
			\end{align}
		By the above work, we may also conclude that for all $m$ large enough, we may set
		$\widetilde{\Omega}_{i(m)}=\U_{i(m)}^{-1}\Omega'_{i(m)}$ and obtain
			\begin{align}
				&\lim_{m\rightarrow{}\infty{}}
				|\widetilde{\Omega}^{c}_{i(m)}|_{g_{i(m)}}=0
				\\
				&\lim_{m\rightarrow\infty}
				|\partial{}\widetilde{\Omega}_{i(m)}|_{g_{i(m)}}=0
				\\
				&\lim_{m\rightarrow\infty}d_{GH}
				\left(
					\left(\widetilde{\Omega}_{i(m)},
						\hat{d}^{g_{i(m)}}_{\widetilde{\Omega}_{i(m)}}
					\right),
					\left(\Omega'_{i(m)},
					\hat{d}^{g_{F_{\infty{}}}}_{\Omega'_{i(m)}}
				\right)
				\right)=0.
			\end{align}			
		where the map $\U_{i(m)}$ gives the estimate on the above estimate on the
		Gromov-Hausdorff distance.
		Therefore, the result now follows from the triangle inequality and 
		Equation \eqref{eq:Lipschitz_Convergence_Length}.

	\end{proof}
	The convergence results stated in the introduction will now follow if we can show that
	the two mentioned families of metrics both lie in $\F(V,R,\lambda,\sigma,\eta,M)$ for 
	some values of $V,R,\lambda,\sigma,\eta$, and $M$. 

	\begin{theorem}
		Let $\sigma,K,V,D>0$, and define $\mathcal{R}=\mathcal{R}(\sigma,K,D)$ to be 
		the collection of Riemannian metrics $g$ on $\T^3$ such that 
			\begin{enumerate}
				\item $\min\{\mathrm{stsys}_1(g),\mathrm{stsys}_{2}(g)\}\geq\sigma$;
				\item $\mathrm{Ric}_{g}\geq-K$;
				\item $\diam_{g}(\T^3)\leq D$.
			\end{enumerate}
		Then, for any sequence of metrics $\{g_{i}\}_{i=1}^{\infty{}}\subset{}
		\mathcal{R}(\sigma,K,V)$ such that 
			\begin{equation}
				\lim_{i\rightarrow{}\infty{}}\|R^{-}_{g_{i}}\|_{L^{2}(g_{i})}=0,
			\end{equation}
		there exists a subsequence $\{g_{i_{j}}\}_{j=1}^{\infty{}}$ and a
		flat metric $g_{F_{\infty}}$ on $\T^{3}$ such that $g_{i_{j}}\rightarrow{}g_{F_{\infty}}$ in the
		sense of Dong-Song.
	\end{theorem}
	\begin{proof}
		Once we show that 
		$\mathcal{R}(\sigma,K,V)\subset{}\F(V_{0},R_{0},\Lambda_{0},
		\sigma_{0},\eta_{0},M_{0})$ for some values 
		$V_{0},R_{0},\Lambda_{0},\sigma_{0},\eta_{0},M_{0}$, then we are
		done.
		Using volume comparison and the diameter bound, we may find an appropriate $V_{0}$.
        Furthermore, we may take $\sigma_{0}=\sigma$, by definition.
		It is then standard theory for smooth manifolds with Ricci curvature lower bounds
		that $\ison_{1}(g)\geq{\Lambda}_{0}$, 
		where $\Lambda_{0}$
		depends only on $K$ and $\mathrm{diam}_g(\T^3)$,
		see for example \cite[Theorem 114]{berger2003panoramic}, \cite[Theorem 3]{gallot1988isoperimetric},
        \cite[Page 294]{Petersen1998IntegralCB}.
        Furthermore, see \cite{DAI20181} for more related results.
        
        Next, we have that $R_{g}\geq -nK$, so we may take $R_0=nKV^{\frac12}_0$.
		It remains to show that we can find $\eta_{0}>0$ and $M_{0}$ such that
		for all $g\in{}\mathcal{R}$ we have that
			\begin{equation}
				\kappa(g,\eta_{0})\leq M.
			\end{equation}
		Fix $a\in{}\T^{3}$ and $\hat{a}_{0}\in{}\pi^{-1}\{a\}$, and let
		$\mathrm{Dir}(\hat{a}_{0})$ be the fundamental domain about
		$\hat{a}_{0}$ given in Lemma \ref{lem:fund_domain_exists}.
		We claim that 
		$\mathrm{diam}(\mathrm{Dir}(\hat{a}_{0}))\leq2 \mathrm{diam}_{g}(\T^{3})$.
		To see this, let $\hat{y}$ be an arbitrary element of $\mathrm{Dir}(\hat{a}_{0})$.
		By its construction we see that for any $\nu\in{}\Z^{3}$ we must have 
		$d(\hat{a}_{0},\hat{y})\leq d(\hat{a}_{\nu},\hat{y})$.
		Let $\gamma$ be any length minimizing geodesic connecting $a$ to $y=\pi(\hat{y})$,
		and let $\widetilde{\gamma}$ be the geodesic lifting $\gamma$ starting at
		$\hat{y}$.
		Then, we see that there is an $a_{\nu_{0}}$ such that 
			\begin{equation}
				d(\hat{a}_{\nu_{0}},\hat{y})=d(a,y)\leq \mathrm{diam}_{g}(\T^{3}).
			\end{equation}
		But, then we must have that
			\begin{equation}
				d(\hat{a}_{0},\hat{y})\leq
				d(\hat{a}_{\nu_{0}},\hat{y})\leq\mathrm{diam}_{g}(\T^{3}).
			\end{equation}
		This gives the claim.
		As such, for all $g\in{}\mathcal{R}$, we have that
			\begin{equation}
				\mathrm{diam}(\mathrm{Dir}(\hat{a}_{0})\leq2D.
			\end{equation}

		Once again, let $\hat{y}$ be an arbitrary element of 
		$\mathrm{Dir}(\hat{a}_{0})$, let $\eta_{0}=\frac{\sigma}{100}{}$,
		and let
			\begin{equation}
				J_{\hat{y}}=\{\nu\in{}\Z^{3}:B\left(\hat{y},\frac{\sigma}{100}\right)
					\cap{}\mathrm{Dir}(\hat{a}_{\nu})\neq\emptyset
				\}.
			\end{equation}
		For all $\nu\in{}J_{\hat{y}}$, let
        \begin{equation}
            \hat{z}_{\nu}\in B\left(\hat{y},\frac{\sigma}{100}\right)\cap
            \mathrm{Dir}(\hat{a}_{\nu}).
        \end{equation}
        Then, we have that 
			\begin{align}
				d(\hat{y},\hat{a}_{\nu})\leq d(\hat{z}_{\nu},\hat{a}_{\nu})+\frac{\sigma}{100}
				&\leq d(\hat{z}_{\nu},\hat{a}_{0})+\frac{\sigma}{100}
				\\
                &\leq d(\hat{y},\hat{a}_0)+\frac{\sigma}{50}.
			\end{align}
		In particular, for all $\nu\in{}J_{\hat{y}}$ we have that
			\begin{equation}
				d(\hat{a}_{\nu},\hat{a}_{0})\leq 2D+\frac{\sigma}{25}.	
			\end{equation}
		Furthermore, as $\stsys_{1}(g)\leq \mathrm{sys}_{1}(g)$, for all $\mu,\nu\in{}J$ we
		must have that
			\begin{equation}
				d(\hat{a}_{\mu},\hat{a}_{\nu})\geq\sigma,
			\end{equation}	
		and so 
		$B(\hat{a}_{\mu},\frac{\sigma}{4})\cap{}B(\hat{a}_{\nu}, \frac\sigma4)=\emptyset$.
		Finally, from the diameter upper bound, and volume lower bound given by the systolic
        inequality and the lower bound on $\sigma$,
        we see that there are constants $p_{0}$ and $P_{0}$ not depending on
		$g\in{}\mathcal{R}$ such that any ball
		of radius $\frac\sigma4$ has volume at least $p_{0}$.
		Thus, we see that
			\begin{equation}
				|J_{\hat{y}}|p_{0}\leq|B\left(\hat{a}_{0},
					2D+\sigma
				\right)|_{\pi^{*}g}\leq P_{0}.
			\end{equation}
		Rearranging terms shows that
			\begin{equation}
				|J_{\hat{y}}|\leq\frac{P_{0}}{p_{0}}{}.
			\end{equation}
		Now, by the doubling property and the volume upper bound,
        we may find a cover of $\mathrm{Dir}(\hat{a}_{0})$
		by at most $G$ balls of radius $\frac\sigma{100}$, where $G$ does not depend
		on $g\in{}\mathcal{R}$.
		Any such cover will also cover $\mathrm{Dir}(\hat{a}_{0})_{\frac\sigma{200}}$.
		In particular, we see that
			\begin{equation}
				\kappa(g,\frac\sigma{200})\leq\frac{GP_{0}}{p_{0}}.
			\end{equation}
		We can now apply Theorem \ref{thm:Dong-Song_conv_for_F} to $\mathcal{R}$ to
		obtain the result.
	\end{proof}
	The result also follows for families of Riemannian metrics which have a uniform
	lower bound in terms of some background metric on $\T^{3}$, as in the following
	theorem.
	\begin{theorem}\label{thm:VADB_Dong_Song}
		Let $g_0$ be a fixed Riemannian metric on $\T^3$, and let $\Lambda,R,V>0$.
		We denote by $\mathcal{V}(g_{0},\Lambda,V)$ the collection of Riemannian metrics
		$g$ on $\T^{3}$ satisfying the following properties:
			\begin{enumerate}
				\item $g\geq g_{0}$;
				\item $\ison_{1}(g)\geq\Lambda$;
				\item $\|R^{-}_{g}\|_{L^{2}(g)}\leq R$;
				\item $|\T^{3}|_{g}\leq V$.
			\end{enumerate}
		Then, for any sequence of metrics $\{g_{i}\}_{i=1}^{\infty{}}\subset{}
		\mathcal{R}(\sigma,K,V)$ such that 
			\begin{equation}
				\lim_{i\rightarrow{}\infty{}}\|R^{-}_{g_{i}}\|_{L^{2}(g_{i})}=0,
			\end{equation}
		there exists a subsequence $\{g_{i_{j}}\}_{j=1}^{\infty{}}$ and a
		flat metric $g_{F_{\infty}}$ on $\T^{3}$ such that $g_{i_{j}}\rightarrow{}g_{F_{\infty}}$ in the
		sense of Dong-Song.
	\end{theorem}
	\begin{proof}
		We need to show that there are $V_{0},R_{0},\Lambda_{0},\sigma_{0},\eta_{0}$,
		and $M_{0}$ such that $\mathcal{V}\subset{}\F$.
		We may immediately take $V_{0}=V,R_{0}=R$, and $\Lambda_{0}=\Lambda$.
		Since $g\geq g_{0}$, we have that
			\begin{equation}
				\min\bigl\{\stsys_{1}(g),\stsys_{2}(g)\bigr\}
				\geq
				\min\bigl\{\stsys_{1}(g_{0}),\stsys_{2}(g_{0})\bigr\}.
			\end{equation}
		Therefore, we may set 
		$\sigma_{0}=\min\bigl\{\stsys_{1}(g_{0}),\stsys_{2}(g_{0})\bigr\}$.

		Let $\eta_{0}>0$ be arbitrary, then there is a $M_{0}$ such that
		$\kappa(g_{0},\eta_{0})\leq M_{0}$.
		Let $\V$ be a test domain for $\kappa(g_{0},\eta_{0})$.
		Since
			\begin{equation}
				\{x:d^{g}(x,\V)<\eta_{0}\}\subset{}
				\{x:d^{g_{0}}(x,\V)<\eta_{0}\},
			\end{equation}
		it follows that $\kappa(g,\eta_{0})\leq\kappa(g_{0},\eta)$.
		Therefore, we may take $M_{0}$ for our final piece of the puzzle,
		and apply Theorem \ref{thm:Dong-Song_conv_for_F} to get the result.
	\end{proof}
 If in addition to a metric lower bound the sequence has volumes converging to the volume of
 the metric lower bound, then we actually have that $g_0$ is flat.
 \begin{corollary}
     Let $g_0$ be a fixed Riemannian metric on $\T^3$, let $\Lambda, R, V>0$,
     and let $\mathcal{V}=\mathcal{V}(g_0,\Lambda,R,V)$ denote the collection of
     Riemannian metrics $g$ on $\T^3$ which satify the following properties:
     \begin{enumerate}
				\item $g\geq g_{0}$;
				\item $\ison_{1}(g)\geq\Lambda$;
				\item $\|R^{-}_{g}\|_{L^{2}(g)}\leq R$;
				\item $|\T^{3}|_{g}\leq V$.
			\end{enumerate}
   Suppose that $g_i$ is a sequence of metrics in $\mathcal{V}$ such that
   \begin{equation}
       |\T^3|_{g_i}\rightarrow |\T^3|_{g_0},
   \end{equation}
   and
   \begin{equation}
       \lim_{i\rightarrow\infty}\|R^{-}_{g_i}\|_{L^2(g_i)}=0.
   \end{equation}
   Then, the metric $g_0$ must be flat.
 \end{corollary}
 \begin{proof}
    		By Theorem \ref{thm:VADB_Dong_Song}, we know that a subsequence of 
		$\{g_{i}\}_{i=1}^{\infty{}}$, also denoted $\{g_{i}\}_{i=1}^{\infty{}}$,
		converges to a flat metric $g_{F}$ in the Dong-Song sense.
		Let $\U_{i}$ and $\widetilde{\Omega}_{i}$ be as in
		Theorem \ref{thm:VADB_Dong_Song}:
		\begin{equation}
			\lim_{m\rightarrow{}\infty{}}\left(|\widetilde{\Omega}_{i}^{c}|_{g_{i}}
			+|\partial{}\widetilde{\Omega}_{i}|_{g_{i}}\right)=0,
		\end{equation}
		$\U_{i}|_{\widetilde{\Omega}_{i}}$ is injective, and it induces
		\begin{equation}
			\lim_{m\rightarrow{}\infty{}}d_{\mathrm{GH}}\left(
			\left(\widetilde{\Omega}_{i(m)},\hat{d}^{g_{i}}_{\widetilde{\Omega}_{i}}
			\right),
			(\T^{3},d^{g_{F}})
			\right)
			=0.
		\end{equation}
		We will now observe several consequences of the inequality $g_{i}\geq g_{0}$:
		\begin{align}
			&\hat{d}^{g_{i}}_{\widetilde{\Omega}_{i}}\geq
			\hat{d}^{g_{0}}_{\widetilde{\Omega}_{i}};
			\label{eq:length_metric_lowerbound}
			\\
			&\lvert{}\widetilde{\Omega}_{i}^{c}{}\rvert_{g_{i}}\geq
			\lvert{}\widetilde{\Omega}_{i}^{c}\rvert_{g_{0}}.
			\label{eq:volume_lower_bound}
		\end{align}
		From \eqref{eq:length_metric_lowerbound} above we have the following sequence
		of inequalities:
		\begin{equation}\label{eq:distance_inequalities}
			\left.d^{g_{0}}
				\right|_{\widetilde{\Omega}_{i}\times{}\widetilde{\Omega}_{i}}
			\leq
			\hat{d}^{g_{0}}_{\widetilde{\Omega}_{i}}
			\leq
			\hat{d}^{g_{i}}_{\widetilde{\Omega}_{i}}.
		\end{equation}
		Furthermore, from \eqref{eq:volume_lower_bound} we have that
		\begin{equation}
			\lim_{i\rightarrow{}\infty{}}
			\lvert{}\widetilde{\Omega}^{c}_{i}{}\rvert_{g_{0}}=0.
		\end{equation}
		Since $g_{0}$ is a smooth Riemannian manifold, balls have a lower bound on
		their volume growth.
		Combined with the above limit, this shows that the inclusion map gives
		\begin{equation}
			\lim_{i\rightarrow{}\infty{}}d_{\mathrm{GH}}\left(
			\left(\widetilde{\Omega}_{i},d^{g_{i}}\right),
			(\T^{3},d^{g_{0}})
			\right)=0.
		\end{equation}
		Next for each $i$, the identity map
		\begin{equation}
			\mathrm{Id}:(\widetilde{\Omega}_{i},\hat{d}^{g_{i}}_{\widetilde{\Omega}_{i}})\rightarrow{}
			(\widetilde{\Omega}_{i},d^{g_{0}})
		\end{equation}
		is a Lipschitz one map.
		Since 
		\begin{equation}
			\lim_{m\rightarrow{}\infty{}}d_{\mathrm{GH}}\left(
			\left(\widetilde{\Omega}_{i(m)},\hat{d}^{g_{i}}_{\widetilde{\Omega}_{i}}
			\right),
			(\T^{3},d^{g_{F}})
			\right)
			=0.
		\end{equation}
		and
		\begin{equation}
			\lim_{i\rightarrow{}\infty{}}
			\lvert{}\widetilde{\Omega}^{c}_{i}{}\rvert_{g_{0}}=0,
		\end{equation}
		it follows from the Arzela-Ascoli theorem that there is a Lipschitz one map
        \begin{equation}
            F:(\T^{3},h)\rightarrow{}(\T^{3},g_{0}).
        \end{equation}
		However, from the definition of $g_{F}$, we have that 
		$\lvert{}\T^{3}{}\rvert_{g_{F}}
		=\lim_{i\rightarrow{}\infty{}}\lvert{}\T^{3}{}\rvert_{g_{i}}$.
		Furthermore, it was an hypothesis that 
		$\lim_{i\rightarrow{}\infty{}}\lvert{}\T^{3}{}\rvert_{g_{i}}=
		\lvert{}\T^{3}{}\rvert_{g_{0}}$.
		Therefore, we must have that the map $F$ is an isometry, and so $g_{0}$
		is flat.

 \end{proof} 

	\section{Appendix}
	
\subsection{Dong-Song Curve Approximation}
The goal of this section is to prove the following approximation lemma. The proof presented here
is a minor modification of the proof found in \cite{Dong-Song-2023}.
\begin{lemma}\label{lem:Dong-Song_Curve_Approximation}
    Let $(M^{3},g)$ be any smooth closed three dimensional Riemannian manifold.
    For any $\varepsilon>0$ there exists a $\delta>0$ such that if $\Omega$ is a connected open
    sub-manifold of $M$ with smooth boundary such that $|\Omega|^c+|\partial\Omega|\leq\varepsilon$,
    then we may find $\widetilde{\Omega}\subset\Omega$, another open sub-manifold with smooth
    boundary, such that 
    $\abs{\widetilde{\Omega}^{c}}+\abs{\partial\widetilde{\Omega}}\leq\varepsilon$, and
    for every $x,y\in\widetilde{\Omega}$ there exists a curve $\gamma\subset\widetilde{\Omega}$
    connecting them such that
        \begin{equation}
            \mathrm{L}\bigl(\gamma\bigr)\leq d(x,y)+\varepsilon.
        \end{equation}
\end{lemma}
\begin{lemma}\label{lem:large_connected_subset}
    Let $\eta,\Lambda>0$ be fixed, let $(M,g)$ be a compact smooth $n-$dimensional 
    Riemannian manifold, and suppose that $\ison_1(M,g)\geq\Lambda$.
    Then, there exists a $\delta>0$ depending only on $\eta$ and $\Lambda$ such that if $E$ is a 
    smooth $n-$dimensional submanifold of $M$ such that $\abs{E}\geq\eta$ and 
    $\abs{\partial E}\leq\delta$, then $E$ has a connected component, say $\Omega$, such that
        \begin{equation}
            \abs{\Omega^c}\leq\frac{1}{\Lambda}\abs{\partial\Omega}\leq\frac{1}{\Lambda}
            \abs{\partial E}.
        \end{equation}
\end{lemma}
\begin{proof}
    Let $E_{i}$ be the connected components of $E$.
    Since $E$ is a smooth $n$ dimensional sub manifold, it follows that
        \begin{equation}
            \partial E=\bigsqcup_i\partial E_i.
        \end{equation}
    Suppose that there is an $i_0$ such that $\abs{E_{i_0}}\geq\frac12\abs{M}$. Then, taking 
    $\Omega=E_{i_0}$, we have from the definition of $\ison_{1}(M,g)$ that
        \begin{equation}
            \abs{\Omega^c}\leq\frac{1}{\Lambda}\abs{\partial\Omega}.
        \end{equation}
    Suppose on the contrary that there is no such $i_0$: for all $i$ we have that $\abs{E_i}<\frac12\abs{M}$. Then, it follows from the definition of $\ison_{1}(M,g)$ that
        \begin{equation}
            \eta\leq\sum_{i}\abs{E_i}\leq\frac{1}{\Lambda}\sum_{i}\abs{\partial E_i}
            =\frac{\abs{\partial E}}{\Lambda}.
        \end{equation}
    Choosing $\abs{\partial E}<\eta\Lambda$ gives a contradiction, and hence we obtain the result.
\end{proof}
\begin{lemma}\label{lem:two_dim_curve_mod}
    Let $(M,g)$ be a smooth, compact, closed, two-dimensional Riemannian manifold.
    Let $A$ be a two-dimensional sub manifold of $M$ with smooth boundary.
    Furthermore, suppose that every component of $\partial A$ bounds a two-dimensional 
    submanifold of $M$, and that $A$ is connected.
    Then, for any $x,y$ in $A$ there is a curve $\gamma$ lying entirely in $A$ connecting $x$ to $y$
    such that
        \begin{equation}
            L(\gamma)\leq d(x,y)+\abs{\partial A}.
        \end{equation}
\end{lemma}
\begin{proof}
    By our hypotheses, we have that $M\setminus A$ is the disjoint and finite union of connected open two dimensional submanifolds $B_i$, with smooth boundaries $\partial B_i$.
    We claim that the boundaries $\partial B_i$ must be connected, and so correspond to the connected components of $\partial A$.
    For any $i$, suppose that $\Sigma_1$ and $\Sigma_2$ are two connected components of $\partial B_i$.
    Such connected components can only bound connected regions.
    Therefore, without loss of generality, we may find disjoint connected regions $\Omega_1$ and $\Omega_2$ with boundaries $\Sigma_1$ and $\Sigma_2$ such that
    $\Omega_1\cap B_{i}=\Omega_{2}\cap B_2=\emptyset$ and such that $\Omega_1\cap A\neq\emptyset$ as well as $\Omega_2\cap A\neq\emptyset$.
    Let $x\in\Omega_1\cap A$ and $y\in\Omega_2\cap A$.
    Then, any curve connecting $x$ to $y$ must pass through $B_i$, which contradicts the assumption that $A$ is connected.

    Let $\gamma_0$ be a unit speed length minimizing geodesic connecting $x$ to $y$.
    We will modify $\gamma_0$ so that it lies entirely in $A$.
    First, we let $s_1$ be the first time that $\gamma_0$ lies in $\overline{B}_1$ and we let $t_1$ be the last time that $\gamma_0$ lies in $\overline{B}_1$.
    Since $\partial B_i$ is connected for all $i$, we may connect $\gamma_0(s_1)\in\partial B_1$ to $\gamma_0(t_1)\in\partial B_1$ by a curve $c_1$ lying in $\partial B_1$ which has length less than $\abs{\partial B_1}$, though is not necessarily unit speed.
    We define $\gamma_1$ to be equal to $\gamma_0$ on $[0,L(\gamma_0)]\setminus[s_1,t_1]$ and equal to $c_1$ on $[s_1,t_1]$.
    We produce $\gamma_{i+1}$ inductively as follows.
    Let $s_{i+1}$ be the first time that $\gamma_i$ lies in $\overline{B}_{i+1}$, and let $t_{i+1}$ be the last time $\gamma_{i}$ lies in $\overline{B}_{i+1}$.
    Since $\partial B_i$ is connected for all $i$, we may connect $\gamma_i(s_{i+1})\in\partial B_{i+1}$ to $\gamma_i(t_{i+1})\in\partial B_{i+1}$ by a curve $c_{i+1}$ lying in $\partial B_{i+1}$ which has length less than $\abs{\partial B_{i+1}}$.
    Then we let $\gamma_{i+1}$ be $\gamma_{i}$ on the interval $[0,L(\gamma_0)]\setminus[s_{i+1},t_{i+1}]$, and equal to $c_{i+1}$ on $[s_{i+1},t_{i+1}]$.
    Since there are only finitely many regions $B_i$ to consider, this process terminates at some $\gamma_{P}$.
    We then have that
        \begin{equation}
            L\left(\gamma_{P}\right)\leq L(\gamma_0)+\sum_{i=1}^{P}\abs{\partial B_i}\leq d(x,y)+\abs{\partial A}.
        \end{equation}
\end{proof}
\begin{lemma}\label{lem:locally_connected_subsets_of_Omega_in_balls}
    Let $(M^3,g)$ be a smooth, compact, closed, three-dimensional Riemannian manifold. For any $L>1$ let $r_L$ be such that balls with radius less than or equal to $r_L$ are geodesically convex, and the exponential map is $L$ bi-Lipschitz.
    Let $r$ be such that $24r\leq r_L$. Then, there is a constant $\delta(r)$, depending only on $r$, such that if $\abs{\Omega^c}+\abs{\partial\Omega}\leq\delta(r)$, then the following is true.
    In every ball $B(a,20r)$ we may find a connected subset $D_a\subset\Omega\cap B(a,20r)$ with the following properties.
        \begin{enumerate}
            \item For all $p,q\in D_a$ there is a curve $\gamma(p,q)=\gamma$ connecting them, which lies entirely in $D_a$.
                  Furthermore, we have that
                    \begin{equation*}
                        L(\gamma)\leq96r+8rL\pi+2r^{-1}\abs{B(a,20r)\cap\partial\Omega}.
                    \end{equation*}
            \item We have that
                    \begin{equation}
                        \abs{B(a,20r)\setminus D_a}\leq2(24^3L^5)\left(\frac{2L^5}{\ison_1(\mathbb{S}^2)}+L^4\right)\abs{B(a,20r)\cap\partial\Omega}.
                    \end{equation}
        \end{enumerate}
\end{lemma}
\begin{proof}
    Let $B(a,20r)$ be an arbitrary ball of radius $20r$, and let $\nu\in T_{a}M$ be an arbitrary unit vector.
    Set $x=\exp_{a}(-4r\nu)$ and $y=\exp_{a}(4r\nu)$.
    For the moment, we will focus our attention to a neighborhood of $x$.
    Let $f:\partial\Omega\rightarrow\R$ be given by $f(z)=d(z,x)$, and observe that $f$ is smooth on $\bigl(B(x,24r)\setminus B(x,r)\bigr)\cap\partial\Omega$.
    From the coarea formula, we have that
        \begin{equation}
            \int_{r}^{2r}\abs{f^{-1}\{t\}}dt=\int_{f^{-1}[r,2r]}\abs{\nabla f}dA_{g}\leq\abs{B(x,2r)\cap\partial\Omega}.
        \end{equation}
    Therefore, using Sard's Lemma and a mean-value inequality, we may find $\sigma_x\in[r,2r]$ such that $\partial B(x,\sigma_x)\cap\partial\Omega=f^{-1}\{\sigma_x\}$ is a smooth submanifold of $\partial\Omega$ and
        \begin{equation}\label{eq:boundary_area_estimate_G_x}
            \abs{\partial B(x,\sigma_x)\cap\partial\Omega}\leq\frac{1}{r}\abs{B(x,2r)\cap\partial\Omega}\leq\frac{1}{r}\abs{B(a,20r)\cap\partial\Omega}.
        \end{equation}
    We will define two sets on $\partial B(x,\sigma_x)$, estimate their volumes, and then study their interaction.
    To begin, let $G_x=\Omega\cap\partial B(x,\sigma_x)$, which has smooth boundary given by
        \begin{equation}
            \partial G_x=\partial\bigl(\Omega\cap\partial B(x,\sigma_x)\bigr)=\partial\Omega\cap\partial B(x,\sigma_x).
        \end{equation}
    From \eqref{eq:boundary_area_estimate_G_x}, we see that $\abs{\partial G_x}\leq\tfrac{1}{r}\abs{B(a,20r)\cap\partial\Omega}$.
    Since $\exp_x$ is $L$ bi-Lipschitz, it follows that 
        \begin{equation}
            \ison_1\bigl(\partial B(x,\sigma_x)\bigr)\geq\frac{\ison_1\bigl(\partial B(0,\sigma_x)\bigr)}{L^{5}}\geq\frac{1}{\sigma_xL^5}\ison_1(\mathbb{S}^2)\geq\frac{1}{2rL^5}\ison_1(\mathbb{S}^2).
        \end{equation}
    We eventually want to apply Lemma \ref{lem:large_connected_subset} to $G_x$.
    In order to do this, it suffices to show that there is an $\eta>0$ such that for all $\abs{\Omega^c}+\abs{\partial\Omega}$ small enough, we have $\abs{G_x}\geq\eta$.

    We may find such an $\eta$ by estimating the number of radial geodesics which intersect $\partial\Omega$.
    Let us define $\mathrm{proj}_{\sigma_x}:B(x,24r)\setminus B(x,\sigma_x)\rightarrow \partial B(x,\sigma_x)$ to be the map projecting radial geodesics onto $\partial B(x,\sigma_x)$.
    Using $\exp_x$ to compare with the Euclidean case, we see that this map has a Lipschitz constant of at most $L^2$.
    Now, consider $H_x=\mathrm{proj}_{\sigma_x}\Bigl(\bigl(B(x,24r)\setminus B(x,r)\cap B(a,20r)\cap\partial\Omega\bigr)\Bigr)$.
    We see that $\abs{H_x}\leq L^4\abs{B(a,20r)\cap \partial\Omega}$.
    For $\omega\in\partial B(x,\sigma_x)$ we shall abuse notation slightly by identifying $\omega$ with the radial geodesic emanating from $x$, say $\gamma_{\omega}$, such that $\gamma_{\omega}(\sigma_x)=\omega$. We will then write $\omega(t)$ to denote $\gamma_{\omega}(t)$.
    
    With this convention in mind, we observe that for any $t\in\left[\sigma_x,16r\right]$ and $\omega\in G_x^c\cap H_x^c$ we have that $\omega(t)$ lies in $\Omega^c$.
    Therefore, using that $\exp_x$ is $L$ bi-Lipschitz, we have the following string of estimates:
        \begin{align}
            \abs{\Omega^c}\geq\abs{\bigl\{\omega(t):\omega\in G_x^c\cap H_x^c;t\in\left[\sigma_x,16r\right]\bigr\}}&\geq\frac{(16r-\sigma_x)}{L^5}\abs{G_x^c\cap H_x^c}
            \\
            &\geq\frac{14r}{L^5}\left(\abs{G_x^c}-\abs{H_x}\right)
            \\
            &\geq\frac{14r}{L^5}\left(\abs{G_x^c}-L^4\abs{B(a,20r)\cap\partial\Omega}\right).
        \end{align}
    Rearranging terms and estimating shows us that
        \begin{equation}
            \abs{\Omega^c}+\frac{14r}{L}\abs{\partial\Omega}\geq\frac{14r}{L^5}\abs{G_x^c}.
        \end{equation}
    Since $\abs{\partial B(x,\sigma_x)}\geq L^{-2}\abs{\partial B(0,\sigma_x)}\geq\tfrac{r^2}{L^2}|\mathbb{S}^2|$, we see that
        \begin{equation}
            \abs{G_x}\geq\frac{r^2}{L^2}|\mathbb{S}^2|-\left(\frac{L^5}{14r}\abs{\Omega^c}+L^4\abs{\partial\Omega}\right)
        \end{equation}
    Letting $\eta=\tfrac{r^2}{2L^2}|\mathbb{S}^2|$, we see that for all $\abs{\Omega^c}+\abs{\partial\Omega}$ small enough,
    depending only on $r$ and $L$, we have that $\abs{G_x}\geq\eta$.
    Therefore, we may apply Lemma \ref{lem:large_connected_subset} to $G_x$ to find a connected subset $\widetilde{G}_x$ with smooth boundary such that
        \begin{equation}
            \abs{\partial\widetilde{G}_x}\leq\abs{\partial G_x}\leq\frac{1}{r}\abs{B(a,20r)\cap\partial\Omega}
        \end{equation}
    and
        \begin{equation}\label{eq:area_tilde-G-compl}
            \abs{G_{x}^c}\leq\frac{\abs{\partial\widetilde{G}_x}}{\ison_1\bigl(\partial B(x,\sigma_x)\bigr)}\leq\frac{2L^5}{\ison_1(\mathbb{S}^2)}\abs{B(a,20r\cap\partial\Omega}.
        \end{equation}
    Finally, we observe that since $\mathbb{S}^2$ is simply connected, and $\widetilde{G}_x$ is 
    connected, we may apply Lemma \ref{lem:two_dim_curve_mod} to $A=\mathrm{cl}\widetilde{G}_x$ 
    to obtain the following: for any $\omega,\zeta\in\mathrm{cl}\widetilde{G}_x$, there exists a 
    curve $\gamma\subset\mathrm{cl}\widetilde{G}_x$ connecting them such that
        \begin{align}
            L(\gamma)\leq d(\omega,\zeta)+\abs{\partial\widetilde{G}_x}&\leq L\sigma_x\pi
	    +\frac{1}{r}\abs{B(a,20,r)\cap\partial\Omega}
            \\
            &\leq2rL\pi+\frac{1}{r}\abs{B(a,20r)\cap\partial\Omega}.
        \end{align}
    
    Let us denote by $D_x$ the set given below:
        \begin{equation}
            D_x=\Bigl(\bigl\{\omega(t):\omega\in \widetilde{G}_x\cap H_x^c;t\in[\sigma_x,24r]\bigr\}\cap B(a,20r)\Bigr)\cup\mathrm{cl}\widetilde{G}_x.
        \end{equation}
    From looking at the definitions, we may see that any two points in $D_x$ may be connected together by a curve in $D_x$ which has length no greater than
    $48r+4rL\pi+r^{-1}\abs{B(a,20r)\cap\partial\Omega}$.

    Let us use the exact same construction in a neighborhood of $y$ to construct the set $D_y$.
    Ultimately, we wish to set $D_a=D_x\cup D_y$, and conclude that for any two points in $D_a$, there is a curve lying in $D_a$ which connects them, and which  has length no greater than
        \begin{equation}\label{eq:local_path_connected_estimate}
            96r+8rL\pi+2r^{-1}\abs{B(a,20r)\cap\partial\Omega}.
        \end{equation}
    In order to do this, it suffices to show that $D_x\cap D_y\neq\emptyset$.
    To do this, we will show that $\abs{D_x\cap D_y}>0$.
    Let us observe that 
        \begin{equation}
            B(a,r)\subset\{\omega(t):\omega\in\partial B(x,\sigma_x);t\in[\sigma_x,5r]\}
        \end{equation}
and 
    \begin{equation}
        B(a,r)\subset\{\nu(t):\nu\in\partial B(y,\sigma_y);t\in[\sigma_y,5r]\}
    \end{equation}
Therefore, we see that
    \begin{align}
        \abs{D_x\cap D_y}\geq&\abs{B(a,r)}
        \\
        &-125L^5\left(\abs{\left(\widetilde{G}_x\cap H^c_{x}\right)^c}+\abs{\left(\widetilde{G}_y\cap H^c_{y}\right)^c}\right)
    \end{align}
So, we see that
    \begin{equation}
        \abs{D_x\cap D_y}\geq\frac{\abs{B(0,r)}}{L^3}-250L^5\left(\frac{2L^5}{\ison_1(\mathbb{S}^2)}+L^4\right)\abs{B(a,20r)\cap\partial\Omega}.
    \end{equation}
Therefore, for all $\abs{\Omega^c}+\abs{\partial\Omega}$ small enough, we see that $D_x\cap D_y\neq\emptyset$, and so $D_a=D_x\cup D_y$ has the desired path connected property given in \eqref{eq:local_path_connected_estimate}.

Now, we need to estimate $\abs{B(a,20r)\cap\Omega\setminus D_a}$.
To do this, observe that $B(a,20r)$ is contained in
    \begin{equation}
        \{\omega(t):\omega\in\partial B(x,\sigma_x);t\in[\sigma_x,24r]\}\cup
        \{\nu(t):\omega\in\partial B(y,\sigma_x);t\in[\sigma_y,24r]\}.
    \end{equation}
Therefore, we see that
    \begin{align}
        \abs{B(a,20r)\setminus D_a}\leq&\abs{\{\omega(t):\omega\in(\widetilde{G}_x\cap H_x^c)^c;t\in[\sigma_x,24r]\}}
        \\
        &+\abs{\{\nu(t):\nu\in(\widetilde{G}_y\cap H_y^c)^c;t\in[\sigma_y,24r]\}}
    \end{align}
Using the estimates above, we see that
    \begin{equation}
        \abs{B(a,20r)\setminus D_a}\leq2(24^3L^5)\left(\frac{2L^5}{\ison_1(\mathbb{S}^2)}+L^4\right)\abs{B(a,20r)\cap\partial\Omega}.
    \end{equation}
\end{proof}
\begin{lemma}[\cite{Dong-Song-2023}]\label{lem:Dong-Song_Hausdorff_Close}
        Let $(M^{3},g)$ be a smooth closed three dimensional Riemannian manifold.
        For every $\varepsilon>0$ there exists a $\delta$ such that if $\Omega\subset M^3$ has smooth boundary, and we have that
            \begin{equation}
                \abs{\Omega}\geq\abs{M}-\delta
            \end{equation}
            \begin{equation}
                \abs{\partial\Omega}\leq\delta,
            \end{equation}
        then there exists a connected subset $\Omega'\subset\Omega$ with smooth boundary such that
            \begin{equation}
                \abs{\Omega'}\geq\abs{M^3}-\varepsilon
            \end{equation}
            \begin{equation}
                \abs{\partial\Omega'}\leq\varepsilon
            \end{equation}
        and
            \begin{equation}
                d_{\mathrm{GH}}\Bigl(\bigl(\Omega',\hat{d}^h_{\Omega'}\bigr),\bigl(\Omega',d^{h}\bigr)\Bigr)\leq\varepsilon.
            \end{equation}
    \end{lemma}
    \begin{proof}
        Fix $L=(1+\tfrac12\varepsilon)$, and fix $r>0$ to be chosen later.
        We will however assume that $24r\leq r_L$, where $r_L$ is the constant appearing in Lemma \ref{lem:locally_connected_subsets_of_Omega_in_balls}.
        We may additionally assume that $\abs{\Omega^c}+\abs{\partial\Omega}$ is small enough so that Lemma \ref{lem:locally_connected_subsets_of_Omega_in_balls} applies to balls of radius $20r$.
        Since $M$ is compact, we may cover $M$ in finitely many balls of the form $B(a_i,r)$;
            \begin{equation}
                1\leq\sum_{i=1}^{P}\chi_{B(a_i,r)}\leq P.
             \end{equation}
        For each $i$ let $D_i\subset\Omega\cap B(a_i,20r)$ be the set given in Lemma \ref{lem:locally_connected_subsets_of_Omega_in_balls}.
        Let $D=\bigcup_{i=1}^{P}D_{i}$, and consider $x,y\in D$ such that $d(x,y)\leq r$.
        Since $x,y\in D$, there exists $i(x)$ and $i(y)$ such that $x\in D_{i(x)}$ and
        $y\in D_{i(y)}$.
        Furthermore, since $\{B(a_{i},r)\}$ covers $M$, there exists $m(x)$ and $m(y)$
        such that $x\in B\left(a_{m(x)}\right)$ and $y\in B\left(a_{m(y)},r\right)$.
        Finally there are $z_x\in B\left(a_{i(x)},20r\right)$ and $z_y\in B\left(a_{i(y)},r\right)$ such that
        $d(x,z_x)=d(y,z_y)=r$, $B\left(z_x,r\right)\subset B\left(a_{m(x)},20r\right)\cap B\left(a_{i(x)},20r\right)$,
        and $B\left(z_y,r\right)\subset B\left(a_{m(y)},20r\right)\cap B\left(a_{i(y)},20r\right)$.
        Since $d(x,y)\leq r$, we also have that $B(x,r)\subset B\left(a_{m(x)},20r\right)\cap B\left(a_{m(y)},20r\right)$.
        As before, we can show that for $\abs{\Omega^c}+\abs{\partial\Omega}$ small enough, we have that
            \begin{align}
                \min\Bigl\{\abs{D_{i(x)}\cap D_{m(x)}},\abs{D_{i(y)}\cap D_{m(y)}},\abs{D_{m(x)}\cap D_{m(y)}}\Bigr\}>0,                
            \end{align}
        and in particular the intersections are not empty.
        As such, we see from Lemma \ref{lem:locally_connected_subsets_of_Omega_in_balls} that there is a path
        $\gamma\subset D$ connecting $x$ to $y$ such that
            \begin{equation}
                \mathrm{L}(\gamma)\leq3r\bigl(96+8L\pi+2r^{-2}\abs{\partial\Omega}\bigr).
            \end{equation}

        Lemma \ref{lem:locally_connected_subsets_of_Omega_in_balls} also gives us the following volume estimate:
            \begin{align}
                \abs{M\setminus D}&\leq\abs{\Omega^c}+2(24^3L^5)\sum_{i=1}^{P}\left(\frac{2L^5}{\ison_1(\mathbb{S}^2)}+L^4\right)\abs{B(a_{i},20r)\cap\partial\Omega}
                \\
                &\leq\abs{\Omega^c}+2P(24^3L^5)\left(\frac{2L^5}{\ison_1(\mathbb{S}^2)}+L^4\right)\abs{\partial\Omega}.
            \end{align}
        
        At this point, $D$ satisfies the requirement of local connectivity, but it does not necessarily have a smooth boundary, nor need this boundary be small.
        Therefore, we must modify $D$.
        Following Dong-Song \cite{Dong-Song-2023}, let $f:\Omega\rightarrow\R$ be given by
            \begin{equation}
                f(x)=\hat{d}_{\Omega}(z,D).
            \end{equation}
        Although the Lipschitz constant of $f$ depends on $\Omega$ and can be rather large, since $\Omega$ is a domain with smooth boundary, the local Lipschitz constant of $f$ is always bounded above by $2$, independently of $\Omega$.
        As such, for any $\eta>0$, we may find a smooth map $\phi:\Omega\rightarrow\R$ such that $\abs{\nabla\phi}\leq 2$ and $\|f-\phi\|_{L^{\infty}}<\tfrac{r}{64}$. Using the coarea formula, we have that
            \begin{equation}
                \int_{\tfrac{r}{32}}^{\frac{r}{16}}\abs{\phi^{-1}\{t\}\cap\Omega}dt=\int_{\Omega\cap\phi^{-1}[\tfrac{r}{32},\tfrac{r}{16}]}\abs{\nabla \phi}dV_{g}\leq2\abs{\Omega\setminus D}.
            \end{equation}
        Therefore, using Sard's lemma and a mean-value inequality, we may find a $t_0\in[\tfrac{r}{32},\tfrac{r}{16}]$ such that $\phi^{-1}\{t_0\}\subset\Omega$ is a submanifold of $M$, possibly with corners.
        Furthermore, we have that
            \begin{equation}
                \abs{\phi^{-1}\{t_0\}}\leq\frac{32}{r}\abs{\Omega\setminus D}.
            \end{equation}
        Observe that $\phi^{-1}\{t_0\}$ bounds the region $\phi^{-1}[0,t_0]$.
        Let us smooth out $\phi^{-1}\{t_0\}$ in such a way that the result $\partial\Omega'$ bounds a region $\Omega'\subset\phi^{-1}[0,t_0]\subset\Omega$ and
            \begin{equation}
                \abs{\partial\Omega'}\leq2\abs{\phi^{-1}\{t_0\}}\leq\frac{64}{r}\abs{\Omega\setminus D}.
            \end{equation}
        Next, observe that $D\subset\Omega'$, since $\|f-\phi\|_{L^{\infty}}<\tfrac{r}{64}$.
        Therefore, we have that
            \begin{equation}
                \abs{{\Omega'}^c}\leq\abs{D^c}.
            \end{equation}
            
        $\Omega'$ is very nearly the set we desire, however there is one more crucial property that it may not satisfy: we want to find a set such that if $x,y$ are elements of the set, and $d(x,y)\leq \tfrac{r}{2}$, then there is a path from $x$ to $y$ in $\Omega'$ which isn't too long.
        Of course, this is true if $x$ and $y$ are in $D$, but for more general $x$ and $y$, there is more work to be done.
        Since $\Omega$ and $\Omega'$ are domains with smooth boundary, the topology generated by $d(\cdot,\cdot)$, the topology generated by $\hat{d}_{\Omega}(\cdot,\cdot)$, and the topology generated by $\hat{d}_{\Omega'}(\cdot,\cdot)$ all agree.
        In particular, we have that $\Omega'$ is compact, and furthermore $\mathrm{cl}\bigl(\Omega'\setminus D\bigr)$ is compact.
        Therefore, we may find a finite collection of points $z_i$ in $\mathrm{cl}\bigl(\Omega'\setminus D\bigr)$ which form an $\tfrac{r}{8}$ net with respect to $\hat{d}_{\Omega'}$.
        In particular, by the definition of $\Omega'$, for each $i$ we may find a curve $c_{i}$ lying in $\Omega$ connecting $z_i$ to some $\hat{z}_i\in D$ which has length less than or equal to $\tfrac{r}{8}$.
        Consider a small cylindrical neighborhood $C_{i}$  of $c_{i}$. We may choose the area of its boundary to be arbitrarily small, and for every point to be within $\tfrac{r}{8}$ of the center curve, by choosing an arbitrarily small radius.
        Furthermore, we may perturb its boundary so that it intersects $\partial\Omega$ and $\partial\Omega'$ transversely.
        We may smooth out
        \begin{equation}
            \Omega'\bigcup_{i} \bigl(C_{i}\cap\Omega\bigr)
        \end{equation}
        to produce a region $\widetilde{\Omega}$ with the following properties.
        First, we may assume that $\abs{\partial\widetilde{\Omega}}\leq2\abs{\partial\Omega'}$.
        Next, we have that $\abs{\widetilde{\Omega}^c}\leq\abs{(\Omega')^c}$.
        Finally, suppose that $x$ and $y$ are in $\widetilde{\Omega}$ such that $d(x,y)\leq r$.
        Let us begin by assuming that $x,y\in\mathrm{cl}\bigl(\widetilde{\Omega}\setminus D\bigr)$. Then there are $z_{x}$ and $z_{y}$ in the $\tfrac{r}{8}$ net such that $\hat{d}_{\Omega'}(x,z_x)$ and $\hat{d}_{\Omega'}(y,z_y)$ are less than $\tfrac{r}{8}$.
        Furthermore, by construction of $\Omega'$, there are $\hat{z}_x$ and $\hat{z}_y$ in $D$ such that $\hat{d}_{\tilde{\Omega}}(z_x,\hat{z}_x)$ and $\hat{d}_{\tilde{\Omega}}(z_y,\hat{z}_y)$ are less than $\tfrac{r}{8}$.
        Thus, we see that
            \begin{align}
                d(\hat{z}_x,\hat{z}_y)&\leq d_{\tilde{\Omega}}(\hat{z}_x,z_x)+d_{\tilde{\Omega}}(z_x,x)+d(x,y)+d_{\tilde{\Omega}'}(y,z_y)+d_{\tilde{\Omega}}(z_y,\hat{z}_y)
                \\
                &\leq r.
            \end{align}
        Let us recall that $L$ denotes the bi-Lipschitz constant of $\exp$ for balls with radius
        less than or equal to $r_L$, which are also geodesically convex by assumption.
        Since $d(\hat{z}_x,\hat{z}_y)\leq r$, it follows from the construction of $D$,
        see Lemma \ref{lem:locally_connected_subsets_of_Omega_in_balls}, that we may find 
        a curve $\gamma$ in $D$ connecting $\hat{z}_x$ and $\hat{z}_y$ which has length bounded as follows:
            \begin{equation}
                L(\gamma)\leq r\bigl(96+8L\pi+2r^{-2}\abs{B(a,20r)\cap\partial\Omega}\bigr).
            \end{equation}
        The other cases are similar, but use fewer applications of the triangle inequality.
        
        To save space in the calculations below, we set $K=K(\abs{\partial\Omega},L,r)$, where
            \begin{equation}
                K(\abs{\partial\Omega},L,r)=96+8L\pi+2r^{-2}\abs{B(a,20r)\cap\partial\Omega},
            \end{equation}
        Let $x,y$ be any two points in $\widetilde{\Omega}$, and let $\gamma$ be a length minimizing geodesic connecting them.
        We may split $\gamma$ into at most $\left\lceil\tfrac{16\diam(M,g)}{r_{L}}\right\rceil$ segments of length $l_i$,
        where $\tfrac{1}{16}r_{L}\leq l_i\leq\tfrac{1}{4}r_{L}$.
        Let us consider any two successive sections of $\gamma$, say $\gamma(x_0,y_0)$ and $\gamma_{i_0+1}(x_{1},y_{1})$, with endpoints $x_j,y_j$ for $j=0,1$, respectively.
        Let $c_0$ be the midpoint of $\gamma(x_0,y_0)$, and consider $\widetilde{\gamma}_0=\exp^{-1}_{c_0}\bigl(\gamma(x_0,y_0)\bigr)\subset T_{c_0}M$.
        Let $\widetilde{C}_r$ be the corresponding cylinder about $\widetilde{\gamma}_0$ of radius $\tfrac{r}{8L}$; it is foliated by curves, as seen below:
            \begin{equation}
                \widetilde{C}_{r}=\left\{\widetilde{\gamma}_{\omega}:\widetilde{\gamma}_{\omega}(t)=t+\omega;\omega\perp\frac{d}{dt}\widetilde{\gamma}_1(t);\abs{\omega}\leq\tfrac{r}{8L}\right\}.
            \end{equation}
        Consider $\exp_{c_0}^{-1}(\partial\widetilde{\Omega})$, and let $H=\{\omega:\widetilde{\gamma}_{\omega}\cap\exp^{-1}_{c_0}(\partial\widetilde{\Omega})\neq\emptyset\}$.
        By projecting the curve $\widetilde{\gamma}_{\omega}$ onto $\omega$, we observe that 
            \begin{align}
                \abs{H}&\leq\abs{\exp^{-1}_{c_0}(\partial\widetilde{\Omega})}
                \\
                &\leq L^2\abs{\partial\widetilde{\Omega}}
                \\
                &\leq2L^2\abs{\partial\Omega'}
                \\
                &\leq\frac{128}{r}\abs{\Omega\setminus D}.
            \end{align}

        Now observe that for $\omega\in H^c$ we have either that $\exp_{c_0}(\widetilde{\gamma}_{\omega})\subset\widetilde{\Omega}$ or $\exp_{c_0}(\widetilde{\gamma}_{\omega})\subset\widetilde{\Omega}^c$.
        Suppose that for all $\omega\in H^c$ that $\exp_{c_0}(\widetilde{\gamma}_{\omega})\subset\widetilde{\Omega}^c$.
        Then, it follows that
            \begin{equation}
                \abs{\widetilde{\Omega}^c}\geq\frac{r_{L}}{16L^5}\left(\pi\left(\frac{r}{8L}\right)^2-\abs{H}\right).
            \end{equation}
        Rearranging terms shows us that
            \begin{equation}
                \abs{\widetilde{\Omega}^c}+\frac{r_{L}\abs{H}}{15L^5}\geq\pi\left(\frac{r}{8L}\right)^2.
            \end{equation}
        This gives a contradiction for $\abs{\Omega^c}+\abs{\partial\Omega}$ small enough depending only on $r$
        and $L$.
        Therefore, thre is at least one $\omega\in H^c$ such that $\exp_{c_0}(\widetilde{\gamma}_{\omega})\subset\widetilde{\Omega}$.
        Let $\omega_0$ denote such an element, and let $\gamma_{\omega_0}(\widetilde{x}_{0},\widetilde{y}_{0})$ denote $\exp_{c_0}(\widetilde{\gamma}_{\omega})$, which has endpoints $\widetilde{x}_{0}$ and $\widetilde{y}_0$
        We may similarly find $\gamma_{1}(\widetilde{x}_1,\widetilde{y}_1)$ with end points $\widetilde{x}_1$ and $\widetilde{y}_1$.
        Then, we have that
            \begin{equation}
                \mathrm{L}\bigl(\gamma_{\omega_j}(\widetilde{x}_{j},\widetilde{y}_{j})\bigr)\leq L d(x_j,y_j),
            \end{equation}
        for $j=0,1$.
        Furthermore, we see that
            \begin{equation}
                d\bigl(\widetilde{y}_0,\widetilde{x}_1\bigr)\leq\frac{r}{4}.
            \end{equation}
        As such, there is a curve connecting them, which lies entirely in $\widetilde{\Omega}$ and which has length
        less than or equal to $Kr$.
        Let us modify $\gamma$ along the section $\gamma_0$ and $\gamma_1$ by replacing $\gamma_j$ with $\widetilde{\gamma}_{\omega_j}$, and connecting their endpoints by the curve above.
        Doing this for all of the segments yields a curve $\tilde{\gamma}$ which lies entirely in $\widetilde{\Omega}$
        and has length bounded above as follows:
            \begin{equation}
                \mathrm{L}\bigl(\widetilde{\gamma}\bigr)\leq\left(\left\lceil\tfrac{16\diam(M,g)}{r_{L}}\right\rceil+2\right)Kr+Ld(x,y),
            \end{equation}
        where the $+2$ comes from the fact that we may have to peturb the endpoints $x,y$ of $\gamma$ in a similar way as above.
        Therefore, for $r$ chosen small enough, and $\abs{\Omega^c}+\abs{\partial\Omega}$ small enough depending only on $r$ and $L$, the above has the following estimate:
            \begin{equation}
                \mathrm{L}\bigl(\widetilde{\gamma}\bigr)\leq(1+\varepsilon)d(x,y).
            \end{equation}
\end{proof}
\subsection{Fundamental Domains}
In this section we give a construction of fundamental domains, which can be found in 
\cite[Section IV.3]{chavel2006riemannian}.
\begin{definition}[\cite{chavel2006riemannian}]
	Let $(M,g)$ be a closed Riemannian manifold, and let $p$ be any point in $M$,
	and denote by $\mathcal{S}_{p}$ the collection of unit vectors in $\mathrm{T}_{p}M$,
	and let $\mathcal{S}M$ denote the unit tangent bundle.
	Finally, let $\pi:\mathrm{T}M:\rightarrow{}M$ be the projection map.
	Let us define a map $c:\mathcal{S}M\rightarrow{}(0,\infty{})$ as follows
	\begin{equation}
		c(\xi)=\sup\left\{t:d\bigl(p,\exp_{\pi(\xi)}(t\xi)\bigr)=t\right\}
	\end{equation}
\end{definition} 
We have the following important result
\begin{theorem}[\cite{chavel2006riemannian}]\label{thm:continuity_of_cut_distance}
	Let $(M,g)$ be a complete Riemannian manifold without boundary.
	Then, the map $c:\mathcal{S}M\rightarrow{}(0,\infty]$ is continuous.
\end{theorem}
\begin{proof}
	Fix $\xi$ in $\mathcal{S}M$ and let $\xi_{k}$ be a sequence approaching $\xi$.
	For notational convenience, let $p=\pi(\xi)$, $p_{k}=\pi(\xi_{k})$, and
	$d_{k}=c(\xi_{k}$.
	
	We will first show that $\limsup_{k\rightarrow{}\infty{}}c(\xi_{k})\leq c(\xi)$.
	Pick a subsequence $d_{k(i)}$ such that 
	$\lim_{i\rightarrow{}\infty{}}d_{k(i)}=\limsup_{k\rightarrow{}\infty{}}d_{k}$.
	If $\limsup_{k\rightarrow{}\infty{}}d_{k}=\infty{}$, then, 
	for any $T>0$ one has that $d_{k(i)}>T$ for all $i$ sufficiently large.
	Since $\xi_{k}\rightarrow{}\xi$, it follows that $\exp_{p_{k(i)}}(T\xi_{k(i)\eta})$ converges
	to $\exp_{p}(T\xi)$.
	Since distance is a continuous function, we see that
	\begin{equation}
		d(p,\exp_{p}(T\xi))=
		\lim_{i\rightarrow{}\infty{}}d(p_{k(i)},T\xi_{k(i)})=T.
	\end{equation}
	Therefore, by the definition of $c$ we see that $c(\xi)\geq T$.
	As $T$ can be chosen to be arbitrarily large, we see that $c(\xi)=\infty$
	
	Next, suppose that $\limsup_{k\rightarrow{}\infty{}}d_{k}=\delta<\infty{}$.
	Then, for any $\varepsilon{}$ we may find an $N$ so that for all $i\geq N$
	we have that $\delta-\varepsilon{}<d_{k(i)}$.
	It follows that
	\begin{equation}
		d\Bigl(p,\exp_{p}\bigl((\delta-\varepsilon{})\xi\bigr)\Bigr)=
		\lim_{i\rightarrow{}\infty{}}d\Bigl(p_{k(i)},
		\exp_{p_{k(i)}}\left(\left(\delta-\varepsilon{}\right)\xi_{k(i)}\right)\Bigr)
                =\delta-\varepsilon{}.
	\end{equation}
        So, we see that $c(\xi)\geq\delta-\varepsilon{}$, where $\varepsilon{}$ was arbitrary.

        Continuity will now follow if we can show that 
        $\liminf_{k\rightarrow{}\infty{}}d_{k}\geq c(\xi)$.
        By way of contradiction, suppose that
        $\liminf_{k\rightarrow{}\infty{}}d_{k}+2\varepsilon{}< c(\xi)$.
        Let $U_{p}$ denote the open star-shaped region about $p=\pi(\xi)$ on which
        $\exp_{p}$ is a diffeomorphism, let
        $T=\liminf_{k\rightarrow{}\infty{}}d_{k}$, and let $\xi_{k(i)}$ be a subsequence
        such that $\lim_{i\rightarrow{}\infty{}}d_{k(i)}=T$.
        Since $\xi_{k(i)}\rightarrow{}\xi$, we have that 
        $\lim_{k\rightarrow{}\infty{}}\exp_{p_{k(i)}}\bigl((T+\varepsilon{})
        \xi_{k(i)}\bigr)=\exp_{p}\bigl((T+\varepsilon{})\xi\bigr)$.
        Since $\exp_{p}(T\xi)\in{}U_{p}$, it follows that for all $i$ large enough, we have that
        $\exp_{p_{k(i)}}\bigl((T+\varepsilon{})\xi_{k(i)}\bigr)\in{}U_{p}$.
        However, we also have for all $i$ sufficiently large that
        $c(\xi_{k(i)})=d_{k(i)}<T+\varepsilon{}$.
        This is a contradiction, and so we see that
        $\liminf_{k\rightarrow{}\infty{}}d_{k}\geq c(\xi)$.
        This establishes the continuity of $c$ on $\mathcal{S}M$.
\end{proof}
\begin{remark}
    Here we avoided the second half of proof presented in \cite{chavel2006riemannian}, since by the results in
    \cite{margerin1993general} conjugate loci are in general ill behaved.
\end{remark}
\begin{definition}
    Given a complete Riemannian manifold $(M,g)$ and a point $p$ in $M$,
    we will let $D_{p}$ denote the following set:
    \begin{equation}
        D_{p}=\{c(\xi)\xi:\xi\in{}\mathcal{S}_{p}M\},
    \end{equation}
    the graph of the continuous function $c:\mathcal{S}_{p}\rightarrow{}(0,\infty{}]$.
    Furthermore, let 
    \begin{equation}
        G_{p}=\{t\xi:t<c(\xi);\xi\in{}\mathcal{S}_{p}\}.
    \end{equation}
\end{definition}
\begin{proposition}
    Let $(M,g)$ be a closed Riemannian manifold, and let $p$ in $M$.
    Then, $D_{p}$ has measure zero.
\end{proposition}
\begin{proof}
    Observe that $\mathcal{S}_{p}$ and $M$ are compact.
    As such, it follows that $c(\mathcal{S}_{p})$ is contained in a compact interval bounded
    away from zero.
    Since $D_{p}$ is the graph of $c$ over $\mathcal{S}_{p}$, it follows from using 
    polar coordinates and the Fubini-Tonelli theorem that
    \begin{equation}
    \int_{T_{p}M}\chi_{D_{p}}=\int_{\mathcal{S}_{p}}
            \int_{a}^{b}n\omega_{n}r^{n-1}\chi_{D_{p}}(r,\theta)drd\theta=
            0.
    \end{equation}
\end{proof}
\begin{lemma}\label{lem:fund_domain_exists}
    Let $(M,g)$ be a closed Riemannian manifold, and let $(\widetilde{M},\pi^{*}g)$ be its
    universal cover.
    For any given $p\in{}M$ and $\widetilde{p}\in{}\widetilde{M}$,
    let $\mathrm{Dir}(\widetilde{p})$ denote the following set
    \begin{equation}
        \mathrm{Dir}(\widetilde{p})=\exp_{\widetilde{p}}\circ(d\pi_{\widetilde{p}})^{-1}
        G_{p}.
    \end{equation}
    Then, we have that $\mathrm{Dir}(\widetilde{p})$ is a fundamental domain.
\end{lemma}
\begin{proof}
    We claim that $\mathrm{Dir}(\widetilde{p})$ is a fundamental domain.
    To begin, suppose that $x,y\in{}\mathrm{int}\mathrm{Dir}(\widetilde{p})$
    are such that $\pi(x)=\pi(y)$.
    This would mean that there are at least two length minimizing geodesics from
    $p$ to $\pi(x)=\pi(y)$, but this contradicts the definition of $G_{p}$.
    Observe that this implies that
    $\overline{\mathrm{Dir}}(\widetilde{p})=\exp_{\widetilde{p}}\circ(d\pi_{\widetilde{p}})^{-1}
    \overline{G}_{p}$, and as a consequence we also have that
    $\partial\mathrm{Dir}(\widetilde{p})\subset{}\exp_{\widetilde{p}}
    \circ(d\pi_{\widetilde{p}})^{-1}D_{p}$.
    Therefore, since $D_{p}$ has measure zero and 
    $\exp_{\widetilde{p}}\circ(d\pi_{\widetilde{p}})^{-1}$ is a smooth map,
    it follows that $\partial{}\mathrm{Dir}(\widetilde{p})$ has measure zero.
    Finally, since $\overline{G}_{p}$ is path connected, it follows that
    $\overline{\mathrm{Dir}}(\widetilde{p})$ is as well.
\end{proof}

	\bibliography{bibliography}
	\bibliographystyle{amsalpha}
\end{document}